\documentclass[11pt]{amsart}

\usepackage[final]{nick}

\thanks{The first and second authors were partially supported by 
EPSRC Grant GR/S28655/01 and JSPS Grant Kiban-B 15340023, 
respectively.}

\newcommand{\nonnegative}{nonnegative }

\newcommand{\Depend}[1]{}

\newenvironment{nestedtheorem}{\em}{}

\newcounter{COUNTER}
\newcommand{\ITEM}{\arabic{COUNTER}. \addtocounter{COUNTER}{1}}
\newcommand{\STARTSTEPS}{\setcounter{COUNTER}{1}}
\newcommand{\STEP}{\textit{Step \ITEM}}
\newcommand{\CASE}{\textit{Case \ITEM}}

\newcommand{\spaceperiod}{\,\,.}
\newcommand{\spacecomma}{\,\,,}

\newcommand{\optionalA}{}
\newcommand{\optionalAn}{}

\newcommand{\supnorm}[2]{{{\left|\!\left|#1\right|\!\right|}_{#2}}}
\newcommand{\matnorm}[1]{{{\left|\!\left|#1\right|\!\right|}}}

\DeclareMathOperator{\singx}{sing}
\newcommand{\sing}[1]{\singx(#1)}

\newcommand{\galt}{\widehat{g}}
\newcommand{\gdel}{{g'}}
\newcommand{\LoopUpSL}[1]{\Lambda_{#1}^{\uparrow}\matSL{2}{\bbC}}
\newcommand{\LoopuMeroSL}[1]{\Lambda_{#1}^{\ast{\mathrm M}}\matSL{2}{\bbC}}

\newcommand{\preliminarySection}{the preliminary section}

\newcommand{\simplefactor}[2]{\mathfrak{g}[#1,\,#2]}
\newcommand{\simplefactorUN}[2]{\mathfrak{g}^\circ[#1,\,#2]}

\newcommand{\converge}[1]{\text{\shortstack{\scriptsize $#1$\\$\to$}}}

\newcommand{\delf}{\psi}
\newcommand{\dels}{\sigma}
\newcommand{\delR}{R}

\newcommand{\LIMIT}[4]%
  {{\lim_{#4}{\supnorm{#1-#2}{#3}}=0}}
\newcommand{\LIMITZ}[4]%
  {{\lim_{#4}{\supnorm{#1}{#3}}=0}}

\newcommand{\CONVERGE}[2]{$#1$ and $#2$ have the convergence}

\DeclareMathOperator{\cchange}{\vartheta}

\newcommand{\ELL}{\mathscr{L}}
\newcommand{\CPLUS}{C_+}
\newcommand{\CPLUSP}{V}
\newcommand{\RRR}{r}
\newcommand{\VVV}{V}

\newcommand{\gaugeConstant}{\kappa}

\newcommand{\delIA}{\calI_A}
\newcommand{\delJA}{\calJ_A}
\newcommand{\delKA}{\calK_A}
\newcommand{\delSA}{\calS_A}

\newcommand{\gfact}{h}

\newcommand{\theoremname}[1]{}

\newenvironment{xtext}{}{}

\title[Delaunay ends]
{Delaunay Ends of Constant Mean Curvature Surfaces}

\author{M. Kilian}
\address{Martin Kilian, Institut f\"ur Mathematik,
Universit\"at Mannheim, 68131 Mannheim, Germany.}
\email{kilian@rumms.uni-mannheim.de}

\author{W. Rossman}
\address{Wayne Rossman, Department of Mathematics,
Kobe University, Rokko Kobe 657-8501, Japan.}
\email{wayne@math.kobe-u.ac.jp}

\author{N. Schmitt}
\address{Nicholas Schmitt, Mathematisches Institut,
Universit\"at T\"ubingen, 72076 T\"ubingen, Germany.}
\email{nick@gang.umass.edu}

\begin{document}

\begin{abstract}
The generalized Weierstrass representation is used to
analyze the asymptotic behavior of a constant mean
curvature surface that arises locally from an
ordinary differential equation with a regular singularity.
We prove that a holomorphic perturbation of an ODE that represents a
Delaunay surface generates a constant mean curvature surface which
has a properly immersed end that is asymptotically Delaunay.
Furthermore, that end is embedded if the Delaunay surface is unduloidal.
\end{abstract}

\maketitle



\typeout{=================intro}
\section*{Introduction}
\begin{xtext}
Delaunay surfaces play a prominent role in the theory of
non-compact complete constant mean curvature (CMC) surfaces because they 
constitute the simplest possible end behavior.
%
%
A famous result by Korevaar, Kusner and
Solomon~\cite{Korevaar_Kusner_Solomon_1989},
building on results of Meeks~\cite{Meeks_1988},
asserts that a properly embedded annular end
of a CMC surface is a Delaunay end.
The study of Delaunay ends
by the conjugate surface methods of Grosse-Brauckmann, Kusner, and 
Sullivan \cite{Grosse-Brauckmann_Kusner_Sullivan_2003,Grosse-Brauckmann_Kusner_Sullivan_2005} require the additional assumption of Alexandrov embeddedness,
and are limited to embedded (unduloidal) Delaunay ends.
The gluing techniques of Mazzeo and Pacard \cite{Mazzeo_Pacard_2001}
are limited to attaching Delaunay ends with small asymptotic necksizes.
The methods used in this paper, based on the generalized Weierstrass
representation
of Dorfmeister, Pedit and Wu~\cite{Dorfmeister_Pedit_Wu_1998},
provide a means to study
both embedded unduloidal and non-Alexandrov-embedded nodoidal type ends
of arbitrary asymptotic necksize.

The generalized Weierstrass representation describes 
constant mean curvature immersions locally via
\emph{holomorphic potentials}. 
The relation between the potential and the immersion involves
a loop group valued differential equation, a loop group 
factorization and a Sym-Bobenko type formula 
\cite{Sym_1985,Bobenko_CMC_1991}. 
Hence the method provides only an indirect relation between the
geometric properties of the induced immersion and its potential.
Nevertheless, the method has been useful in proving the existence of many 
new classes of non-compact constant mean curvature surfaces 
with non-trivial topology. Particular progress has been 
made when the surface is homeomorphic to an $n$-punctured 
Riemann sphere \cite{Kilian_McIntosh_Schmitt_2000, Kilian_Schmitt_Sterling_2004, Schmitt_2006, Schmitt_Kilian_Kobayashi_Rossman_2006}. 
In this case, the punctures correspond to 
poles of the potential. 
Graphics of these surfaces \cite{Schmitt_2002} 
have long suggested that simple poles with appropriate residues 
yield the asymptotic end behavior
of a Delaunay surface. 
We prove this correlation between simple poles and Delaunay ends. 

More specifically, given a holomorphic potential $A\,dz/z$ of a 
Delaunay surface,
consider a holomorphic perturbation
$\xi = A\,dz/z + \Order(z^{0}) dz$.
Our main result,
stated precisely as \autoref{thm:main1}, and generalized
in \autoref{thm:main2},
is:
\begin{theorem*}
An annular constant mean curvature immersion induced by
a holomorphic perturbation of a Delaunay potential
is $C^\infty$-asymptotic to a
half-Delaunay surface.
In particular, it is properly immersed.
Moreover, if the half-Delaunay surface is embedded,
then the end of the immersion is properly embedded.
\end{theorem*}

The surface induced by a perturbed Delaunay potential
may gain topology or geometric complexity ---
see for example the
$n$-noids~\cite{Schmitt_Kilian_Kobayashi_Rossman_2006, Rossman_Schmitt_2006, Schmitt_2006}
and higher genus examples with
ends~\cite{Kilian_Kobayashi_Rossmann_Schmitt_2005}.
Nonetheless at $z=0$, the perturbed surface is
asymptotic to the underlying Delaunay surface.

%
The convergence of the surfaces is obtained
by showing that their moving frames and metrics converge.
More specifically,
let $\Phi_0$ and $\Phi$ be the respective solutions to
the ODEs $d\Phi_0=\Phi_0 A\,dz/z$
and its perturbation $d\Phi=\Phi\xi$.
The convergence of the ratio of $\Phi$ to $\Phi_0$
is shown using the holomorphic gauge relating them
at the regular singularity $z=0$.
%
The periodicity of the Delaunay surface provides
growth rate estimates on the positive part of $\Phi_0$
by Floquet analysis.
This leads to the convergence
of their unitary and positive factors,
in turn implying $C^1$-convergence of the surfaces.
A bootstrap argument strengthens this to $C^\infty$-convergence.

In Part 2 we deal with the situation in which the initial condition to the ODE 
does not extend holomorphically from the $r$-circle to the unit circle, even 
though the monodromy of the solution is still unitary. We show that in this 
setting the solution has acquired singularities that arise from 
Bianchi-B\"acklund transforms. Thus the second part accommodates the 
additional singularities that appear from dressing by simple factors, and 
proves that adding bubbles to a surface with a Delaunay end preserves this 
Delaunay end.


%


We have benefited greatly on this long-term project from discussions with
Franz Pedit, Ian McIntosh, Martin Schmidt and many others.
We thank J. Dorfmeister and
H. Wu~\cite{Dorfmeister_Wu_2002}
and S.-P. Kobayashi~\cite{Kobayashi_2006}
for making available their unpublished works.
\end{xtext}

\begin{figure}[ht]
\centering
\includegraphics[width=4in]{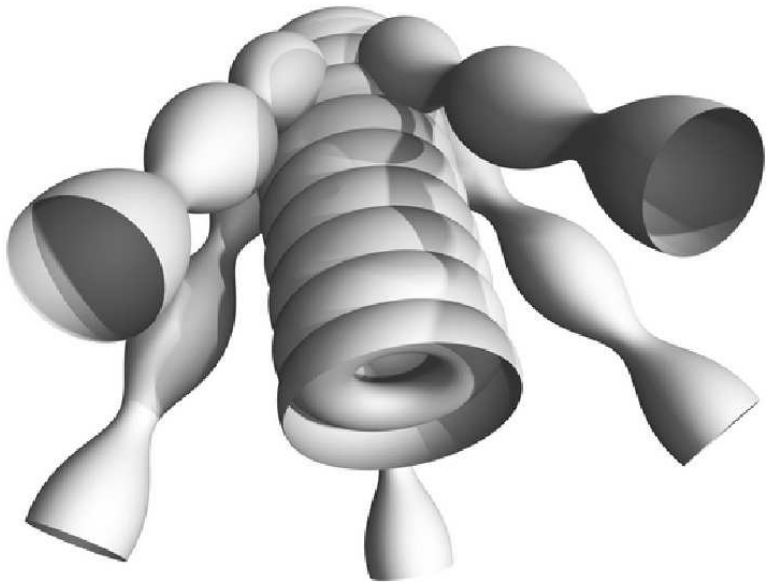}
\caption{
A CMC immersion of the six-punctured sphere with asymptotically Delaunay ends
and pyramidal symmetry~\cite{Schmitt_2006}.
Five of the six ends are unduloidal; the sixth is nodoidal
with large negative weight.
This image was created with 
CMCLab~\cite{Schmitt_2002b}, 
a freely available software implementation
of the generalized Weierstrass representation. 
}
\end{figure}



\section*{The generalized Weierstrass representation}
\typeout{=================preliminary}
\subsection*{\texorpdfstring{The $r$-Iwasawa factorization}%
  {The r-Iwasawa factorization}}
We will use the following subsets of $\bbC$:
\begin{align*}
\bbS^1 &= \{\lambda\in\bbC\suchthat \abs{\lambda} = 1\}
\AND
\bbC^\ast = \bbC\setminus\{0\}\spacecomma\\
\calC_r &= \{\lambda\in\bbC\suchthat \abs{\lambda\,}=r\}\spacecomma
  \quad r\in(0,\,1]\spacecomma\\
\calD_r &= \{\lambda\in\bbC\suchthat \abs{\lambda\,}<r\}\spacecomma
  \quad
  \calD_r^\ast = \calD_r\setminus\{0\}\spacecomma
  \quad r\in(0,\,1]\spacecomma\\
\calA_{s,r} &=
 \{\lambda\in\bbC\suchthat s < \abs{\lambda\, } < r\}\spacecomma
\quad 0<s<r\spacecomma\\
\calA_r &= \calA_{r,1/r}\spacecomma\quad r\in(0,\,1)\spaceperiod
\end{align*}
Given a domain $\calU\subset\bbC\setminus\{0\}$
which is invariant under the map
$\lambda\mapsto 1/\ol{\lambda}$,
and a holomorphic map $X:\calU\to\mattwo{\bbC}$,
define the holomorphic map $X^\ast:\calU\to\mattwo{\bbC}$ by
\begin{equation*}
X^*(\lambda) := 
\overline{X(1/\bar{\lambda})}^{\,t}.
\end{equation*}
We will use the following loop groups:
\begin{itemize}
\item
For $r\in(0,\,1]$,
$\LoopSL{r}$ is the group of analytic maps $\calC_r\to\matSL{2}{\bbC}$.
\item
For $r\in(0,\,1)$,
$\LoopuSL{r}\subset\LoopSL{r}$ is the subgroup of loops each of
which is the boundary of a holomorphic map $Y:\calA_r\to\matSL{2}{\bbC}$
satisfying the condition $Y^\ast = Y^{-1}$.
$\LoopuSL{1}\subset\LoopSL{1}$ is the subgroup of loops $X$
satisfying the condition $X^\ast = X^{-1}$.
\item
For $r\in(0,\,1]$,
$\LoopposSL{r}\subset\LoopSL{r}$ is the subgroup
of loops for which each loop is the boundary of a holomorphic
map $\calD_r\to\matSL{2}{\bbC}$.
\item
Let $\calT\subset\matSL{2}{\bbC}$ denote the group of upper triangular
matrices whose diagonal elements are in $\bbR_{>0}$.
For $r\in(0,\,1]$,
$\LooppSL{r}\subset\LoopposSL{r}$ is the subgroup
of loops $X$ such that $X(0)\in\calT$.
\item
For $r\in(0,\,1)$,
$\LoopUpSL{r}\subset\LoopSL{r}$
is the group of analytic maps $X\in\LoopSL{r}$
such that $X$ is the boundary of a holomorphic map
$\calA_{r,\,1}\to\matSL{2}{\bbC}$.
\item
For $r\in(0,\,1)$,
$\LoopuMeroSL{r}\subset\LoopSL{r}$
is the group of analytic maps $X\in\LoopSL{r}$
such that $X$ is the boundary of a meromorphic map
$\calA_{r}\to\matSL{2}{\bbC}$
satisfying $X^\ast = X^{-1}$ away from its poles.
For $r=1$, $\LoopuMeroSL{1} = \LoopuSL{1}$.
\end{itemize}

For $k\in\bbZ_{\ge 0}$,
we define the $C^k$-topology on each of these loop groups
with respect to the loop parameter $\lambda$ in their
respective domains $\calC_r$, $\calA_r$ or $\calD_r$.
The $C^\infty$-topology is the intersection of the $C^k$-topologies.
%
%
%
The  multiplication map 
\begin{equation*}
  \LoopuSL{r} \times \LooppSL{r} \to \LoopSL{r}
\end{equation*}
is a $C^\infty$
diffeomorphism~\cite{McIntosh_Nonlinearity_1994, Burstall_Pedit_1995}.
The unique factorization of a loop 
$\Phi \in \LoopSL{r}$ into 
\begin{equation*}
  \Phi = \Uni{r}{\Phi}\cdot\Pos{r}{\Phi}
\end{equation*}
is the \emph{$r$-Iwasawa factorization} of $\Phi$.
We call $\Uni{r}{\Phi}$ the \emph{r-unitary factor}
and $\Pos{r}{\Phi}$ the \emph{r-positive factor} of $\Phi$.
The \emph{QR-factorization}
is the Iwasawa factorization of constant loops
$\matSL{2}{\bbC}\to\matSU{2}{}\times\calT$.

Dressing a loop $X\in\LoopSL{r}$ by $C\in\LoopSL{r}$
is left-multiplication of $X$ by $C$,
followed by projection of the $r$-Iwasawa
factorization to the unitary group.
We denote the dressed loop by $\dress{r}{C}{X}$ or $\Uni{r}{CX}$.

\subsection*{The generalized Weierstrass representation}

\begin{xtext}
The generalized Weierstrass representation~\cite{Dorfmeister_Pedit_Wu_1998}
represents harmonic maps in terms of certain
holomorphic  $1$-forms with values in a loop algebra 
(\emph{holomorphic potentials}). 
\end{xtext}
This representation constructs all constant mean curvature (CMC) surfaces
in the 3-dimensional Euclidean, spherical and hyperbolic
spaceforms~\cite{Schmitt_Kilian_Kobayashi_Rossman_2006},
and is as follows for Euclidean 3-space:

1. Let $\Sigma$ be a Riemann surface. With $r\in(0,\,1]$, let $\xi$ be an
\emph{$r$-potential}, that is, a
$\Loopsl{r}$-valued
differential form on $\Sigma$
which is the boundary of a meromorphic differential
on $\calD_r$, with a pole only at $\lambda=0$, which is simple and is only
in the upper-right entry.
To avoid branch points in the induced surface, assume also that the coefficient
of $\lambda^{-1}$ in the series expansion of the upper-right entry of $\xi$
in $\lambda$ at $\lambda=0$ is never zero on $\Sigma$.

2. Let $\Phi$ be a solution to the ordinary differential equation
$d\Phi=\Phi\xi$ on the universal cover $\widetilde\Sigma$ of $\Sigma$.
We call $\Phi$ the \emph{holomorphic frame}.

3. Then $F=\Uni{r}{\Phi}$ is the extended frame of some CMC immersion
$f:\widetilde\Sigma\to\bbR^3$.

4. The Sym formula
\begin{equation*}
f=\Sym{r}{\Phi} = -2 H^{-1} F' F^{-1}
\end{equation*}
computes the associate family of CMC immersion $f$ with
constant mean curvature $H\in\bbR^\ast$ from the frame $F$,
with associate family parameter $\lambda\in\bbS^1$.
Here, the prime denotes differentiation with respect to $\theta$,
where $\lambda=e^{i\theta}\in\bbS^1$.

The immersion $f$ of $\widetilde\Sigma$ descends to an immersion of $\Sigma$
at $\lambda=1$
if and only if every element $M_F$ of the monodromy group of the
extended frame $F$
satisfies
\begin{equation*}
M_F(1)=\pm\id
\AND
M_F'(1)=0\spaceperiod
\end{equation*}

Given $d\Phi=\Phi\xi$ and an analytic \emph{gauge}
$g=g(z,\,\lambda):\Sigma\to\LoopSL{r}$,
then $\Psi:=\Phi g$ satisfies the equation $d\Psi = \Psi\eta$,
where
\begin{equation*}
\eta = \gauge{\xi}{g} := g^{-1}\xi g + g^{-1} dg\spaceperiod
\end{equation*}

\part{Delaunay asymptotics}
\label{part1}

\begin{xtext}
In Part~\ref{part1} we show that
under the assumption of unitary monodromy,
an immersion constructed from a perturbed Delaunay potential
is asymptotic to a Delaunay immersion.
In Part~\ref{part2} we generalize this result
to the setting of $\LoopSL{r}$ for arbitrary $r$.
\end{xtext}

\typeout{=================outline1}
\section*{Outline of results}
\label{sec:outline1}

\setcounter{COUNTER}{1}

\begin{xtext}


\autoref{Sec:delaunay} discusses the construction of
the family of Delaunay immersions
via the generalized Weierstrass representation.
The generalized Weierstrass potential
for the Delaunay immersion is of the form $A\,dz/z$, where
$A$ is an $\matsl{2}{\bbC}$-valued
\emph{Delaunay residue}
given in \autoref{thm:delres}.
We compute the Iwasawa factorization of the
holomorphic Delaunay frame $\exp(A\log z)$
in \autoref{thm:delframe}.
This is an extension of a result
in~\cite{Schmitt_Kilian_Kobayashi_Rossman_2006}.
These factors are used in \autoref{sec:delgrowth} to compute the
growth rate of the positive Iwasawa factor,
and in \autoref{Sec:simple-factor-dressing}
to compute asymptotics of dressed Delaunay frames.

In \autoref{sec:delgrowth},
we estimate the growth rate $\tau$ of the positive Iwasawa factor of
$\exp(A\log z)$ as $z\to 0$
(\autoref{thm:delgrowth1}).
This growth rate result is used in \autoref{sec:delframe-asymptotics}.

Given a holomorphic perturbation
\begin{equation*}
\xi = A\frac{dz}{z} + \Order(z^{0}) dz
\end{equation*}
of the Delaunay potential $A\,dz/z$ producing a closed
once-wrapped Delaunay surface,
let $\Phi$ satisfy $d\Phi = \Phi\xi$
on an $r$-circle with unitary monodromy around $z=0$.
%

%
In \autoref{sec:delframe-asymptotics}
we show that the unitary and positive factors of $\Phi$
are asymptotic to those of a holomorphic Delaunay frame
(\autoref{thm:frame1}).

In \autoref{sec:immersion1},
we use this convergence to obtain $C^\infty$-convergence
of the CMC end to a Delaunay surface.
The convergence of the positive part implies that of the metric,
and together with the frame convergence,
a bootstrap argument on the Gauss equation gives the
$C^\infty$-convergence.
If the base Delaunay surface is embedded, then the asymptotic
end is embedded and has exponential convergence.
These results are summarized in \autoref{thm:main1}.

In Part~2 we generalize these results to the case of
dressed holomorphic frames.
\end{xtext}

\section{The Delaunay frame and its growth}


\label{Sec:delaunay}
\typeout{=================delaunay}
\subsection{The Delaunay residue}
\label{sec:delres}

\begin{xtext}
A Delaunay surface is described by a holomorphic potential 
on $\bbC^*$ of a very simple kind.
Our description of these potentials is in the setting
of~\cite{Schmitt_Kilian_Kobayashi_Rossman_2006}.
\end{xtext}

\begin{proposition}
\label{thm:delres}
\theoremname{Delaunay residue}
Let $A:\bbS^1\to\matsl{2}{\bbC}$ be analytic.
Then the following are equivalent:

\emph{(i)}
$A^* = A$, and $A$ is the boundary of a meromorphic map
$A:\calD_1\to\matsl{2}{\bbC}$
such that $A$ is holomorphic on $\calD_1\setminus\{0\}$,
the upper-right entry of $A$ has a simple pole or no pole at $0$,
and the other entries of $A$ do not have poles at $0$.

\emph{(ii)}
There exist $a,b\in\bbC$ and $c\in\bbR$ such that 
\begin{equation} \label{eq:delres}
A = \begin{pmatrix}
c & a\lambda^{-1} + \ol{b} \\ b + \ol{a}\lambda & -c
\end{pmatrix}.
\end{equation}
\end{proposition}

\begin{proof}
If $A$ is of the form \eqref{eq:delres}, then (i) clearly holds.

Conversely, suppose (i) holds.
Since $A$ extends meromorphically to $\calD_1$, and
$A^* = A$, then $A$ extends meromorphically to $\CPone$.
The entries of $A$ are then meromorphic functions on 
$\CPone$, from which it follows that $A$ must be of the form \eqref{eq:delres}.
\end{proof}

\begin{xtext}
After a rigid motion, the harmonic Gauss map 
of any Delaunay 
surface is framed by the $r$-unitary part of the
solution
$\Phi(z) = \exp(z\,A)$
of the ODE $d\Phi = \Phi A\, dz$, for some $A$ as in 
\eqref{eq:delres} (see~\cite{Burstall_Kilian_2005,Schmitt_Kilian_Kobayashi_Rossman_2006}).
This prompts us to make the following
\end{xtext}

\begin{definition} 
\label{def:delres}
\theoremname{Delaunay residue}
A \emph{Delaunay residue} is a meromorphic $\matsl{2}{\bbC}$-valued
matrix map as in \eqref{eq:delres},
with $a,\,b\in\bbC^\ast$ and $c \in \bbR$.
Let $\mu:\bbC^\ast\to\bbC$ be an eigenvalue of
$A$ satisfying $\Real\mu\ge 0$.
\Comment{%
The term ``Delaunay residue'' is perhaps misleading:
it suggests that the matrix produces a closed Delaunay surface
at $\lambda=1$ in the Sym formula.
In fact, we make no such assumptions, either that the
surface is closed, or that it is not a twizzler.}
\end{definition}

\Depend{
\begin{proof}
This definition uses \autoref{thm:delres}.
\end{proof}
}

\begin{xtext}
It can be shown~\cite{Kilian_McIntosh_Schmitt_2000, Burstall_Kilian_2005}
that up to rigid motions, 
any Delaunay surface in $\bbR^3$ can be obtained by an off-diagonal Delaunay 
residue with real non-zero parameters $a,\,b \in \bbR^*$
satisfying the closing condition $a+b=1/2$, and 
that the resulting \emph{necksize} of the surface 
depends on the product $ab$: 
when $ab >0$, the resulting surface is an unduloid, 
when $ab < 0$, it is a nodoid, and
when $a=b$, the resulting surface is a round cylinder
(see~\cite{Kilian_2004, Schmitt_Kilian_Kobayashi_Rossman_2006}).

For $a,\,b\in\bbC^\ast$,
we denote as the \emph{vacuum} the case $\abs{a}=\abs{b}$ and $c=0$.
\Comment{Section \autoref{sec:del-gauge} discusses gauge-equivalence
of Delaunay potentials.}
\end{xtext}

\typeout{=================delaunay-frame}
\subsection{The Delaunay frame}
\label{sec:delframe}

\begin{xtext}
The unitary frames for all CMC tori are computed
in~\cite{Bobenko_CMC_1991, Bobenko_Tori_1991} in terms of theta functions.
In the case of spectral genus 1,
the Iwasawa factorization of $\exp( (x+iy) A) )$
can be expressed in terms of elliptic functions and elliptic integrals,
and has a Floquet form whose period is that of the induced Delaunay surface.
The computation of this factorization in \autoref{thm:delframe}
requires the following sets and functions.

\end{xtext}

\begin{notation}
\label{not:mu-domain}
\theoremname{Delaunay residue domain}
Given a Delaunay residue $A$,
let $\nu_1,\,\nu_2\in\bbC$ be the zeros of $\det A$,
with $\abs{\nu_1}\le\abs{\nu_2}$.
Let $p\in\bbS^1$ be the point on the straight line segment
with endpoints $\nu_1$ and $\nu_2$,
or $p=\nu_1=\nu_2$ in the case of the vacuum.
Let $\alpha=-p$.
Define the following subsets of $\bbC$
(see \autoref{fig:domain}):%
\begin{equation*}
\delIA = \{rp\in\bbC\suchthat 0\le r < \infty\},\quad
\delJA = \{rp\in\bbC\suchthat \abs{\nu_1}\le r\le \abs{\nu_2}\},\quad
\delKA = (\delIA\setminus\delJA)\cup\{\nu_1,\,\nu_2\}
\spaceperiod
\end{equation*}
The set of \emph{resonance points} for $A$ is
\begin{equation}
\label{eq:S_A}
\delSA
 =
 \{\lambda\in\bbC^\ast\suchthat \mu(\lambda)\in\half\bbZ^\ast\}\spaceperiod
\end{equation}
\Comment{%
With the coefficients $a,\,b,\,c$ of the Delaunay residue
$A$ as in \autoref{thm:delres},
then $I_A\cap\bbS^1 = -(ab)/\abs{ab}$.}
%
\end{notation}

\Depend{
\begin{proof}
This notation uses \autoref{def:delres}.
\end{proof}
}

\begin{notation}
\label{not:delv}
\theoremname{Delaunay v and sigma}
Let $A$ be a Delaunay residue
and let $a,\,b,\,c$ be its coefficients as in \autoref{thm:delres}.
Define $v:\bbR \to \bbR_{>0}$
as the elliptic function satisfying
\begin{equation}
\label{eq:delv}
  (v')^2 =  -v^4 + 4(\abs{a}^2+\abs{b}^2+c^2)v^2 - 16 \abs{ab}^2\spacecomma
  \quad v(0)=2\abs{b}\spacecomma
\end{equation}
taking the non-constant solution except in the case of the vacuum.
When $c \neq 0$, $v'(0)$ is taken to have the same sign as $-c$.
The function $v$ is the restriction to $\bbR$ of an elliptic function
on $\bbC$ with a real and a pure imaginary period;
let $\rho\in\bbR_+$ be its real period.
\Comment{See \autoref{sec:delv-notes} for details on $v$.
\begin{itemize}
\item
\autoref{lem:jacobi-ode} defines the Jacobi elliptic functions.
\item
\autoref{lem:jacobi-identities} gives some identities involving
the Jacobi elliptic functions.
\item
\autoref{thm:delv-elliptic} writes $v$ explicitly as an elliptic
function and computes its period as a complete elliptic
integral.
\end{itemize}
The metric of the Delaunay surface with extended frame
$F$ is $4\abs{abH^{-1}}v^{-1}$.}

Define the elliptic integral of the third kind
$\psi:\bbR\times(\bbC\setminus\delJA)\to\bbC$ by
\begin{equation}
\label{eq:delf}
\delf(x,\,\lambda) =
\int_0^x
  \frac{2\,dt}{1 + (4 \ol{a}\ol{b} \lambda)^{-1} v^2(t)}
\spacecomma
\end{equation}
and define $\sigma:\bbC\setminus\delJA\to\bbC$ by
$\sigma(\lambda) = \psi(\rho,\,\lambda)$.
\Comment{See \autoref{sec:psi-notes} for more on $\delf$.
\autoref{thm:elliptic3} writes $\delf$ as an elliptic integral of the
third kind.
This Mathematica code
\textattachfile[%
  subject={Delaunay elliptic function and integral},
  description={Mathematica code computing elliptic function and integral.}
  ]{attachments/elliptic.m}{elliptic.m}
computes $v$ and $\dels$.}

\end{notation}

\begin{theorem}
\label{thm:delframe}
\theoremname{Delaunay factorization}
Let $A$ be a Delaunay residue and let
$\Phi = \exp((x+iy)A)$.
Then there exists an analytic map
$\delR:\bbR\times(\bbC\setminus\delJA)\to\matSL{2}{\bbC}$
satisfying $\delR(x+\rho,\,\lambda) = \delR(x,\,\lambda)$,
such that, restricting to $\bbC\setminus\delJA$,
\begin{subequations}
\label{eq:delFB}
\begin{align}
\label{eq:delF}
\Uni{1}{\Phi} &= \exp( (x+iy-\rho^{-1}\dels x) A) \delR
\spacecomma
\\
\label{eq:delB}
\Pos{1}{\Phi} &= \delR^{-1} \exp( \rho^{-1}\dels x A)
\spaceperiod
\end{align}
\end{subequations}
\end{theorem}

\begin{proof}
Let prime denote the derivative with respect to $x$.
With $h=\diag((b/\abs{b})^{1/2},\,(b/\abs{b})^{-1/2})$,
define
$\delR(x,\,\lambda):\bbR\times(\bbC\setminus\delJA)\to\matSL{2}{\bbC}$
by
\begin{gather*}
R(x,\,\lambda) = \exp( (\rho^{-1}\dels(\lambda)x - \delf(x,\,\lambda)) A )S(x,\,\lambda)
\spacecomma\\
S_1(x,\,\lambda) = 
\begin{pmatrix}
v^2 + 4\ol{a}\ol{b}\lambda  & v' + 2cv \\
  0 &  2(b+\ol{a}\lambda)v
\end{pmatrix}h
\spacecomma
\quad
S(x,\,\lambda) = (\det S_1)^{-1/2}S_1
\spaceperiod
\end{gather*}
The square root
$(\det S_1)^{1/2}$ can be taken to be a single-valued analytic
function in $\lambda$ on $\bbR\times(\bbC\setminus\delJA)$,
its sign chosen so that $S(0,\,\lambda)=\id$.
%
Then $R$ is periodic in $x$ with period $\rho$ because each of the two
factors defining it are.

On $\bbC\setminus\delJA$, define
$F = \exp( (1-\rho^{-1}\sigma)x A)R$ and
$B = R^{-1}\exp(\rho^{-1}\sigma x A)$.
Then
\begin{subequations}
\begin{align}
\label{eq:F-ODE}
F^{-1}F' &= \gauge{ ((1-\psi')A) }{S} =
\begin{pmatrix}
0 &  - \frac{v}{2} + \frac{2 a b}{\lambda v} \\
 \frac{v}{2} -\frac{2 \ol{a}\ol{b}\lambda}{v} & 0
\end{pmatrix} .\, h =:\theta
\spacecomma\quad
F(0,\,\lambda)=\id\spacecomma
\\
\label{eq:B-ODE}
B'B^{-1} &= -\gauge{(-\psi'A)}{S}
=
 \begin{pmatrix}
  -\frac{v'}{2v} & v \\ \frac{4 \ol{a}\ol{b}\lambda}{v} & \frac{v'}{2v}\end{pmatrix} .\, h =:\eta
\spacecomma\quad
B(0,\,\lambda)=\id
\spaceperiod
\end{align}
\end{subequations}
\Comment{ This Mathematica code
\textattachfile[%
  subject={The Delaunay frame},
  description={Mathematica code checking Delaunay frame calculations.}
  ]{attachments/delframe.m}{delframe.m}
checks the above two equations.}

Since $\theta$ is analytic in $\lambda$ on $\bbC^\ast$,
and $\theta^\ast = -\theta$,
then $F$ extends to a map $\bbR\to\LoopuSL{1}$.
Similarly, $\eta$ is analytic in $\lambda$ on $\bbC$,
and $\left.\eta\right|_{\lambda=0}$ is upper-triangular with
real diagonal entries,
so $B$ extends to a map $\bbR\to\LooppSL{1}$.
Hence $\Phi = \exp(iyA)F\cdot B$ is the $1$-Iwasawa factorization of $\Phi$.
\Depend{
This lemma uses \autoref{not:mu-domain}.}
\end{proof}

\Comment{For the case of the vacuum, see \autoref{thm:vacuum}.}

\begin{corollary}
\label{lem:del-period}
\theoremname{Delaunay periodicity}
Let $A$ be a Delaunay residue
and let $\sigma$ be as in \autoref{not:delv} for this $A$.
Let $\exp((x+iy)A) = F\cdot B$ be the $1$-Iwasawa factorization.
Then
 $\exp(\sigma A)$ extends to a holomorphic function of $\lambda$ on $\bbC$,
and the following quasiperiodicity
formulas hold for all $x\in\bbR$ and $n\in\bbZ$:
\begin{subequations}
\begin{align}
\label{eq:periodF}
F(x + n\rho,\,y) &= \exp( n(\rho-\dels)A )F(x,\,y)\spacecomma
\quad
\lambda\in\bbC^\ast\spacecomma\\
\label{eq:periodB}
B(x + n\rho) &= B(x)\exp( n \dels A )\spacecomma
\quad
\lambda\in\bbC
\spaceperiod
\end{align}
\end{subequations}
\end{corollary}

\begin{proof}
Since $\Phi=FB$, and $\Phi$ and $B$ are holomorphic in $\lambda$
on $\calD_1^\ast$, then $F$ on $\bbS^1$ extends holomorphically in $\lambda$
to $\calD_1^\ast$. Since $F^\ast = F^{-1}$, then $F$ extends holomorphically
in $\lambda$ to $\bbC^\ast$. Since $\Phi$ is holomorphic in $\lambda$
on $\bbC^\ast$, then $B$ extends holomorphically in $\lambda$ on
$\bbC$.
(This can also be seen from the ODE~\eqref{eq:B-ODE} satisfied by $B$.)
But by~\eqref{eq:delB},
$B(\rho)$ and $\exp(\sigma A)$ are equal
on $\bbC\setminus\delJA$.
The quasiperiodicity formulas for $F$ and $B$ follow
by~\eqref{eq:delFB}.
\Depend{%
The theorem depends on \autoref{thm:delframe}.}
\end{proof}

%

\typeout{=================tau}
\subsection{The Delaunay growth exponent}
\label{sec:tau}
\begin{xtext}
Our proof of the convergence of the $r$-Iwasawa factors
of a holomorphic Delaunay frame $\Phi=\exp(A\log z)$ and a perturbation of it
requires growth bounds on $\Pos{r}{\Phi}$.
Being an exponential, $\Phi$ grows exponentially, and its growth rate
is the absolute value of the real part of an eigenvalue $\mu$ of
the Delaunay residue $A$.
On $\bbS^1$ ($r=1$), the factor $\Pos{r}{\Phi}$ has the same growth behavior
as $\Phi$, because $\Uni{r}{\Phi}$ does not grow.
For $r<1$, the Floquet behavior of $\Pos{r}{\Phi}$ again
implies exponential growth, its rate determined by
the eigenvalues of its value after one period.

We begin by studying the real part $\tau$ of these eigenvalues,
showing it is less than the real part of the eigenvalue $\mu$ of $A$.
\end{xtext}

%
\begin{lemma}
\label{thm:tau}
\theoremname{tau lemma}
Let $A$ be a Delaunay residue
and let $\mu$ be its eigenvalue as in~\eqref{def:delres}.
Let $\sigma$ be as in \autoref{not:delv} for this $A$.
Then the function $\tau=\rho^{-1}\Real\mu\sigma$
extends to a single-valued continuous function on $\bbC$,
which is real analytic and harmonic on $\bbC\setminus\delKA$.
Moreover, with $\delKA$ as \autoref{not:mu-domain},
$\tau=0$ on $\delKA$,
$\tau>0$ on $\calD_1\setminus\delKA$,
and $\tau=\mu$ on $\bbS^1$.
\end{lemma}

\begin{proof}
\STARTSTEPS
We first consider the nonvacuum case.

\STEP
We show that
  $\sigma$ extends holomorphically as a function of $\lambda$ along any curve
  in $\bbC\setminus\{\nu_1,\,\nu_2\}$.
Let $J(t,\,\lambda)dt$ be the integrand defining $\delf$ in~\eqref{eq:delf}.
Then $J$ can be extended to a meromorphic function of $t$
on $\bbC$, by considering $v$ as an elliptic function on $\bbC$.
Let $\gamma$ be a curve with endpoints $0$ and $\rho$ which
is homotopic to the straight line $[0,\,\rho]$ in
$\bbC\setminus\sing{v}$.
Let
$\Lambda=\{\lambda\in\bbC\suchthat J(t,\,\lambda)\ne\infty\text{ along }\gamma\}$.
Then we can define $\tilde{\sigma}(\lambda)$ on $\Lambda$
as the integral of $J$ along $\gamma$.
Since $\tilde\sigma$ is analytic, and is equal to $\sigma$
on the intersections of their domains $\Lambda\cap(\bbC\setminus\delJA)$,
then $\tilde\sigma$ is an analytic extension of $\sigma$.
This provides a construction for analytically extending
$\sigma$ along any curve in $\bbC\setminus\{\nu_1,\,\nu_2\}$.

\STEP 
Since for all $x\in\bbR$, $B(x,\,\lambda)=\Pos{1}{\exp(xA(\lambda))}$ is
a holomorphic function of $\lambda$ on $\bbC$,
then $B(\rho,\,\lambda)$ is holomorphic in $\lambda$.
Hence $\cosh(\mu\sigma) = \half\tr B(\rho,\,\lambda)$
is analytic on $\bbC$.
An argument shows
that at $\lambda\in\bbC$,
if $\cosh(\mu\sigma)\not\in\{\pm 1\}$, $\mu\sigma$ is analytic at $\lambda$,
and if $\cosh(\mu\sigma)\in\{\pm 1\}$,
$\mu\sigma$ is of the form $\pi i k+\mu g$ for some $k\in\bbZ$ and
some holomorphic function $g$ near $\lambda$.
\Comment{See \autoref{thm:mu-sigma-behavior}.}

\STEP
Near $\lambda=0$, $\sigma$ is analytic. Hence if
any branch of $\mu\sigma$ is
analytically extended along a closed once-wrapped curve around $0$,
with respective values $p_0$ and $p_1$ at the beginning and end of the curve,
then $p_1=-p_0$.
Likewise, if any branch of $\mu\sigma$ is
analytically extended along a closed curve around $\nu_1$,
with respective values $q_0$ and $q_1$ at the beginning and end of the curve,
it follows from step 2 that $q_1\equiv -q_0\mod 2\pi i$.

Putting these together, if $\mu\sigma$ is analytically
extended along a closed once-wrapped curve around $0$ and $\nu_1$,
with respective values $r_0$ and $r_1$ at the beginning and end of the curve,
then $r_1\equiv r_0\mod 2\pi i$.
Hence $\tau$ is a single-valued real analytic function on
$\bbC\setminus\delKA$.
It is harmonic there because it is locally the real part of an
analytic function.

\STEP
To show that $\tau$ is continuous on $\delKA\setminus\{\nu_1,\,\nu_2\}$
with value $0$, write $\mu$ and $\sigma$ in terms of their
real and imaginary parts
$\mu = \mu_1 + i\mu_2$ and $\sigma = \sigma_1+i \sigma_2$.
Then $\Real\mu\sigma = \mu_1\sigma_1 -\mu_2\sigma_2$.
Since $\sigma_2$ and $\mu_1$ are $0$ on $\delKA$,
$\tau$ is continuous on $\delKA\setminus\{\nu_1\,\,\nu_2\}$
with value $0$ there.
\Comment{
Since $\sigma_2$ is $0$ on $\delKA$, then
$\mu_2\sigma_2$ is continuous there.
Hence $\Real\mu\sigma$ is continuous on $\delKA$ with value $\mu_1\sigma_1$,
which is $0$ on $\delKA$ because $\mu_1$ is $0$ there.}

\STEP
We now show that $\tau$ is continuous at $\nu_1$ and $\nu_2$
with value $0$.
By step 2, with $k\in\{1,\,2\}$,
we have $\tau = \Real(c + \mu g) = \Real(\mu g)$ near $\nu_k$,
for some $c\in \pi i \bbZ$ and holomorphic function $g$.
Since $\mu(\nu_k)$ is continuous at $\nu_k$ and $\mu(\nu_k)=0$, then
$\tau$ extends continuously to $\nu_k$ with value $0$ there.

\STEP
To show that $\tau=\mu$ on $\bbS^1$,
since $F$ in~\eqref{eq:periodF}
takes values in $\matSU{2}{}$ on $\bbS^1$,
then by~\eqref{eq:periodF}, so does $\exp( (\rho-\dels) A)$.
Since $A$ is tracefree and is
hermitian on $\bbS^1$, then $\rho-\dels$ is pure imaginary
on $\bbS^1$. Hence $\tau=\mu$ on $\bbS^1$.

\STEP
Since $\tau$ is continuous on $\ol{\calD_1}$
and harmonic on $\calD_1\setminus\delKA$,
then by the maximum principle for harmonic functions,
$\tau$ on $\calD_1$ attains its minimum on the boundary $\bbS^1\cup\delKA$
of $\ol{\calD_1}\setminus\delKA$.
Since $\tau$ is $0$ on $\delKA$, then
$\tau$ is strictly positive on $\calD_1\setminus\delKA$.

For the vacuum case,
let $a,\,b$ be the coefficients of the Delaunay
residue as in \autoref{def:delres},
and let $\alpha$ be as in \autoref{not:mu-domain}.
Then $\tau=\Real\rho^{-1}\mu\sigma$ satisfies all the properties of the theorem,
where
\begin{equation}
\label{eq:sigma-vacuum}
\mu = \abs{b}\alpha^{-1/2}\lambda^{-1/2}(\lambda+\alpha)
\AND
\rho^{-1}\mu\sigma = 2\abs{b} \alpha^{-1/2}\lambda^{1/2}
\spaceperiod
\qedhere 
\end{equation}
\end{proof}


\begin{figure}[ht]
\centering
\subfigure[%
The $\lambda$-plane with the zeros $\nu_1$ and $\nu_2=1/\ol{\nu_1}$
of the eigenvalues of a typical Delaunay residue.
]{\label{fig:domain}%
\includegraphics[width=2.5in]{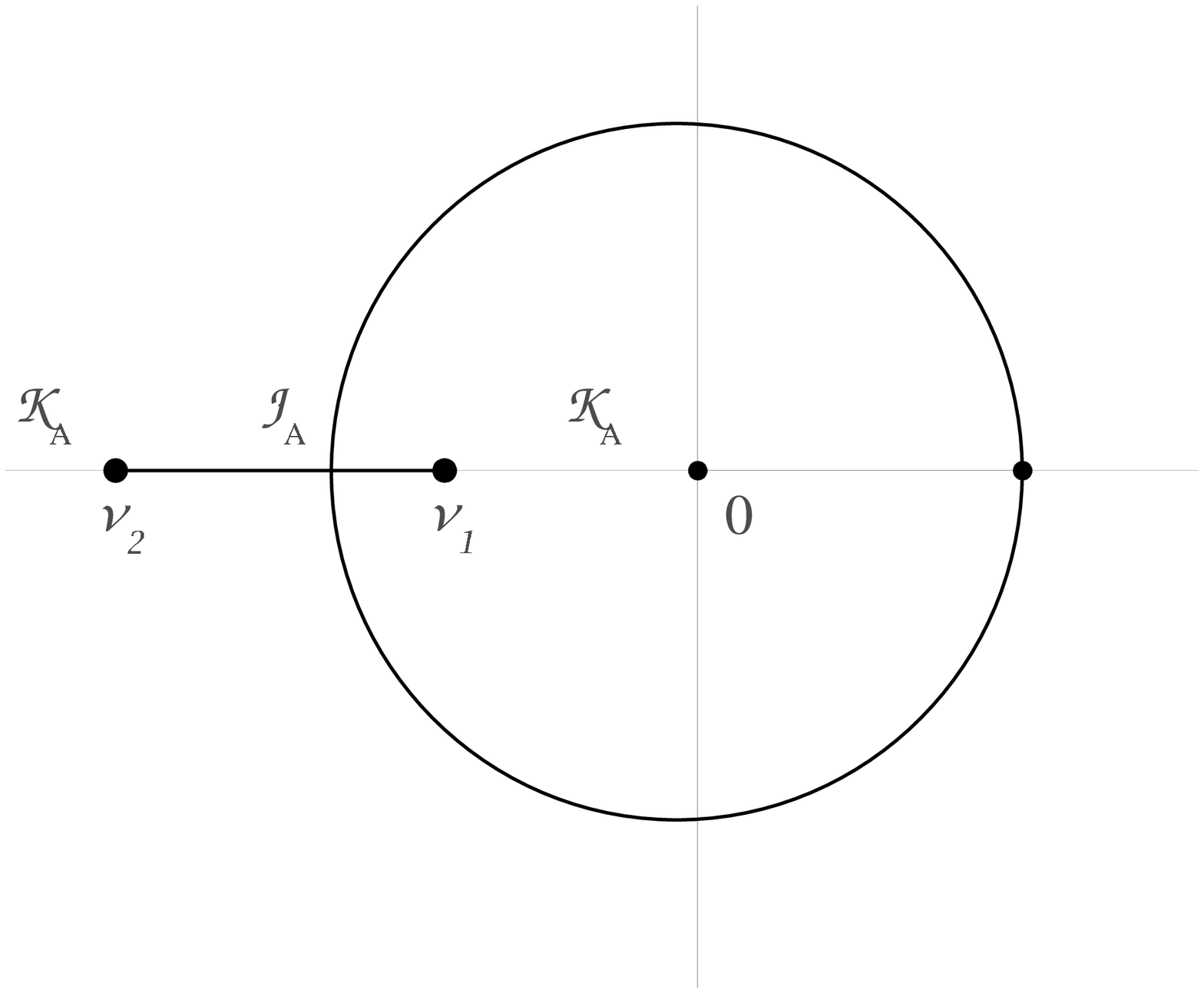}}
\hspace{0.1in}
\subfigure[Graphs of $\Real \mu$ (above) and $\tau$ (below)
over a half-disk for a typical Delaunay frame.
$\Real \mu$ represents the growth of $\exp(A\log z)$,
$\tau$ the growth of its positive factor,
and their difference the growth of its unitary factor.
%
]{\label{fig:graph}%
\includegraphics[width=2.5in]{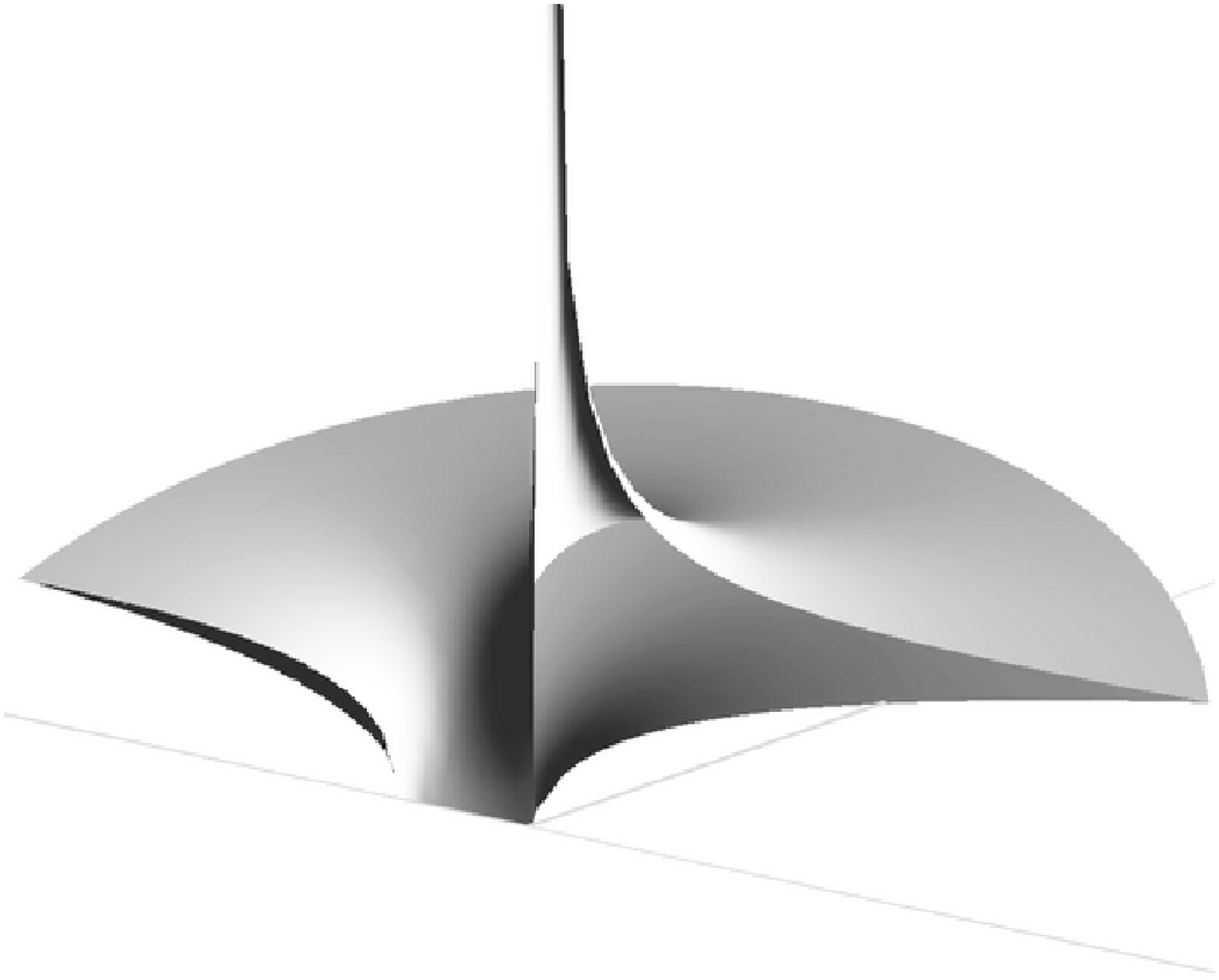}}
\caption{}
\label{fig:diagrams}
\end{figure}


\begin{lemma}
\label{thm:mu-tau}
\theoremname{mu-tau lemma}
Let $A$ be a Delaunay residue, let $\mu$
be its eigenvalue as in \autoref{def:delres},
and let $\tau$ be as in \autoref{thm:tau} for this $A$.
Then
$\tau \in (0,\,\Real\mu]$ on $\calD_1\setminus\delKA$.
\end{lemma}

\begin{proof}
In the case of the vacuum,
the result follows by~\eqref{eq:sigma-vacuum};
hence we assume the nonvacuum case.
(See \autoref{fig:graph} for the graphs of
$\Real(\mu)$ and $\tau$ for a typical Delaunay residue.)

We first show that $\Real(\mu)-\tau$ is \nonnegative on $\calD_r^\ast$
for some $r$ near zero,
and then apply the maximum principle for harmonic functions
to conclude it is strictly positive on $\calD_1\setminus\delKA$.

\STARTSTEPS
\STEP
We prove the following claim:
Let $r>0$ and let $f:\calD_r\to\bbC$ be a holomorphic function.
If $f$ has the symmetry $f(\ol{\lambda})= \ol{f(\lambda)}$,
and $f(0)\in\bbR_+$,
then there exists $r_1\in(0,\,r]$ such that
$g:=\Real(\lambda^{-1/2}f)$ is \nonnegative
on $\calD_{r_1}^\ast$, where the square root is chosen to have \nonnegative
real part.

To prove the claim, let $z=x+iy = \lambda^{1/2}$, and let
\begin{equation*}
h
 := \frac{g}{\Real \lambda^{-1/2}}
 = \Real f(z^2) + \frac{y}{x}\Imag f(z^2)
\spaceperiod
\end{equation*}
The function $k(x,\,y)=\Imag f(z^2)$ is real analytic in $x$ and $y$.
By the symmetry of $f$,
$f(-y^2)$ is real, so for all $y$, $k(0,\,y) = \Imag f(-y^2) = 0$.
It follows that $k(x,\,y)/x$ is real analytic in $x$ and $y$.
Then $\lim_{z\to 0} yk(x,\,y)/x = 0$,
so $\lim_{z\to 0} h = f(0)\in\bbR_+$.
Hence there exists $r_1\in(0,\,r]$ such that
$h$ is strictly positive in $\calD_{r_1}$.
Since $\Real \lambda^{-1/2}$ is strictly positive on 
$\calD_{r_1}\setminus \bbR_{\le 0}$,
then so is $g$.
But $g$ is $0$ on $\bbR_-$,
so $g$ is \nonnegative on $\calD_{r_1}^\ast$.

\STEP
Let $\sigma$ be as in \autoref{thm:tau}
on $\calD_{\abs{\nu_1}}$, so $\tau = \Real\sigma$.
With $\alpha$ as in \autoref{not:mu-domain},
define $f:\calD_{\abs{\nu_1}}\to\bbC$ by
\begin{equation*}
f = \lambda^{1/2}\left(\mu(\alpha\lambda) - \sigma(\alpha\lambda)\right)
\spaceperiod
\end{equation*}
Then $f$ is holomorphic on $\calD_{\abs{\nu_1}}$,
$f(\ol{\lambda}) = \ol{f(\lambda)}$ and $f(0)\in\bbR_+$.
By the above claim applied to $f$,
there exists $r_1\in(0,\abs{\nu_1}]$ such that
$\Real \lambda^{-1/2} f = \Real(\mu(\alpha\lambda))-\tau(\alpha\lambda)$
is \nonnegative on $\calD_{r_1}^\ast$.
Hence
$\Real(\mu)-\tau$ is \nonnegative $\calD_{r_1}^\ast$.

\STEP
For any $s\in(0,\,r_1]$, define
$V_s$ as the union of $\calC_s$ and the straight line segment along $\delKA$
from $\nu_1$ to this circle.
Then for all $s<r_1$,
$\Real(\mu)-\tau$ is harmonic on the open region $R_s$
between $V_s$ and $\bbS^1$,
is $0$ on $\bbS^1$ and $V_s\cap\delKA$, and is \nonnegative on $\calC_s$
by step 2.
By the maximum principle for harmonic functions,
$\Real(\mu)-\tau$ is strictly positive on $R_s$.
Since this is true for any $s\in(0,\,r_1)$,
then $\Real(\mu)-\tau$ is strictly positive on $\calD_1\setminus\delKA$.
%
\end{proof}

\typeout{=================delaunay-growth1}
\subsection{Growth of the Delaunay positive part}
\label{sec:delgrowth}

\begin{xtext}
Above we investigated the function $\tau$, the real part of the
eigenvalue of the value of $\Pos{r}{\exp(A\log z)}$ after one
period, showing that it is less than the real part of the eigenvalue
$\mu$ of the corresponding Delaunay residue.
The Floquet behavior of $\Pos{r}{\exp(A\log z)}$, detailed in
\autoref{thm:delframe}, implies that $\tau$, and hence $\mu$,
bounds the exponential growth of $\Pos{r}{\exp(A\log z)}$
(\autoref{thm:delgrowth}).
\end{xtext}

\begin{notation}
\label{not:norm}
\theoremname{Matrix norms}
For $v \in \bbC^2$, we denote the 
{\em vector norm} by 
$\abs{v} = \sqrt{\transpose{\ol{v}}v}$, and 
for $M\in\mattwo(\bbC)$, we set 
\begin{equation*}
\matnorm{M}  = \max_{\abs{v}=1} \abs{Mv}\spaceperiod
\end{equation*}
For a map $M:\calR\to\mattwo(\bbC)$ on a subset $\calR\subset\bbC$
we set
\begin{equation*}
\supnorm{M}{D} = \sup_{\lambda\in\calR} \,\matnorm{M(\lambda)}\spaceperiod
\end{equation*}

\Comment{%
Section \autoref{sec:norm} considers this norm in detail.
\begin{itemize}
\item
See \autoref{def:norm} for the definitions of vector and matrix norm.
\item
See \autoref{lem:max-norm} for a proof that the $\sup$ norm is a
matrix norm.
\item
See \autoref{lem:induced-2-norm} for a computation of the $\sup$ matrix norm
in terms of singular values.
\item
See \autoref{thm:norm-of-inverse} for some properties of this matrix
norm with respect to unitary matrices and inverses.
\end{itemize}}

Given $r\in(0,1]$, a subset $\calR\subset\bbC$
containing $\calC_r$, and a loop $X\in\LoopSL{r}$
which extends to a map $Y:\calR\to\matSL{2}{\bbC}$,
by an abuse of notation
we will write $\supnorm{X}{\calR}$ for $\supnorm{Y}{\calR}$.
\end{notation}

\begin{lemma}
\label{thm:exp}
\theoremname{exp}
Let $X:\calR\to\matsl{2}{\bbC}$ be a continuous map
on a domain $\calR\subseteq\bbC$,
and let $\mu:\calR\to\bbC$ be any eigenvalue function of $X$.
Then there exists a continuous function $c:\calR\to\bbR_+$ such that
$\matnorm{\exp X}
 \le c e^{\abs{\Real \mu}}.$
\end{lemma}

\begin{proof}
The result follows from the formula
$\exp X = \cosh(\mu)\id + \mu^{-1}\sinh(\mu) X$
and the estimates
$\abs{\cosh \mu}\le e^\abs{\Real{\mu}}$ and
$\left|\mu^{-1}\sinh \mu\right|\le e^{\abs{\Real{\mu}}}$.
\Comment{See \autoref{thm:exp-full} for a full proof.
The theorem uses~\autoref{not:norm}.}
\end{proof}

\begin{lemma}
\label{thm:delgrowth}
\theoremname{Undressed Delaunay growth}
Let $A$ be a Delaunay residue,
and let $\mu$ be its eigenvalue as in \autoref{def:delres}.
Let $\tau$ be as in \autoref{thm:tau}.
Then there exists a continuous function $c:\ol{\calD_1}\setminus\{0\}\to\bbR_+$
such that for all
 $(z,\,\lambda)\in\{0<\abs{z}<1\}\times\ol{\calD_1}\setminus\{0\}$,
\begin{equation}
\label{eq:del-growth}
\matnorm{\Pos{1}{\exp(A\log z)}}
 \le c\abs{z}^{-\tau}
 \le c\abs{z}^{-\Real\mu}\spaceperiod
\end{equation}
\end{lemma}

\begin{proof}
For $x+iy\in\bbC$, define $B(x+iy,\,\lambda)=\Pos{1}{\exp((x+iy)A)}$.
With $\rho$ as in \autoref{lem:del-period},
define $x_0:\bbR\to[0,\,\rho)$ and $n:\bbR\to\bbZ$ as the unique
functions such that $x = x_0 + n \rho$.
Let $\dels$ be as in \autoref{not:delv} for $A$.
By \autoref{lem:del-period},
\begin{equation}
\label{eq:del-growth1}
  B(x+iy)
 =
  B(x_0+n\rho)
 = B(x_0)\exp(n \dels A)\spaceperiod
\end{equation}
This quasiperiodicity of $B$ determines its growth rate as follows.

Define the continuous function $c_1:\bbC^\ast\to\bbR_+$
by
$c_1(\lambda) = \max_{x\in[0,\,\rho)}\matnorm{B(x,\,\lambda)}$.
Then
$\matnorm{B(x_0,\,\lambda)} \le c_1(\lambda)$
on $[0,\,\rho)\times\bbC^\ast$.
With $\tau$ as in \autoref{thm:tau} for $A$,
by \autoref{thm:exp} there exists
a continuous function $c_2:\bbC^\ast\to\bbR_+$
such that on $\bbC^\ast$,
\begin{equation*}
\matnorm{\exp(n \delf(\rho,\,\lambda) A(\lambda))}
 \le
 c_2 e^\abs{n\rho\tau(\lambda)}
 \le
 c_2 e^\abs{-x_0\tau(\lambda)} e^\abs{x\tau(\lambda)}\spaceperiod
\end{equation*}
Define the continuous function $c_3:\bbC^\ast\to\bbR_+$ by
$c_3(\lambda) =\max_{x\in[0,\,\rho)}e^{\abs{x\tau(\lambda)}}$.
Then
\begin{equation}
\label{eq:del-growth3}
\matnorm{\exp(n \delf(\rho,\,\lambda) A(\lambda))}
 \le
 c_2(\lambda)c_3(\lambda) e^\abs{\tau(\lambda)x}
\spaceperiod
\end{equation}

The quasiperiodicity~\eqref{eq:del-growth1},
the choice of $c_1$, and~\eqref{eq:del-growth3}
yield the estimate
\begin{equation*}
\matnorm{B(x+iy,\,\lambda)}
\le
c_1(\lambda)c_2(\lambda)c_3(\lambda) e^\abs{\tau(\lambda)x}\spaceperiod
\end{equation*}

Since $0<\abs{z}<1$, then $x<0$.
But $\abs{z}=e^x$,
so $\abs{z}^{-\tau} = e^{-\tau x} = e^{\abs{\tau x}}$.
The growth estimate~\eqref{eq:del-growth}
follows with $c=c_1c_2c_3$.

The second inequality in~\eqref{eq:del-growth} follows by
\autoref{thm:mu-tau}.

\Comment{%
Note that the estimate on $B$ in terms of $\tau$ holds
on $\bbC$, not just on the unit disk.}
\end{proof}

\subsection{Dressed Delaunay frames}
\begin{xtext}
The growth bound on the positive $r$-Iwasawa factor of the holomorphic
Delaunay frame $\Phi=\exp(A\log z)$,
computed above in \autoref{thm:delgrowth},
is preserved by dressing (right-multiplying) $\Phi$
by an $r$-loop $C$,
provided that the unitarity
of its monodromy is preserved.
In the special case that $C$ can be extended analytically to $\calA_{r,\,1}$,
the resulting frame $C\Phi$ again induces a Delaunay immersion
(\autoref{thm:triv-del-dress}),
and the growth of its positive factor
is bounded by the same bound as that for $\Pos{r}{\Phi}$
(\autoref{thm:delgrowth1}).

The general case of dressing by an arbitrary loop $C$ which preserves the unitarity
of the monodromy is characterized in \autoref{sec:delaunay-dressing},
where it is shown that $C\Phi$ induces a multibubbleton,
a Bianchi-B\"acklund transformed Delaunay immersion.
The same growth bound applies to the positive part of $C\Phi$
for this larger class of dressing matrices.
\end{xtext}

\begin{lemma}
\label{thm:triv-del-dress}
\theoremname{Trivial Delaunay dressing}
Let $A$ be a Delaunay residue.
Let $r\in(0,\,1]$.
Let $C\in\LoopUpSL{r}$, and
assume $C\exp(2\pi i A)C^{-1}\in\LoopuSL{r}$.
%
Let $C = C_u\cdot C_+$ be the $r$-Iwasawa factorization of $C$.
Then
\begin{enumerate}
\item
$C_+AC_+^{-1}$ extends meromorphically to $\CPone$
and is a Delaunay residue.
\item
$\Sym{r}{C\exp(A\log z)}$ and
$\Sym{r}{\exp(A\log z)}$ are Delaunay surfaces differing by a rigid motion.
\end{enumerate}
\end{lemma}

\begin{proof}
Let $M= \exp(2\pi i A)$
and let $A_1= C_+A C_+^{-1}$.
Since $CMC^{-1}\in\LoopuSL{r}$,
then $C_+MC_+^{-1} = \exp(2\pi i A_1)\in\LoopuSL{r}$.
With $\mu$ a local analytic eigenvalue of $A$ or $A_1$,
the formula
\begin{equation*}
\exp(2\pi i A_1)
 = \cos(2\pi \mu)\id + \mu^{-1}\sin(2\pi \mu)A_1
\end{equation*}
shows that $A_1$ extends meromorphically to $\bbS^1$.
Since $\exp(2\pi i A_1)\in\LoopuSL{r}$, then
$A_1^\ast = A_1$ away from its poles.
Write
\begin{equation*}
A_1 = \begin{pmatrix}x & y \\ y^\ast & -x\end{pmatrix}
\end{equation*}
for some meromorphic functions $x$ and $y$
in a neighborhood of $\bbS^1$ satisfying $x = x^\ast$.
Since $A$ is holomorphic on $\bbS^1$, then so is
\begin{equation*}
-\det A = -\det A_1 = x x^\ast + y y^\ast
\spaceperiod
\end{equation*}
Hence $x x^\ast + y y^\ast$ is bounded on $\bbS^1$.
Since $x x^\ast$ and $y y^\ast$ are each \nonnegative on $\bbS^1$,
then each is bounded on $\bbS^1$,
so each is holomorphic on $\bbS^1$.
It follows that $x$, $y$ and $y^\ast$ are holomorphic on $\bbS^1$,
so $A_1$ is holomorphic on $\bbS^1$.
Since $A_1$ is holomorphic on $\calD_1$ and satisfies $A_1^\ast = A_1$,
then $A_1$ is a Delaunay residue by \autoref{thm:delres}.

To show that $C\exp(A\log z)$ induces a Delaunay immersion,
note that
\begin{equation*}
C\exp(A\log z) = C_u\exp(A_1\log z)C_+
\spacecomma
\end{equation*}
so
\begin{equation*}
\Uni{r}{C\exp(A\log z)} = C_u\Uni{r}{\exp(A_1\log z)}
\spaceperiod
\end{equation*}
Hence $\Sym{r}{C\exp(A\log z)}$ and
 $\Sym{r}{\exp(A_1\log z)}$
differ by a rigid motion.
By Lemma 6 in~\cite{Schmitt_Kilian_Kobayashi_Rossman_2006},
$\Sym{r}{\exp(A_1\log z)}$ and $\Sym{r}{\exp(A\log z)}$
are Delaunay surfaces differing by a rigid motion.
\Comment{This lemma says there's a rigid motion and a coordinate change.}
\end{proof}

\begin{xtext}
\end{xtext}

\begin{theorem}
\label{thm:delgrowth1}
\theoremname{Dressed Delaunay growth I}
Let $A$ be a Delaunay residue,
and let $\mu$ be its eigenvalue as in \autoref{def:delres}.
Let $\tau$ be as in \autoref{thm:tau}.
Let $r\in(0,\,1]$.
Let $C\in\LoopUpSL{r}$
be such that the dressed monodromy
$C\exp(2\pi i A)C^{-1}$ around $z=0$ is in $\LoopuSL{r}$.
Then there exists a continuous function $c:\calD_1^\ast\to\bbR_+$
such that for all $(z,\,\lambda)\in\{0<\abs{z}<1\}\times\calD_1^\ast$,
\begin{equation*}
\matnorm{\Pos{r}{C\exp(A\log z)}}
 \le c\abs{z}^{-\tau}
 \le c\abs{z}^{-\Real\mu}\spaceperiod
\end{equation*}
\end{theorem}

\begin{proof}
Let $C = C_u\cdot C_+$ be the $r$-Iwasawa factorization of $C$.
By \autoref{thm:triv-del-dress},
$A_1 = C_+AC_+^{-1}$ is a Delaunay residue.
Then
\begin{equation*}
C\exp(A\log z) = C_u \exp(A_1\log z)C_+
\spacecomma
\end{equation*}
so
\begin{equation}
\label{eq:delgrowth1-1}
\Pos{r}{C\exp(A\log z)} = \Pos{r}{\exp(A_1\log z)}C_+
\spaceperiod
\end{equation}
It follows from $\det A = \det A_1$ that
the function $\tau$ in \autoref{thm:tau} for $A$
is the same as that for $A_1$.
By \autoref{thm:delgrowth},
there exists a continuous function $c_1:\ol{\calD_1}\setminus\{0\}\to\bbR_+$
such that
on $\{0<\abs{z}<1\}\times\calD_1^\ast$,
\begin{equation}
\label{eq:delgrowth1-2}
\matnorm{\Pos{r}{\exp(A_1\log z)}}
\le
 c_1\abs{z}^{-\tau}
 \le
 c_1 \abs{z}^{-\Real\mu}.
\end{equation}
Let $c_2 = \matnorm{C_+}$ on $\calD_1$.
Then by~\eqref{eq:delgrowth1-1} and~\eqref{eq:delgrowth1-2},
\begin{equation*}
\matnorm{\Pos{r}{C\exp(A\log z)}}
\le
c_1c_2\matnorm{\Pos{r}{\exp(A_1\log z)}}
\le
 c_1c_2\abs{z}^{-\tau}
 \le
 c_1c_2 \abs{z}^{-\Real\mu}
\spaceperiod
\end{equation*}
The result follows with $c=c_1c_2$.
\end{proof}

\section{The perturbed Delaunay frame}

\label{Sec:perturbed-delframe}
\typeout{=================zap}
We will use the following notation throughout the next
several sections.
Let $A$ be a Delaunay residue and let $\mu$ be
an eigenvalue of $A$ with \nonnegative real part.
Let $\delSA$ be the set of resonance points for $A$, defined in~\eqref{eq:S_A}.
Let $\Sigma\subset\bbC$ be a neighborhood of $0\in\bbC$
and let $\Sigma^\ast=\Sigma\setminus\{0\}$.
Choose $r\in(0,\,1]$.

\begin{definition}
\label{def:setup}
\theoremname{Perturbed Delaunay r-potential}
\mbox{}

\begin{enumerate}
\item
A \emph{perturbed Delaunay $r$-potential}
is an $r$-potential $\xi$ on $\Sigma^\ast$ of the form
\begin{equation*}
\xi=Az^{-1}dz + \Order(z^0)dz
\spaceperiod
\end{equation*}

\item
An \emph{$r$-gauge}
is an analytic map $\Sigma\to\LoopposSL{r}$.
\Comment{Section \autoref{sec:del-gauge} discusses various gauges
of Delaunay and perturbed Delaunay potentials.
\autoref{thm:nnoid} gives a local gauge which gauges the standard
$n$-noid potentials into perturbed Delaunay potential form.}

\item Given neighborhoods $\Sigma,\,\Sigma'\subset\bbC$ of $0\in\bbC$,
a \emph{coordinate change} $\cchange$ is a holomorphic map
$\cchange:\Sigma'\to\Sigma$ which satisfies $\cchange(0)=0$
and has a holomorphic inverse $\cchange(\Sigma)\to\Sigma'$.
\end{enumerate}
\end{definition}

\subsection{The \texorpdfstring{$z^AP$}{ZAP} lemma}
\label{sec:zap}

\begin{xtext}
Given a linear matrix ODE $d\Phi=\Phi\xi$
for which $\xi$ has a simple pole at $z=0$ and residue $A$,
a standard result
in the theory of regular singularities~\cite{Hartman_1982}
states that under certain conditions on the eigenvalues of $A$,
there exists a solution of the form $\Phi=z^A P = \exp(A\log z)P$,
where $P$ extends holomorphically to $z=0$.
\autoref{thm:parametrized-zap}
summarizes these results for our context,
in which $\xi$ depends analytically on a parameter $\lambda$.
We call the decomposition~\eqref{eq:parametrized-zap-decomposition}
the \emph{$z^AP$ decomposition}.
\end{xtext}

\begin{xtext}
The coefficients of the $z^AP$ gauge $P$ can be computed
in terms of a linear map $\ELL_n$,
whose definition and properties are given in the next lemma.
\end{xtext}


\begin{lemma}
\label{thm:L1}
\theoremname{Properties of linear map I}
Let $A\in\matgl{2}{\bbC}$
and let $\mu_1,\,\mu_2\in\bbC$ be the eigenvalues of $A$.
For $n\in{\bbZ_{\ge 0}}$, define
the linear map $\ELL_n:\matgl{2}{\bbC}\to\matgl{2}{\bbC}$ by
\begin{equation*}
\ELL_n(X) = nX+[A,\,X]\spaceperiod
\end{equation*}
%
Then:
\begin{enumerate}
\item
The eigenvalues of $\ELL_n$ are $n$, $n$, $n+\mu_1-\mu_2$, $n-\mu_1+\mu_2$.
Hence $\ELL_n$ is invertible if and only if $n\ne 0$ and
$\mu_1-\mu_2\not\in\{n,\,-n\}$.
\item
Suppose $n\in\bbZ_+$.
Let $R = n\id + A - \widehat{A}$,
where $\widehat{A}$ denotes the adjugate of $A$.
Then $\ELL_n$ is invertible if and only if $R$ is invertible,
and in this case, $\ELL_n^{-1}$ is given by
\begin{equation}
\label{eq:zap_Linv}
n \ELL_n^{-1}(X) = X - R^{-1}[A,\,X]
\spaceperiod
\end{equation}

\item
For any $X\in\matgl{2}{\bbC}$,
$\tr(\ELL_n X) = n\tr X$.
Hence if $\ELL_n$ is invertible,
then for any $Y\in\matgl{2}{\bbC}$,
$\tr Y = n\tr(\ELL_n^{-1} Y)$.

\item
Let $n\in\bbZ_+$.
If $A$ is holomorphic (respectively~meromorphic) on some domain $\calR$,
then $\ELL_n$ is holomorphic (respectively~meromorphic) on $\calR$.
In this case, if $\mu_1-\mu_2$ is not identically $n$ or $-n$, then
$\ELL_n^{-1}$ extends meromorphically to $\calR$.
\end{enumerate}
\end{lemma}

\begin{proof}
Statement (i) follows from the fact that the
eigenvalues of $\ad_A$ are $0,\,0,\,\mu_1-\mu_2,\,-\mu_1+\mu_2$.
\Comment{see \autoref{thm:eigenvalues-AX} for a proof.}

To prove (ii), since the eigenvalues of $R$ are
$n+\mu_1-\mu_2$, $n-\mu_1+\mu_2$,
the by (i), $R$ is invertible if and only if $\ELL_n$
is invertible. In this case,
\begin{align*}
\ELL_n\left(X-R^{-1}[A,\,X]\right)
&=
nX + \left(\id - nR^{-1} - A R^{-1}\right)[A,\,X] + R^{-1}[A,\,X]A \\
&=
nX + R^{-1}\left(\left(R - n\id - RAR^{-1}\right)[A,\,X] + [A,\,X]A\right) \\
&=
nX + R^{-1}\left(-\widehat{A}[A,\,X] + [A,\,X]A\right)
=
nX
\spaceperiod
\end{align*}
\Comment{See \autoref{lem:L-inverse}.}

Statement (iii) is clear from the definition of $\ELL_n$
and the fact that $\tr[A,\,X]=0$.

Identifying $\matgl{2}{\bbC}$ with $\bbC^4$,
the entries of the $4\times 4$ matrices for $\ELL_n$ and $\ELL_n^{-1}$
are rational functions of the entries of $A$,
and hence are meromorphic. This proves statement (iv).
\end{proof}

\Comment{
The setup for \autoref{thm:parametrized-zap} is as follows:
\begin{itemize}
\item
Let $A$ be a Delaunay residue (\autoref{def:delres}).
\item
Let $\mu$ be its eigenvalue as in \autoref{def:delres}.
\item
Let $r\in(0,\,1]$.
\item
With $\delSA$ as in~\eqref{eq:S_A},
assume $\calC_r\cap \delSA\ne\emptyset$.
\item
Let $\Sigma\subset\bbC$ be a neighborhood of $z=0$.
\item
Let $\Sigma^\ast = \Sigma\setminus\{0\}$.
\item
Let $\xi$ be an perturbed Delaunay $r$-potential (\autoref{def:setup})
on $\Sigma^\ast$ of the form
\begin{equation*}
\xi=Az^{-1}dz + \Order(z^0)dz
\spaceperiod
\end{equation*}
\item For (ii):
Let $\widetilde{\Sigma^\ast}\to\Sigma^\ast$ be the universal
cover of $\Sigma^\ast$.
\item
Let $\Phi:\widetilde{\Sigma^\ast}\to\LoopSL{r}$
be a holomorphic solution to the equation $d\Phi = \Phi\xi$
\end{itemize}
}

\begin{lemma}
\label{thm:parametrized-zap}
\theoremname{ZAP lemma}
Let $\calR = \bbC^\ast \setminus\delSA$.

\begin{enumerate}
\item
There exists a holomorphic solution
$P:\Sigma\times\calR\to\matSL{2}{\bbC}$ to the gauge equation
\begin{equation}
\label{eq:parametrized-zap-gauge}
\gauge{(A z^{-1}dz)}{P}=\xi\spacecomma\quad
P(0,\lambda)=\id\spaceperiod
\end{equation}

\item
Let $r\in(0,\,1]$ and assume $\calC_r\cap\delSA=\emptyset$.
Let $\Phi:\widetilde{\Sigma^\ast}\to\LoopSL{r}$
be a holomorphic solution to the equation $d\Phi = \Phi\xi$
on the universal cover
$\widetilde{\Sigma^\ast}\to\Sigma^\ast$ of $\Sigma^\ast$.
Then there exists
$C\in\LoopSL{r}$ such that
\Comment{on $\widetilde{\Sigma^\ast}\times\calC_r$,}
\begin{equation}
\label{eq:parametrized-zap-decomposition}
\Phi(z,\,\lambda) = C(\lambda) \exp(A(\lambda)\log z)P(z,\,\lambda)\spaceperiod
\end{equation}
\item
Moreover, if $\xi$ satisfies
$\xi=Az^{-1}dz + \Order(z^n)dz$,
then $P$ satisfies
$P = \id + \Order(z^{n+1})$.
\end{enumerate}

\end{lemma}

\begin{proof}
(i).
By the pointwise version of this lemma
(see Theorem 10.1 in~\cite{Hartman_1982}),
at each $\lambda_0\in\calR$,
there exists a unique solution $P(z,\,\lambda_0)$ to the gauge
equation~\eqref{eq:parametrized-zap-gauge}.%
\Comment{See \autoref{lem:pointwise-zap}.}
To show that $P$ is holomorphic in $\lambda$ on $\calR$, let
\begin{equation*}
\xi = A z^{-1}dz + \sum_{k=0}^\infty B_k z^k dz
\end{equation*}
and
\begin{equation*}
P = \sum_{k=0}^\infty P_k z^k\spacecomma
\quad
P_0=\id
\end{equation*}
be the respective series expansions for $\xi$ and $P$ in $z$ at $z=0$.
Then for all $k\in\bbZ_+$,
$\ELL_k$ is invertible at $\lambda_0$
by \autoref{thm:L1}(i)
and the assumption that $\calR\cap\delSA = \emptyset$,
The coefficients $P_k$ (see~\cite{Hartman_1982}) are given by
$P_k = \ELL_k^{-1}(C_k)$, where
\begin{equation*}
C_k = \sum_{i+j=k-1} P_i B_j\spacecomma
\quad
k\in\bbZ_+\spaceperiod
\end{equation*}
By \autoref{thm:L1}(iv), each $P_k$ is holomorphic in $\lambda$
at $\lambda_0$.
It follows by the absolute convergence of power series
that $P$ is holomorphic in $\lambda$ at $\lambda_0$.
\Comment{
Note that by the uniqueness of the solution to~\eqref{eq:zap-gauge-X},
$P$ is single-valued on $\calR$ even though $\calR$ is not simply connected.}

(ii)
Since $P$ satisfies~\eqref{eq:parametrized-zap-gauge}, then
the map $\Phi_0:\widetilde{\Sigma^\ast}\times\calC_r\to\matSL{2}{\bbC}$
defined by
\begin{equation*}
\Phi_0(z,\,\lambda)
 =
 \exp(A(\lambda)\log z)P(z,\,\lambda)
\end{equation*}
satisfies
$d\Phi_0 = \Phi_0\xi$.
Since $\Phi$ satisfies the same linear ODE as $\Phi_0$,
then the map $C=\Phi\Phi_0^{-1}$ is $z$-independent,
and is an element of $\LoopSL{r}$.

(iii)
For any $k\in\{1,\dots,n-1\}$ we have
$B_0=\dots=B_{k-1}=0$, so $C_k=0$.
Since $\ELL_k$ is invertible on $\calR$,
then by the formula for $P_k$ in \autoref{thm:parametrized-zap},
$P_k = \ELL_k^{-1}(C_k) = \ELL_k^{-1}(0) = 0$.
\Depend{%
This lemma uses \autoref{def:setup}.}
\end{proof}

\typeout{=================gauge}
\subsection{Gauging the potential}
\label{sec:gauge}

\begin{xtext}
With $A$ a Delaunay residue,
let $\xi$ be a perturbed Delaunay $r$-potential
\begin{equation*}
Az^{-1}dz + \Order(z^{0})dz\spaceperiod
\end{equation*}
Let $n\in\bbZ_+$ be such that
\[
\tfrac{n}{2} < \max_{\calC_r}\Real\mu < \tfrac{n+1}{2}\spaceperiod
\]
The left inequality,
by \autoref{thm:gauge-theorem},
guarantees that $\xi$ can be gauged to the form
\begin{equation*}
A z^{-1}dz + \Order(z^{n})dz\spaceperiod
\end{equation*}
The right inequality,
by \autoref{thm:frame1},
implies the convergence of the unitary and positive
factors of $\exp(\log z A)$ and its perturbation.



\end{xtext}

\begin{xtext}
The following preliminary lemma
determines a necessary and sufficient condition on a map $X$
such that $\calL_n^{-1}(X)$ extends holomorphically
to the pole of $A$.
This condition is used to insure that the gauge constructed in
\autoref{thm:gauge-lemma} is a positive gauge, that is,
that it extends holomorphically to $\lambda=0$.
The entries of a matrix
$X\in\mattwo(\bbC)$ are denoted by $X_{ij}$ with $i,\,j\in\{1,\,2\}$.
\end{xtext}

\begin{lemma}
\label{thm:L-holo}
\theoremname{Holomorphic image of L-inverse}
Let $\calR\subset\bbC$ be a neighborhood of $p\in\bbC$.
Let $A:\calR\to\matgl{2}{\bbC}$ and $X:\calR\to\matgl{2}{\bbC}$
be meromorphic on $\calR$, with orders of entries
\begin{equation}\begin{split}
\label{eq:L-holo1}
&\ord_{p}A_{11}\ge 0\spacecomma\quad
\ord_{p}A_{12} = -1\spacecomma\quad
\ord_{p}A_{21} = 0\spacecomma\quad
\ord_{p}A_{22}\ge 0\\
&\ord_{p}X_{11}\ge 0\spacecomma\quad
\ord_{p}X_{12}\ge -1\spacecomma\quad
\ord_{p}X_{21}\ge 0\spacecomma\quad
\ord_{p}X_{22}\ge 0
\spaceperiod
\end{split}\end{equation}
Let $n\in\bbZ_+$ and let $\ELL_n$ as in \autoref{thm:L1}
defined with respect to $A$.

Then ${\ELL_n}^{-1}(X)$,
which extends meromorphically to $p$ by \autoref{thm:L1}(iv),
is holomorphic at $p$ if and only if
$A_{12}X_{21} + A_{21}X_{12}$ is holomorphic at $p$.
\Comment{ This Mathematica code
\textattachfile[%
  subject={Holomorphic Inverse of L},
  description={Mathematica code checking gauge calculations.}
  ]{attachments/gauge.m}{gauge.m}
checks this calculation.}
\end{lemma}

\begin{proof}
By~\eqref{eq:L-holo1}, we have $\ord_p\det A = -1$.
A calculation using~\eqref{eq:zap_Linv}
shows that the orders of the entries of
$Y=\ELL_n^{-1}(X)$ satisfy
\begin{equation*}
\ord_{p}Y_{11}\ge 0\spacecomma\quad
\ord_{p}Y_{12}\ge -1\spacecomma\quad
\ord_{p}Y_{21}\ge 0\spacecomma\quad
\ord_{p}Y_{22}\ge 0\spaceperiod
\end{equation*}
Hence $\ELL_n^{-1}(X)$ is holomorphic at $p$ if and only if its upper-right
entry $Y_{12}$ is holomorphic at $p$.
A calculation shows
\begin{equation*}
\ord_p{Y_{12}} = \ord_p (A_{12}X_{21}+A_{21}X_{12})
\spaceperiod
\end{equation*}
Hence $Y$ is holomorphic at $p$ if and only if $A_{12}X_{21}+A_{21}X_{12}$ is.
\end{proof}

\Comment{
The setup for \autoref{thm:gauge-lemma} is as follows:
\begin{itemize}
\item
Let $A$ be a Delaunay residue as in \autoref{def:delres}.
\item
Let $\mu$ be an eigenvalue of $A$ as in \autoref{def:delres}.
\item
Suppose
\begin{equation*}
\min_{\lambda\in\calC_r}\Real\mu(\lambda)
 \le \tfrac{n}{2} <
 \max_{\lambda\in\calC_r}\Real\mu(\lambda)
\spaceperiod
\end{equation*}
\item
Let $\Sigma^\ast\subset\bbC$ be a punctured neighborhood of $0\in\bbC$.
\item
Let $r\in(0,\,1]$.
\item
Fix $n\in\bbZ_+$.
\item
Let $\xi$ be a perturbed Delaunay $r$-potential (\autoref{def:setup})
on $\Sigma^\ast$ of the form
\begin{equation*}
\xi=Az^{-1}dz + \Order(z^{n-1})dz
\spaceperiod
\end{equation*}
\end{itemize}
Note that the gauge $g$ constructed in \autoref{thm:gauge-lemma}
is generically a map into $\LoopposSL{r}$, not into $\LooppSL{r}$,
because it is complex analytic in $z$.
}

\begin{lemma}
\label{thm:gauge-lemma}
\theoremname{Gauge}
Let $\xi$ be a perturbed Delaunay $r$-potential
\begin{equation*}
\xi=Az^{-1}dz + \Order(z^{n-1})dz
\end{equation*}
for some $n\in\bbZ_+$,
and suppose
\begin{equation}
\label{eq:gauge-inequality}
\min_{\lambda\in\calC_r}\Real\mu(\lambda)
 \le \tfrac{n}{2}
 < \max_{\lambda\in\calC_r}\Real\mu(\lambda)
\spaceperiod
\end{equation}
Then there exists
a neighborhood $\Sigma'\subset\bbC$ of $0\in\bbC$,
an $r$-gauge
$g:\Sigma'\to\LoopposSL{r}$,
and \optionalA $\gaugeConstant\in\bbC$
such that
\begin{equation}
\label{eq:zap_gauge_newxi}
\gauge{\xi}{g} = Az^{-1}dz + \gaugeConstant A z^{n-1} dz + O(z^n)dz
\spaceperiod
\end{equation}
\end{lemma}

\begin{proof}
Let $\ELL_n$ be as in \autoref{thm:L1},
defined with respect to $A$ on $\calD_r^\ast$.
The inequalities~\eqref{eq:gauge-inequality}
imply
$\calD_r^\ast\cap\{\lambda\in\bbC^\ast\suchthat\mu(\lambda)=n/2\}=\emptyset$.
\Comment{See Lemma~\ref{thm:mu-set}.}
Hence by \autoref{thm:L1}(i),
$\ELL_n$ is invertible on $\calD_r^\ast$.

Define $B$ by
\begin{equation*}
\xi=Az^{-1}dz + Bz^{n-1}dz + \Order(z^n)dz
\spaceperiod
\end{equation*}
Then
\begin{equation*}
\ord_{\lambda=0}A_{12} = -1\spacecomma\quad
\ord_{\lambda=0}A_{21} = 0\spacecomma\quad
\ord_{\lambda=0}B_{12} \ge -1\spacecomma\quad
\ord_{\lambda=0}B_{21} \ge 0\spacecomma
\end{equation*}
so
\begin{equation*}
\ord_{\lambda=0}\frac{B_{12}}{A_{12}} \ge 0\spacecomma\quad
\ord_{\lambda=0}\frac{B_{21}}{A_{21}} \ge 0\spaceperiod
\end{equation*}
Hence the limit
\begin{equation*}
\gaugeConstant := \frac{1}{2} \lim_{\lambda\to 0}
  \left(\frac{B_{12}}{A_{12}} + \frac{B_{21}}{A_{21}}\right)
\end{equation*}
exists in $\bbC$.

Let $X = \gaugeConstant A - B$.
A calculation shows that with this choice of $\kappa$,
$A_{12}X_{21} + A_{21}X_{12}$ is holomorphic at $0$.
Hence by \autoref{thm:L-holo},
$\ELL_n^{-1}(\gaugeConstant A - B)$ on $\calD_r^\ast$
extends holomorphically to $\lambda=0$.
Hence we have the holomorphic map $C:\calD_r\to\matSL{2}{\bbC}$
defined by
\begin{equation}
\label{eq:gaugeC}
C = \ELL_n^{-1}(\gaugeConstant A - B)\spaceperiod
\end{equation}

Define the holomorphic map $h:\Sigma\times\calD_r\to\matGL{2}{\bbC}$
by $h = \id + C z^{n}$.
Since $\det h$ will not in general be identically $1$,
we define $g=h+\Order(z^{n+1})$
with $\det g=1$ as follows.
%
Since $h_{22} = C_{22}z^n$, and $\ol{\calD_r}$ is bounded away from
the poles of $C_{22}$,
then we have the uniform convergence $h_{22}\to 1$ on $\ol{\calD_r}$
as $z\to 0$.
Hence there exists a neighborhood $\Sigma'\subset\Sigma$ of $0\in\Sigma$
on which $h_{22}>0$.
Define the holomorphic map
$g:\Sigma'\to\LoopposSL{r}$ by
\begin{equation*}
g= h + \begin{pmatrix}q & 0\\0 & 0\end{pmatrix}\spacecomma\quad
q = \frac{1-\det h}{h_{22}}\spaceperiod
\end{equation*}
A calculation shows $\det g=1$.

In order to show that $g$ satisfies the gauge equation~\eqref{eq:zap_gauge_newxi},
we show that $g = h + \Order(z^{n+1})$.
We have $\det h = 1 + \tr(C)z^n + \Order(z^{n+1})$.
Since $\phi A-B$ is tracefree,
then $C$ is tracefree by \autoref{thm:L1}(iii).
Hence $\det h = 1 + \Order(z^{n+1})$.
Since $h_{22}=1+\Order(z^n)$, then $q=\Order(z^{n+1})$.
Hence $g = h + \Order(z^{n+1})$.

To complete the proof, we will use the following formula.
For any $n\in\bbZ_+$, $r$-gauge $g = \id + g_n z^n + \Order(z^{n+1})$
and $r$-potential $\xi = A z^{-1}dz + B_{n-1}z^{n-1} + \Order(z^n)$,
we have
\begin{equation}
\label{eq:gauge-formula}
\gauge{\xi}{g}=A z^{-1} dz + (\ELL_n(g_n)+B_{n-1})z^{n-1}dz + \Order(z^n)dz
\spaceperiod
\end{equation}
The formula can be shown by expanding 
$\gauge{\xi}{g}=g^{-1}\xi g+g^{-1}dg$ in $z$.

The result~\eqref{eq:zap_gauge_newxi}
follows from~\eqref{eq:gaugeC},
together with formula~\eqref{eq:gauge-formula}
replacing $g_n$ with $C$ and $B_{n-1}$ with $B$.
\end{proof}


\Comment{
The setup for \autoref{thm:coordinate-change} is as follows:
\begin{itemize}
\item
Let $A$ be a Delaunay residue.
\item
Let $\Sigma\subset\bbC$ be a neighborhood of $0\in\bbC$ with
conformal coordinate $z$.
\item
Let $r\in(0,\,1]$.
\item
Let $n\in\bbZ_+$.
\item
Let $\gaugeConstant\in\bbC$.
\item
let $\xi$ be a $r$-potential on $\Sigma\setminus\{0\}$ satisfying
\begin{equation*}
\xi = Az^{-1}dz + \gaugeConstant A z^{n-1} dz + \Order(z^n)dz
\end{equation*}
\end{itemize}
}

\begin{lemma}
\label{thm:coordinate-change}
\theoremname{Coordinate change}
Let $\xi$ be a perturbed Delaunay $r$-potential (\autoref{def:setup})
\begin{equation}
\label{eq:coordinate-change1}
\xi = Az^{-1}dz + \gaugeConstant A z^{n-1} dz + \Order(z^n)dz
\end{equation}
for some $n\in\bbZ_+$ and $\gaugeConstant\in\bbC$.
Then there exists
a coordinate change $\cchange:\Sigma'\to\bbC$
on a neighborhood $\Sigma'\subset\bbC$ of $0\in\bbC$
with conformal coordinate $w$,
such that on $\Sigma'\setminus\{0\}$,
\begin{equation}
\label{eq:coordinate-change2}
\cchange^\ast\xi = A w^{-1}dw + \Order(w^n)dw\spaceperiod
\end{equation}
\end{lemma}

\begin{proof}
Let $\sigma:\Sigma\to\bbC$ be any analytic function satisfying
\begin{equation*}
\sigma(z) = z + (\gaugeConstant/n) z^{n+1} + O(z^{n+2})\spaceperiod
\end{equation*}
Then $\sigma(0)=0$ and $\sigma'(0)=1$, so
there exists a coordinate change $\cchange:\Sigma'\to\Sigma$
on a neighborhood $\Sigma'\subset\bbC$ of $0$
which is the inverse of $\sigma$.
Because $\cchange$ is a conformal diffeomorphism,
for any $j\in\bbZ$ and any differential
$\omega=\Order(z^j)dz$, we have
$\cchange^\ast \omega = \Order(w^j)dw$.
A calculation shows that on $\cchange(\Sigma')$,
\begin{equation*}
\left(z^{-1}+\gaugeConstant z^{n-1}\right)dz
 =
\left(\sigma^{-1} + \Order(\sigma^n)\right)d\sigma
\spaceperiod
\end{equation*}
\Comment{%
The calculation is as follows:
\begin{equation*}
\frac{d\sigma}{dz} = 1 + k(n+1)/n z^n + \Order(z^{n+1})
\spaceperiod
\end{equation*}
Also,
\begin{equation*}
(z^{-1}+\kappa z^{n-1})\sigma = 1 + k(n+1)/n z^n + \Order(z^{n+1})
\spaceperiod
\end{equation*}
Putting these two equations together yields the result.
}
Hence
\begin{equation*}
\cchange^\ast((z^{-1}+\gaugeConstant z^{n-1})dz) = w^{-1}dw + \Order(w^{n})dw\spaceperiod
\end{equation*}
This with~\eqref{eq:coordinate-change1}
implies the result~\eqref{eq:coordinate-change2}.
\end{proof}


\Comment{
The setup for \autoref{thm:gauge-theorem} is as follows:
\begin{itemize}
\item
Let $A$ be a Delaunay residue (\autoref{def:delres}).
\item
Let $\mu$ be its eigenvalue as in \autoref{def:delres}.
\item
Let $n\in\bbZ_+$.
\item
Let $r\in(0,\,1]$.
\item
Assume
\begin{equation*}
\min_{\lambda\in\calC_r}\Real \mu(\lambda) \le 1/2
\AND
n/2 < \max_{\lambda\in\calC_r}\Real \mu(\lambda)
\spaceperiod
\end{equation*}
\item
Let $\Sigma\subset\bbC$ be a neighborhood of $0\in\bbC$ with
conformal coordinate $z$.
\item
Let $\gaugeConstant\in\bbC$.
\item
let $\xi$ be a $r$-potential on $\Sigma\setminus\{0\}$ satisfying
\begin{equation*}
\xi=Az^{-1}dz + \Order(z^{0})dz\spaceperiod
\end{equation*}
\end{itemize}
}


We now apply \autoref{thm:gauge-lemma} and \autoref{thm:coordinate-change}
iteratively to transform $\xi$ to the form
$Aw^{-1}dw + O(w^n)dw$.

\begin{theorem}
\label{thm:gauge-theorem}
\theoremname{Gauge and coordinate change}
Let $\xi$ be a perturbed Delaunay $r$-potential (\autoref{def:setup})
\begin{equation*}
\xi=Az^{-1}dz + \Order(z^{0})dz\spaceperiod
\end{equation*}
Assume for some $n\in\bbZ_+$ that
\begin{equation}
\label{eq:gauge-inequalities}
\min_{\lambda\in\calC_r}\Real \mu(\lambda) \le \half
\AND
\tfrac{n}{2} < \max_{\lambda\in\calC_r}\Real \mu(\lambda)
\spaceperiod
\end{equation}
Then there exists a coordinate change $\cchange:\Sigma'\to\Sigma$
on a neighborhood $\Sigma'\subset\bbC$ of $0\in\bbC$
with coordinate $w$,
and an $r$-gauge $g:\Sigma'\to\LoopposSL{r}$,
such that
\begin{equation*}
\cchange^\ast(\gauge{\xi}{g}) = Aw^{-1}dw + O(w^n)dw\spaceperiod
\end{equation*}

Moreover,
let $\Phi$ satisfy $d\Phi=\Phi\xi$.
Then $\Psi = \cchange^\ast(\Phi g)$ satisfies
$d\Psi = \Psi\left( \cchange^\ast(\gauge{\xi}{g}) )\right)$,
and the monodromies of $\Phi$ and $\Psi$ on a loop around $z=0$ are equal.
\end{theorem}

\begin{proof}
For any conformal maps
 $\cchange_1:\Sigma'\to\Sigma$,
 $\cchange_2:\Sigma''\to\Sigma'$
and $r$-gauges
$g_1:\Sigma\to\LoopposSL{r}$,
$g_2:\Sigma'\to\LoopposSL{r}$,
we have the composition formula
\begin{equation}
\label{eq:gauge-composition-formula}
\cchange_2^\ast( \gauge{(\cchange_1^\ast(\gauge{\xi_1}{g_1}))}{g_2} )
=
\theta^\ast(\gauge{\xi}{g})
\spacecomma
\end{equation}
where $\cchange:\Sigma''\to\Sigma$ is the conformal map and
$g:\Sigma\to\LoopposSL{r}$ is the $r$-gauge
defined by
\begin{equation*}
\cchange = \cchange_1\circ\cchange_2
\AND
g = g_1\cdot \bigl( {\bigl(\cchange_1^{-1}\bigr)}^\ast g_2 \bigr)
\spaceperiod
\end{equation*}
The result then follows using the composition
formula~\eqref{eq:gauge-composition-formula}
by induction on $n$,
applying \autoref{thm:gauge-lemma}
and \autoref{thm:coordinate-change} repeatedly.
The conditions~\eqref{eq:gauge-inequalities} insure
that condition~\eqref{eq:gauge-inequality} holds for each step.
\end{proof}

\begin{lemma}
\label{thm:shift}
\theoremname{Shift lemma}
Let $r\in(0,\,1]$, $X\in\LoopSL{r}$, and $Y\in\LoopposSL{r}$.
Then
\begin{enumerate}
\item
There exists $U\in\matSU{2}{}$ such that
\begin{equation}
\label{eq:shift}
\Uni{r}{XY} = \Uni{r}{X}U
\AND
\Pos{r}{XY} = U^{-1}\Pos{r}{X}Y\spaceperiod
\end{equation}
\item
Let $\xi$ be an $r$-potential on a Riemann surface $\Sigma$.
Let $\Phi$ satisfy
$d\Phi= \Phi\xi$.
Let $g:\Sigma\to\LoopposSL{r}$
be a gauge.
Then $\Sym{r}{\Phi}=\Sym{r}{\Phi g}$.
\end{enumerate}
\end{lemma}

\begin{proof}
To prove (i),
let $Y_+=\Pos{r}{Y}$ and let $U_1$ be the unitary factor
in the QR-decomposition of $Y(0)$.
Then $Y = U_1 Y_+$.
Let $X_u=\Uni{r}{X}$ and $X_+=\Pos{r}{X}$, so $X=X_uX_+$.
Let $U$ be the unitary factor in the QR-decomposition of $X_+(0)U_1$.
Then $X_+U_1=U\Pos{r}{X_+U_1}$.
Hence
\[
\Pos{r}{XY} = \Pos{r}{X_+Y} = \Pos{r}{X_+U_1Y_+} = \Pos{r}{X_+U_1}Y_+
= U^{-1}X_+U_1Y_+
= U^{-1}X_+Y
\spaceperiod
\]
This proves the second equality in~\eqref{eq:shift}, and the first equality
follows from it.

To prove (ii), by (i) there exists $U\in\matSU{2}{}$ such that
$\Uni{r}{\Phi g} = \Uni{r}{\Phi}U$.
It follows that $\Sym{r}{\Phi}$ = $\Sym{r}{\Phi g}$.
\end{proof}

\typeout{=================cauchy}
\subsection{Cauchy integral formula for vector-valued maps}

\begin{xtext}
The following technical lemma
shows from the
convergence of a family of holomorphic functions
the convergence of all its derivatives in a strictly contained
subdomain.
The proof is standard and uses the Cauchy integral formula.
\end{xtext}%
We use the notation $\calR'\subset\subset \calR$ to mean that
$\calR'$ is bounded away from $\boundary{\calR}$.
Let
$X^{(n)}$ denote the $n$\!'th derivative of $X$ with respect to $\lambda$.

\begin{lemma}
\label{thm:cauchy}
\theoremname{Cauchy Limit Lemma}
Let $V$ be a finite dimensional vector space over $\bbC$,
and let $\abs{\,\cdot\,}$ be a vector norm on $V$.
With $\Sigma\subset\bbC$ a punctured neighborhood of $0\in\bbC$,
and $\calR\subset\bbC$ a bounded domain with smooth boundary,
let $X = X(z,\,\lambda):\Sigma\times \calR\to V$
be continuous in $z$ and holomorphic in $\lambda$.
Suppose $\lim_{z\to 0}\sup_{\lambda\in \calR}\abs{X}=0.$
Then for all $n\in\bbZ_{\ge 0}$,
and every subset $\calR'\subset\subset \calR$,
\begin{equation*}
\lim_{z\to 0}\sup_{\lambda\in \calR'}\abs{X^{(n)}}=0
\spaceperiod
\end{equation*}
\end{lemma}

\begin{proof}
With $n=\dim V$, we may assume that $V=\bbC^n$.
Since all vector norms on a finite dimensional vector space are equivalent,
we may assume that the norm $\abs{\,\cdot\,}$ is given by
$\abs{(Y_1,\dots,Y_n)}=\max_{j}{\abs{Y_{j}} }$.
Let $\calR'\subset\subset \calR$.
By the Cauchy integral formula, 
for all $n\in\bbZ_{\ge 0}$ and all $\lambda\in \calR'$,
\begin{equation*}
X^{(n)}(z,\,\lambda) =
 \frac{n!}{2\pi i}\int_{\boundary{\calR}} \frac{X(z,\,\nu)}{(\nu-\lambda)^{n+1}}d\nu\spacecomma
\end{equation*}
so
\begin{equation*}
{\abs{X^{(n)}(z,\,\lambda)}}
 \le
\frac{n!}{2\pi}\max_{j}\int_{\boundary{\calR}}
 \frac{\abs{X_{j}(z,\,\nu)}}{\abs{\nu-\lambda}^{n+1}}\abs{d\nu}
\le
\frac{n!}{2\pi}\int_{\boundary{\calR}}
\frac{{\abs{X(z,\,\nu)}}}{\abs{\nu-\lambda}^{n+1}}\abs{d\nu}\spaceperiod
\end{equation*}
Since $\lambda\in \calR'$ and $\calR'$ is bounded away from $\boundary{\calR}$,
then there exists \optionalA $c_1\in\bbR_{>0}$ such that
for all $\nu\in\boundary{\calR}$,
$\abs{\nu-\lambda}^{-(n+1)} < c_1$.
Hence with $c_2 = c_1 n!/(2\pi)$,
\begin{equation*}
\abs{X^{(n)}(z,\,\lambda)}
\le
c_2\int_{\boundary{\calR}}
{\abs{X(z,\,\nu)}}\abs{d\nu}\spaceperiod
\end{equation*}
Since this holds for all $\lambda\in \calR'$, and the right-hand side is
independent of $\lambda$, so
\begin{equation*}
\sup_{\lambda\in \calR'}\abs{X^{(n)}(z,\,\lambda)}
\le
c_2\int_{\boundary{\calR}}
\abs{X(z,\,\nu)}\abs{d\nu}\spaceperiod
\end{equation*}
The result follows by taking limits.
\end{proof}

\typeout{=================frame-asymptotics1}
\subsection{Asymptotics of the perturbed Delaunay frame}
\label{sec:delframe-asymptotics}

\begin{xtext}
We now have the tools to show that the $r$-Iwasawa factors
of a holomorphic frame constructed
from a perturbed Delaunay potential converge to
the respective factors of a holomorphic
frame constructed from a Delaunay potential
(\autoref{thm:frame1}).
This convergence implies that the immersions and their
metrics converge (\autoref{thm:delasym-summary-1}).
A bootstrap argument in \autoref{Sec:immersion}
strengthens this to $C^\infty$-convergence of the immersions.
\end{xtext}

\begin{lemma}
\label{thm:monodromyFB}
Let $\xi$ be a perturbed Delaunay potential, let
$d\Phi = \Phi\xi$,
and assume that the monodromy $M$ of $\Phi$ around $z=0$
satisfies $M\in\LoopuSL{r}$.
With $\Phi_0=\exp(A \log z)$,
let $\Phi=C\Phi_0 P$ be the decomposition as in
\autoref{thm:parametrized-zap}.
Let $C\Phi_0=F_0\cdot B_0$ and $\Phi=F\cdot B$ be the respective
$r$-Iwasawa factorizations on $\widetilde{\Sigma^\ast}$.
Then $F_0^{-1}F$, $B$ and $B_0$ are lifts of unique single-valued maps
on $\Sigma^\ast$.
\end{lemma}

\begin{proof}
Let $\tau:\widetilde\Sigma^\ast\to\widetilde\Sigma^\ast$ be a deck transformation
corresponding to the curve defining the monodromy $M$ of $\Phi$.
Then by the definition of $M$,
for all $z\in\widetilde\Sigma^\ast$ we have
$\tau^\ast\Phi = M\Phi$, so
\[
(\tau^\ast F) (\tau^\ast B) = M F B\spacecomma
\]
or
\begin{equation}
\label{eq:monodromyFB1}
(\tau^\ast B) B^{-1} = {(\tau^\ast F)}^{-1} M F\spaceperiod
\end{equation}
The left hand side of~\eqref{eq:monodromyFB1} is in $\LooppSL{r}$, and
since by hypothesis, $M\in\LoopuSL{r}$, the right hand side
of~\eqref{eq:monodromyFB1}
is in $\LoopuSL{r}$. Hence each side is equal to $\id$,
so $\tau^\ast B = B$, and $B$ is the lift of a single-valued
map on $\Sigma^\ast$.

Since $\Phi = C\exp(A\log z) P$, and $\tau^\ast P = P$, then
\[
M = C \exp(2\pi i A) C^{-1}\spaceperiod
\]
Since $\Phi_0 = C\exp(A\log z)$, then $M$ is also the monodromy of $\Phi_0$.
Hence $\tau^\ast\Phi_0 = M\Phi_0$. The same argument as above,
with $F_0$ and $B_0$ replacing $F$ and $B$ respectively,
shows that $\tau^\ast B_0 = B_0$, so $B_0$ is the lift of a single-valued
map on $\Sigma^\ast$.

Since
\begin{equation}
\label{eq:monodromyFB2}
F_0^{-1}F = B_0 P B^{-1}\spaceperiod
\end{equation}
and the right hand side of~\eqref{eq:monodromyFB2}
is invariant under the action of $\tau^\ast$,
then so is the left hand side of~\eqref{eq:monodromyFB2},
so $F_0^{-1}F^{-1}$ is the lift of a
single-valued map on $\Sigma^\ast$.
\end{proof}

\Comment{
The setup for \autoref{thm:frame1} is as follows:
\begin{itemize}
\item
Let $A$ be a Delaunay residue as in \autoref{def:delres}.
\item
Let $\mu$ be an eigenvalue of $A$ as in \autoref{def:delres}.
\item
Assume $\calA_{r,1}\cap\delSA=\emptyset$.
\item
Suppose
\begin{equation*}
\max_{\lambda\in\calC_r} \Real\mu(\lambda) < \tfrac{n+1}{2}
\spaceperiod
\end{equation*}
\item
Let $r\in(0,\,1)$.
\item
Let $\Sigma\subset\bbC$ be a neighborhood of $0\in\bbC$.
\item
Let $\xi$ be an $r$-potential (\autoref{def:setup}) on $\Sigma^\ast$
satisfying
\begin{equation*}
\xi=Az^{-1}dz + \Order(z^{n})dz\spaceperiod
\end{equation*}
\item
Let $\widetilde{\Sigma^\ast}\to\Sigma^\ast$ be a universal cover of
$\Sigma^\ast$.
\item
Let $\Phi:\widetilde{\Sigma^\ast}\to\LoopUpSL{r}$
satisfy $d\Phi=\Phi\xi$.
\item
Assume that the monodromy $M$ of $\Phi$ at $z=0$
satisfies $M\in\LoopuSL{r}$.
\end{itemize}
}

\begin{theorem}
\label{thm:frame1}
\theoremname{Frame asymptotics I}
Let $r\in(0,\,1)$ and assume $\calC_r\cap\delSA = \emptyset$.
Let $\xi$ be a perturbed Delaunay $r$-potential
\begin{equation}
\label{eq:delasym-xi}
\xi=Az^{-1}dz + \Order(z^{n})dz\spaceperiod
\end{equation}
Suppose
\begin{equation}
\label{eq:delasym-max}
\max_{\lambda\in\calC_r} \Real\mu(\lambda) < \tfrac{n+1}{2}
\spaceperiod
\end{equation}
Let $\Phi:\widetilde{\Sigma^\ast}\to\LoopSL{r}$
satisfy $d\Phi=\Phi\xi$
on the universal cover $\widetilde{\Sigma^\ast}\to\Sigma^\ast$
of $\Sigma^\ast$,
and assume that the monodromy $M$ of $\Phi$ around $z=0$
satisfies $M\in\LoopuSL{r}$.
Let $\Phi_0=\exp(A \log z)$
and let $\Phi = C\Phi_0 P$ be the $z^AP$-decomposition of $\Phi$
as in \autoref{thm:parametrized-zap}.
Assume for some $0<s_1<r<s_2<1$ there exists
a continuous function $c:\calA_{s_1,s_2}\to\bbR_{+}$ such that
$\Pos{r}{C\Phi_0}$ extends analytically to $\calA_{s_1,s_2}$ and
\begin{equation}
\label{eq:delasym-estimate}
\supnorm{\Pos{r}{C\Phi_0}}{\calA_{s_1,s_2}}
 \le c\abs{z}^{-\Real{\mu}}
\spaceperiod
\end{equation}
Then
\CONVERGE{\Phi}{C\Phi_0}
\begin{equation}\begin{split}
\label{eq:delasym}
&\LIMIT{ \left(\Uni{r}{C\Phi_0}\right)^{-1}\Uni{r}{\Phi} }{\id}{\calA_{\RRR}}{z\to 0}\spacecomma\\
&\LIMIT{ \Pos{r}{\Phi}\left(\Pos{r}{C\Phi_0}\right)^{-1} }{\id}{\calD_{\RRR}}{z\to 0}
\spaceperiod
\end{split}\end{equation}
\end{theorem}

\begin{proof}
\STARTSTEPS
\STEP
Let $\upsilon_1 = \max_{\lambda\in\calC_r} \Real\mu(\lambda)$.
There exist $r_1,\,r_2\in\bbR_+$ such that $s_1\le r_1 < r < r_2 \le s_2$,
$\calA_{r_1,r_2}\cap\delSA=\emptyset$ and
$\Pos{r}{C\Phi_0}$ extends analytically to $\calA:=\calA_{r_1,r_2}$.
By~\eqref{eq:delasym-estimate},
\begin{equation}
\label{eq:delasym-B}
\supnorm{\Pos{r}{C\Phi_0}}{\ol{\calA}} \le c\abs{z}^{-\upsilon_1}
\spaceperiod
\end{equation}

\STEP
Let $C\Phi_0=F_0\cdot B_0$ and $\Phi=C\Phi_0P=F\cdot B$ and be the respective
$r$-Iwasawa factorizations on $\widetilde{\Sigma^\ast}$.
Since by hypothesis $M\in\LoopuSL{r}$,
by \autoref{thm:monodromyFB},
$B$, $B_0$ and $F_0^{-1}F$
on $\widetilde{\Sigma^\ast}$ descend to single-valued
analytic maps on $\Sigma^\ast$.

Let $m=n+1$.
By \autoref{thm:parametrized-zap}(iii),
$P$ has the expansion on $\Sigma'$
\begin{equation*}
P(z,\,\lambda) = \id + \sum_{k=m}^\infty P_k(\lambda) z^k\spaceperiod
\end{equation*}
Then
\begin{equation*}
B_0 P B_0^{-1} - \id = \sum_{k=m}^{\infty} B_0 P_k B_0^{-1} z^k
\spaceperiod
\end{equation*}
By~\eqref{eq:delasym-B},
\begin{equation*}
\supnorm{B_0}{\calA} \le c \abs{z}^{-\upsilon_1}
\text{ and }
\supnorm{B_0^{-1}}{\calA} \le c \abs{z}^{-\upsilon_1}
\spacecomma
\end{equation*}
so
\begin{equation*}
\supnorm{B_0 P B_0^{-1} - \id}{\calA} \le
 \sum_{k=m}^{\infty}
 \supnorm{B_0}{\calA}
 \supnorm{P_k}{\calA}
\supnorm{B_0^{-1}}{\calA} \abs{z}^k
\le
 c^2\sum_{k=m}^\infty \supnorm{P_k}{\calA} \abs{z}^{k-2\upsilon_1}\spaceperiod
\end{equation*}
Since by~\eqref{eq:frame-asymptotics-inequality},
$\upsilon_1 < m/2$, then
$k - 2\upsilon_1 > 0$ for all $k\in\bbZ_{\ge m}$.
Hence
\begin{equation*}
\LIMIT{B_0 P B_0^{-1}}{\id}{\calA}{z\to 0}
\spaceperiod
\end{equation*}
%
By \autoref{thm:cauchy},
$B_0PB_0^{-1}\to \id$ as $z\to 0$
as a map in $C^\infty(\calC_r,\,\matSL{2}{\bbC})$.

\STEP
We have the $r$-Iwasawa factorization
\begin{equation*}
B_0 P B_0^{-1} = (F_0^{-1}F)\cdot(B B_0^{-1})\spaceperiod
\end{equation*}
%
%
As noted in \preliminarySection,
the $r$-Iwasawa factorization is a homeomorphism
in the $C^\infty$-topologies on its domain and target spaces.
The result~\eqref{eq:delasym} follows.
\Comment{The $r$-Iwasawa factorization is in fact a diffeomorphism
of Lie groups, but we do not use this here.}
\end{proof}


\section{Asymptotics of immersions}

\label{Sec:immersion}
\typeout{=================geometry}
\subsection{Convergence of geometric data}
\label{sec:geometry}

\begin{xtext}
In \autoref{thm:frame1}
we showed that the unitary and positive part of a holomorphic
frame constructed from a perturbed Delaunay potential
converge to those of a Delaunay holomorphic frame.
From this we obtain convergence
of the immersions, metrics,
and moving frames and normals.
%
This convergence is strengthened to $C^\infty$-convergence
in \autoref{sec:immersion1}.
\end{xtext}

\begin{lemma}
\label{thm:geometry}
\theoremname{Convergence of geometric data}
Let $\Sigma^\ast\subset\bbC$ be a punctured neighborhood of $0\in\bbC$
with coordinate $z$.
Let $r\in(0,\,1]$.
Let $\xi_1$ and $\xi_2$ be $r$-potentials.
Let $\alpha_k$ be the coefficient of $\lambda^{-1}$ in
the upper-right entry $\xi_k$ in its $\lambda$ expansion,
and suppose
\begin{equation}
\label{eq:geom-alpha-limit}
\lim_{z\to 0}(\alpha_1/\alpha_2) = 1
\spaceperiod
\end{equation}
For $k\in\{1,\,2\}$,
let $\Phi_k$ satisfy
$d\Phi_k = \Phi_k \xi_k$.
Let $\calR\subset\bbC$ be a neighborhood of $\lambda_0\in\bbS^1$.
Suppose
\begin{equation}\begin{split}
\label{eq:geom}
&\LIMIT{\Uni{r}{\Phi_1}^{-1}\Uni{r}{\Phi_2}}{\id}{\calR}{z\to 0}\\
&\lim_{z\to 0}
\left.\Pos{r}{\Phi_1}\left(\Pos{r}{\Phi_2}\right)^{-1}\right|_{\lambda=0}
= \id
\spaceperiod
\end{split}\end{equation}
Let $f_k = \Sym{r}{\Phi_k}$ be the immersions obtained from $\Phi_k$,
and let $v_k^2\abs{dz}^2$ be the respective metrics of $f_k$.
Then:
\begin{enumerate}
\item
The immersions converge, that is,
$\lim_{z\to 0}\left.\left(f_2-f_1\right)\right|_{\lambda_0} = 0$.
\item
The metrics converge in ratio; that is, $\lim_{z\to 0}(v_2/v_1) = 1$.
\item
The moving frames $G_k$ for $f_k$ converge; that is,
$\lim_{z\to 0}\left.(G_1^{-1}G_2)\right|_{\lambda_0} = \id$.
\item
The difference of the normals of $f_1$ and $f_2$ converges to $0$
as $z\to 0$.
\end{enumerate}

Moreover, let $x+iy=\log z$
and let $w_k^2(dx^2+dy^2)$ be the metric for $f_k$ in the
$x$ and $y$ coordinates.
If $w_1$ is bounded and bounded away from $0$ on $\Sigma$,
then
the immersions have $C^1$-convergence in the coordinates $x$ and $y$;
that is,
\begin{equation}
\label{eq:geometry-df}
\lim_{z\to 0}\left.\left(\tfrac{d}{dx}f_2-\tfrac{d}{dx}f_1\right)\right|_{\lambda_0} = 0
\AND
\lim_{z\to 0}\left.\left(\tfrac{d}{dy}f_2-\tfrac{d}{dy}f_1\right)\right|_{\lambda_0} = 0
\spaceperiod
\end{equation}
\end{lemma}

\begin{proof}
\STARTSTEPS
For $k\in\{1,\,2\}$, let
 $F_k=\Uni{r}{\Phi_k}$ and
 $B_k=\Pos{r}{\Phi_k}$.
Let prime denote differentiation with respect to $\theta$, where
$\lambda=e^{i\theta}$.
%

\STEP
To show the convergence of the immersions, 
by the Sym formula for $\bbR^3$ in \preliminarySection,
we have
\begin{equation*}
f_k=-2H^{-1}F_k'F_k^{-1}\spacecomma\quad k\in\{1,\,2\}
\spacecomma
\end{equation*}
where $H\in\bbR_+$ is the constant mean curvature.
Hence
\begin{equation}
\label{eq:end13}
f_2 - f_1
 =
 -2H^{-1} F_1 \bigl(F_1^{-1} F_2\bigr)' F_2^{-1}\spaceperiod
\end{equation}
Let $\calR'\subset\subset\calR$ be a neighborhood of $\lambda_0$.
Since $F_1^{-1}F_2-\id$ is holomorphic on $\calR$,
by~\autoref{thm:cauchy},
\begin{equation}
\label{eq:F1invF2}
\lim_{z\to 0}\supnorm{(F_1^{-1}F_2)'}{\calR'} = 0
\spaceperiod
\end{equation}
The convergence of the immersions
$\lim_{z\to 0}\left.(f_2-f_1)\right|_{\lambda_0}=0$
follows from~\eqref{eq:end13} and~\eqref{eq:F1invF2}
and the fact that
$\matnorm{F_1(\lambda_0)}= 1 = \matnorm{F_2^{-1}(\lambda_0)}$.

\STEP
To show the convergence of the metrics, since $F_kB_k=\Phi_k$,
then
$B_k$ satisfies the gauge equation $\gauge{(F_k^{-1}dF_k)}{B_k} = \xi_k$.
An examination of the coefficient of $\lambda^{-1}$ in this gauge equation
above shows that
\begin{equation}
\label{eq:B-metric-formula}
v_k = 2\abs{H^{-1}\alpha_k}\rho_k^2
\spacecomma
\end{equation}
\Comment{See \autoref{thm:dpw-computations}.}
where $\rho_k$ is
the constant term of the upper-right entry of
$B_k$ in its $\lambda$ expansion.
Equation~\eqref{eq:geom-alpha-limit} implies
$\lim_{z\to 0}(\rho_1/\rho_2) = 1$,
and the result follows.

\STEP
The moving frames $G_k$ for $f_k$ are defined by the equations
\begin{equation}
\label{eq:geom-moving-frames}
(f_k)_x = v_k G_k \left(\begin{smallmatrix}0 & dz\\-d\ol{z} & 0\end{smallmatrix}\right) G_k^{-1}\spacecomma
\quad
k\in\{1,\,2\}
\spaceperiod
\end{equation}
The extended frames $F_k$ and the moving frames $G_k$ are related
by a gauge $g_k$ defined by
\[
G_k = F_kg_k\spacecomma\quad
g_k = \diag(p_k,\,p_k^{-1})\spacecomma\quad
p_k^2 = i H \alpha_k/(\lambda\abs{H\alpha_k})
\quad
k\in\{1,\,2\}
\spaceperiod
\]
\Comment{See \autoref{thm:dpw-computations}.}
Then
\begin{equation}
\label{eq:geom-moving-frame-ratio}
G_1^{-1}G_2 - \id = g_1^{-1}(F_1^{-1}F_2 - g_1 g_2^{-1})g_2
\spaceperiod
\end{equation}
By~\eqref{eq:geom-alpha-limit}, it follows
that $\lim_{z\to 0}g_1(\lambda_0)g_2^{-1}(\lambda_0)=1$.
Since $\matnorm{g_k(\lambda_0)}$ is finite,
the convergence of the moving frames (iii)
follows from~\eqref{eq:geom-moving-frame-ratio}.
The convergence of the normals follows from the
convergence 
of the moving frames.
\Comment{See \autoref{thm:normals}.}

\STEP
To show the convergence of the derivatives of the
immersions~\eqref{eq:geometry-df}, let
\[
e_1=(\begin{smallmatrix}0 & 1\\ -1& 0\end{smallmatrix})
\AND
e_2=(\begin{smallmatrix}0 & i\\ i& 0\end{smallmatrix})\spaceperiod
\]
Writing $x=x_1$ and $y=x_2$,
by~\eqref{eq:geom-moving-frames} we have
$(f_k)_{x_j} = w_k G_k e_j G_k^{-1}$.
Hence
\begin{equation*}
\tfrac{d}{dx_j}(f_2-f_1)= (w_2-w_1)G_2e_jG_2^{-1} + w_1 G_1 [G_1^{-1}G_2,\,e_j]G_2^{-1}
\spaceperiod
\end{equation*}
Since $\lim_{z\to 0}(v_2/v_1)=1$
and $v_k=w_k e^{-x/2}$, then
$\lim_{z\to 0}(w_2/w_1)=1$.
Since by assumption, $w_1$ or $w_2$ is bounded away from $0$,
then $\lim_{z\to 0}(w_2-w_1)=0$.
The first equation of~\eqref{eq:geometry-df}
follows, since
$\matnorm{G_1(\lambda_0)e_jG_1^{-1}(\lambda_0)}$ is finite,
$w_1$ is bounded,
$\matnorm{G_1(\lambda_0)}$ and
$\matnorm{G_2^{-1}(\lambda_0)}$ are finite,
and $\lim_{z\to 0}\matnorm{G_1^{-1}(\lambda_0)G_2(\lambda_0)}=0$ by~\eqref{eq:F1invF2}.
\end{proof}


\begin{theorem}
\label{thm:delasym-summary-1}
\theoremname{Frame and metric convergence}
Let $A$ be a Delaunay residue satisfying
\begin{equation}
\label{eq:min-eigenvalue-condition}
\min_{\lambda\in\calC_r} \Real\mu(\lambda) \le \half
\spaceperiod
\end{equation}
Let $\xi$ be a Delaunay $r$-perturbation of $A$.
Let $\Phi$ satisfy $d\Phi=\Phi\xi$
with unitary monodromy at $z=0$.
Assume $\Phi$ is in $\LoopUpSL{r}$,
and $(\calC_r\cup\calA_{r,\,1})\cap\delSA=\emptyset$.
Let $f$ be the CMC immersion induced by $\Phi$.
Then there exists a Delaunay immersion $f_0$
with the same necksize as that induced by $\Phi_0:=\exp(A\log z)$,
such that as $z\to 0$,
in the coordinates $x+iy=\log z$,
the ratio of their metrics converges to $0$,
and $f-f_0$ converges $C^1$ to $0$.
\end{theorem}

\begin{proof}
We shall require the following fact:
as functions of $r$,
$\max_{\lambda\in\calC_r}\Real\mu$ is strictly decreasing on $(0,\,1]$,
and $\min_{\lambda\in\calC_r}\Real\mu$ is $0$ on $(0,\,\abs{\nu_1}]$
and is strictly increasing on $[\abs{\nu_1},\,1]$.
(The boundary of the upper graph in \autoref{fig:graph}
illustrates this behavior.)

Let
 $\upsilon_0 = \min_{\lambda\in\bbC_r}\Real\mu(\lambda)$ and
 $\upsilon_1 = \max_{\lambda\in\bbC_r}\Real\mu(\lambda)$.
Since $\upsilon_1\ge\upsilon_0\ge 1/2$, there exists \optionalA
$n\in\bbZ_+$ such that
\begin{equation*}
\tfrac{n}{2} \le \upsilon_1 < \tfrac{n+1}{2}
\spaceperiod
\end{equation*}
If $\upsilon_1 = n/2$,
by the above fact,
$r$ can be replaced with a smaller value
and $\upsilon_0$ and $\upsilon_1$ by the corresponding values
so that this new $\nu_1$ satisfies
\begin{equation}
\label{eq:frame-asymptotics-inequality}
\tfrac{n}{2} < \upsilon_1 < \tfrac{n+1}{2}
\spaceperiod
\end{equation}
Then $\upsilon_0$ is still $\le 1/2$ by the above fact.
By \autoref{thm:gauge-theorem}
and \autoref{thm:shift}(ii),
we may assume
after a gauge and coordinate change that $\xi$ has the
form~\eqref{eq:delasym-xi}.

%
%
Since $\Phi$ is in $\LoopUpSL{r}$,
we have the $z^AP$ decomposition $\Phi = C \exp(A \log z) P$
on $\calA_{r,1}\setminus\delSA$.
Since by assumption
$(\calC_r\cup\calA_{r,\,1})\cap\delSA=\emptyset$,
then $C$ and $P$ extend to $\calC_r\cup\calA_{r,1}$.

By \autoref{thm:delgrowth1},
there exists a continuous function $c:\calD_1^\ast\to\bbR_+$
such that for all $(z,\,\lambda)\in\{0<\abs{z}<1\}\times\calD_1^\ast$,
\begin{equation*}
\label{eq:frame-asymptotics-B-estimate}
\matnorm{\Pos{r}{C\Phi_0}}
 \le c\abs{z}^{-\Real\mu}\spaceperiod
\end{equation*}

%
By \autoref{thm:frame1},
\CONVERGE{\Phi}{C\Phi_0}~\eqref{eq:delasym}.
Hence taking $\Phi_1=\Phi_0$ and $\Phi_2=\Phi$ in \autoref{thm:geometry},
condition~\eqref{eq:geom} holds.
By \autoref{thm:triv-del-dress},
$C\Phi_0$ induces a Delaunay immersion
which differs from $\Phi_0$ by a rigid motion.
Taking $\xi_1=A\,dz/z$ and $\xi_2=\xi$ in \autoref{thm:geometry},
condition~\eqref{eq:geom-alpha-limit} of \autoref{thm:geometry} holds
because $\xi_2$ is a holomorphic perturbation of $\xi_1$.

With $w_1$ and $w_2$ as in \autoref{thm:geometry},
a computation shows that with $v$ as in \autoref{thm:delframe},
by~\eqref{eq:B-metric-formula},
the metric of the Delaunay immersion is
$w_1 = 4\abs{abH^{-1}}v^{-1}$.
Hence $w_1$ is periodic and nonzero in the coordinate $x=\log\abs{z}$,
and is hence bounded and bounded away from $0$ on $\Sigma$.
The assertion follows by \autoref{thm:geometry}.
\end{proof}

\typeout{=================immersion-asymptotics}
\subsection{Asymptotics of conformal immersions}
\label{sec:immersion1}

\begin{xtext}
Now that we have shown that the frames and metrics of the
immersion produced by a perturbed Delaunay potential
converge to those of a Delaunay immersion,
we apply a bootstrap argument on the Gauss equation
to obtain $C^\infty$-convergence of the immersions
up to rigid motion.
For details and formulas, see~\cite{Bobenko_CMC_1991}.
\end{xtext}
The notation $f\converge{C^n}g$ means
$\lim_{z\to 0}\supnorm{f-g}{C^n}=0$,
where the $C^n$ norm is taken over some implied compact domain.


\begin{theorem}
\label{thm:conformal-immersion-asymptotics}
\theoremname{Bootstrap}
Let $\Omega\subset\bbR^2$ be a bounded domain,
and let $\Omega'\subset\subset\Omega$ be a strictly contained subset.
Let $f:\ol{\Omega}\to\bbR^3$ and the sequence
$f_j:\ol{\Omega}\to\bbR^3$, $j\in\bbZ_+$
be a sequence of
conformal immersions in $C^\infty(\Omega)\cap C^0(\ol{\Omega})$.
Let $e^u,\,H,\,Q\,dz^2$ and $e^u_j,\,H_j,\,Q_j dz^2$
be the respective conformal factors, mean curvature functions,
and Hopf differential
of $f$ and the $f_j$.
Suppose that
$f_j\converge{C^1} f$,
$H_j\converge{C^0} H$,
$Q_j\converge{C^{\infty}} Q$
and $u_j\converge{C^0} u$ on $\Omega$.
Then $f_j\converge{C^\infty} f$ on $\Omega'$.
\end{theorem}

\begin{proof}
\STARTSTEPS
We will use the following
gradient estimate for Poisson's equation~\cite{Gilbarg_Trudinger_1977}.
Let $\Sigma$ be a bounded domain in $\bbR^2$ and
$\Sigma'\subset\subset\Sigma$.
Then there exists \optionalA $c\in\bbR_+$ depending only on
$\dist(\Sigma',\,\boundary{\Sigma})$
such that for all $w\in C^2(\Sigma)\cap C^0(\ol{\Sigma})$,
\begin{equation}
\label{eq:poisson}
\sup_{\Sigma'}\abs{\del_x w}
+\sup_{\Sigma'}\abs{\del_y w}
 \le
 c\bigl(\sup_{\Sigma}\abs{w}+\sup_{\Sigma}\abs{\Delta w}\bigr)\spaceperiod
\end{equation}
Fix $k\in\bbN_{\ge 1}$ and choose domains
$\Omega_1,\dots,\Omega_k$ with
$\Omega'=\Omega_k\subset\subset \dots \subset\subset\Omega_1\subset\subset\Omega$.

\STEP
We apply the estimate~\eqref{eq:poisson} $k$ times
to show that
$u_j\converge{C^k}u$ and $H_j\converge{C^{k-1}}H$ on $\Omega'$

The first iteration is as follows.
For the functions $u,\,H,\,Q$, define
\begin{equation*}
\Psi[u,\,H,\,Q] = -\half H^2 e^w + 2 \abs{Q}^2 e^{-w}\spaceperiod
\end{equation*}
Then the Gauss equations
for $f$ and the $f_j$ are respectively
\begin{equation*}
\Delta u = \Psi[u,\,H,\,Q]\spacecomma\quad 
\Delta u_j = \Psi[u_j,\,H_j,\,Q_j]\spaceperiod
\end{equation*}
Hence
\begin{equation*}
\Delta(u_j-u) = \Phi[u_j,\,H_j,\,Q_j] - \Phi[u,\,H,\,Q]\spaceperiod
\end{equation*}
Since
$u_j\converge{C_0}u$,
$H_j\converge{C_0}H$ and
$Q_j\converge{C_0}Q$ on $\Omega$
then 
\begin{equation*}
\Phi[u_j,\,H_j,\,Q_j] - \Phi[u,\,H,\,Q]\ \converge{C_0}\ 0\spacecomma
\end{equation*}
so
$\Delta(u_j-u)\converge{C_0}0$  on $\Omega$.
Applying the estimate~\eqref{eq:poisson}
with $\Sigma=\Omega$, $\Sigma'=\Omega_1$,
and $w=u_j-u$, we obtain
$u_j\converge{C_1}u$ on $\Omega_1$.

The Codazzi equations
for $f$ and the $f_j$ are respectively
\begin{equation*}
H_z = 2 e^{-u}Q_{\ol{z}}\spacecomma\quad
(H_j)_z = 2 e^{-u_j}(Q_j)_{\ol{z}}\spaceperiod
\end{equation*}
Since $u_j\converge{C^0} u$ and $Q_j\converge{C^1} Q$ on $\Omega$,
it follows that $H_j\converge{C^1} H$ on $\Omega$.

Using the Gauss equation again in the same way, since $u_j\converge{C^1}u$,
$H_j\converge{C_1}H$ and $Q\converge{C_1}Q$ on $\Omega_1$,
then $\Delta(u_j-u)\converge{C_1}0$.
Hence
$\Delta \del_x(u_j-u)\converge{C_0}0$
and
$\Delta \del_y(u_j-u)\converge{C_0}0$.
Applying the estimate~\eqref{eq:poisson} again,
with $\Sigma=\Omega_1$, $\Sigma'=\Omega_2$,
and $w$ take to be first $\del_x(u_j-u)$ and then $\del_y(u_j-u)$,
we obtain
$u_j\converge{C_2}u$ on $\Omega_2$.

Repeating the argument $k-2$ more times shows that
$u_j\converge{C_k}u$
and $H_j\converge{C^{k-1}}H$ on $\Omega_k=\Omega'$.
Since $k$ is arbitrary, the above argument shows
$u_j\converge{C^\infty}u$
and $H_j\converge{C^\infty}H$
on $\Omega'$.

\STEP
Let $e_1,\,e_2,\,e_3$
be a positively oriented orthonormal basis for $\matsu{2}{}$.
Let $F$ and $F_j$ be the respective moving frames for $f$ and $f_j$
with values in $\matSU{2}{}$, defined respectively by
\begin{equation}
\label{eq:imm-f}
\del_x f = e^u Fe_1F^{-1}\spacecomma\quad
\del_y f = e^u Fe_2F^{-1}\spacecomma
\AND
\del_x f_j = e^u F_je_1F_j^{-1}\spacecomma\quad
\del_y f_j = e^u F_je_2F_j^{-1}\spaceperiod
\end{equation}
By assumption, $f_j\converge{C^1} f$, from which it follows that
$F_j\converge{C^0}F$.

Write the Lax pairs for $F$ and $F_j$ respectively
\begin{equation*}
\del_x F = F U\spacecomma\quad \del_y F = F V\spacecomma
\AND
\del_x F_j = F_j U_j\spacecomma\quad \del_y F_j = F_j V_j\spaceperiod
\end{equation*}
Since
 $u_j\converge{C^\infty}u$,
 $H_j\converge{C^\infty}H$ and
 $Q_j\converge{C^\infty}Q$,
then
$U_j\converge{C^\infty}U$ and
$V_j\converge{C^\infty}V$.
Since the Lax pair equations are linear and $F_k\converge{C^0}F$,
it follows that $F_j\converge{C^\infty} F$
by smooth dependence of solutions on parameters.
By~\eqref{eq:imm-f},
we obtain the convergence $f_j\converge{C^{\infty}} f$.
\end{proof}

\subsection{Embedded ends}

\begin{xtext}
Under appropriate assumptions,
if an immersion $f$ of a cylinder converges $C^2$ to an embedding of the
cylinder, then $f$ is embedded.
\end{xtext}

Let $\Sigma_{\rm Cyl}$ be the cylinder%
\Comment{This theorem is written in terms of $x>0$, not $x<0$, because
it's more natural. However, throughout the rest of the paper we have
convergence $x\to -\infty$ (since we want $x+iy=\log z$, and $x\to -\infty$
as $z\to 0$.}
\begin{equation}
\label{eq:cylinder}
\Sigma_{\rm Cyl}
 =
 \{(x,\,y)\in\bbR^2\suchthat x> 0 \text{ and }  0\le y\le 2\pi \}\spacecomma
\end{equation}
with the boundary lines $y=0$ and $y=2\pi$ identified.

\begin{lemma}
\label{thm:embedded}
\theoremname{Embeddedness}
With
$\Sigma_{\rm Cyl}$ as in~\eqref{eq:cylinder},
let $f:\Sigma_{\rm{Cyl}}\to\bbR^3$
be a properly embedded surface with bounded curvature
which satisfies the following condition:
\begin{equation}
\label{eq:spread}
\parbox{5in}{
For all $\delta>0$ there exists \optionalAn $\epsilon>0$ such that
for all $p,\,q\in \Sigma_{\rm{Cyl}}$,\\
if $\abs{p-q}_{\Sigma_{\rm{Cyl}}}>\delta$, then 
$\abs{f(p)-f(q)}_{\bbR^3}>\epsilon$.}
\end{equation}
Let $g:\Sigma_{\rm{Cyl}}\to\bbR^3$ be an immersion which converges 
to $f$ as $x\to+\infty$, in the $C^2$-topology of $\Sigma_{\rm{Cyl}}$.
Then there exists \optionalAn $x_0\in\bbR_+$ such that $g$ restricted to
$\{x>x_0\}$ is properly embedded.
\end{lemma}

\begin{proof}
The $C^0$-convergence of $g$ to $f$ and the properness of $f$ implies
that $f$ is proper.

Since $f$ has bounded curvature, there exists
$\delta>0$ such that for any two points $p,\,q\in\bbR^2$,
if $\abs{p-q}_{\Sigma_{\rm{Cyl}}}<\delta$, 
then the total curvature of the image under $f$ 
of the shortest straight line from $p$ to $q$ 
in $\Sigma_{\rm{Cyl}}$ is at most $\pi/4$.
By the $C^2$-convergence 
of $g$ to $f$, we can choose \optionalAn $x_1>0$ so that the total curvature 
of the image under $g$ of those straight lines is less 
than $\pi/2$ when 
the $x$ coordinates of $p$ and $q$ are greater than $x_1$.  

With the above choice of $\delta$,
let $\epsilon$ be as in~\eqref{eq:spread}.
By the $C^0$-convergence 
of $g$ to $f$, we can choose an $x_2$ so that $\abs{f-g}_{\bbR^3}<\epsilon/4$ 
at any point for which $x>x_2$.  Set $x_0=\max\{x_1,x_2\}$.  

Now suppose that $g$ is not embedded on $\Sigma_{\rm{Cyl}} 
\cap \{x>x_0\}$, so there exist
distinct $p,\,q\in \Sigma_{\rm{Cyl}} \cap \{x>x_0\}$ such 
that $g(p)=g(q)$.
Then the proof divides into two cases.

\textit{Case 1: $\abs{p-q}_{\Sigma_{\rm{Cyl}}}\ge\delta$.}
Then $\abs{f(p)-f(q)}_{\bbR^3}>\epsilon$.
Thus $\abs{g(p)-g(q)}_{\bbR^3}>\epsilon/2$,
which contradicts $g(p)=g(q)$.

\textit{Case 2: $\abs{p-q}_{\Sigma_{\rm{Cyl}}}<\delta$.}
Let $\gamma:[0,\,1]\to \Sigma_{\rm{Cyl}}$ be the parametrized straight line from
$p$ to $q$.
Then $g\circ\gamma$ is a closed curve in $\bbR^3$ which is smooth
except at one point.
If $g\circ\gamma$ is smoothed out at that point, it has total curvature
$\ge 2\pi$.%
\Comment{This is probably due to Milnor.}
But the turning angle of the tangent at the non-smooth point is $\le\pi$, 
hence the total curvature of the smoothed out 
$g\circ\gamma$ is $\leq 3\pi/2$, a 
contradiction.
\end{proof}

\typeout{=================main1}
\subsection{Asymptotics of Delaunay immersions}

\begin{xtext}
We conclude Part~1 with its culminating theorem,
showing that a CMC end obtained from a perturbed Delaunay
potential is asymptotic to a half-Delaunay surface.
This theorem restricts itself to the case that
the holomorphic frame extends analytically to $\calA_{r,1}$.
\end{xtext}

\begin{theorem}
\label{thm:main1}
\theoremname{Main theorem I}
Let $A$ be a Delaunay residue as in~\eqref{eq:delres}
satisfying condition~\eqref{eq:min-eigenvalue-condition}.
Let $r\in(0,\,1]$ and assume $\calA_{r,\,1}\cap\delSA=\emptyset$.
On the punctured unit disk $\Sigma^\ast=\{z\in\bbC\suchthat 0<\abs{z}<1\}$,
let
\begin{equation*}
\xi = A\frac{dz}{z} + \Order(z^0)dz
\end{equation*}
be a perturbed Delaunay $r$-potential.
Let $f_0$ and $f$ be the immersions of $\Sigma^\ast$
induced by the generalized Weierstrass representation at
$\lambda=1$
by $A\,dz/z$ and $\xi$ respectively,
so $f_0$ is a Delaunay immersion.
Assume that $f$ is obtained from a holomorphic $r$-frame
with values in $\LoopUpSL{r}$ whose monodromy at $z=0$
is in $\LoopuSL{r}$.

With $\Sigma_{\rm Cyl}$ as in ~\eqref{eq:cylinder},
let $\phi:\Sigma_{\rm Cyl}\to\Sigma^\ast$ be the map $\phi(x,\,y)=e^{x+iy}$.
Then some rigid motion of $f\circ\phi$ converges to $f_0\circ\phi$
in the $C^\infty$-topology of $\Sigma_{\rm{Cyl}}$ as $x\to -\infty$.
Furthermore, if $f_0$ is embedded, then $f$ is properly embedded.
\end{theorem}

\begin{proof}
Let $\rho$ be the period of the Delaunay surface $f_0$.
Let
$\Sigma_{\rm{Rect}}=[0,\,\rho]\times[0,\,2\pi]\subset\bbR^2$.
Fix a small $\epsilon>0$.
Then $f\circ\phi$ and $f_0\circ\phi$ can be lifted to
the $\epsilon$-neighborhood $\calN_\epsilon(\Sigma_{\rm{Rect}})$
of $\Sigma_{\rm{Rect}}$.
We define a sequence $g_k$ of $C^\infty$ functions on $\Omega$
by
\begin{equation*}
g_k(x,\,y) = (f\circ\phi)((x+k\rho)+i y)
\end{equation*}
on $\Sigma_{\rm{Rect}}$ and
extend $g_k$ real analytically to
$\calN_\epsilon(\Sigma_{\rm{Rect}})$.

Let $u,\,H,\,Q$ and $u_k,\,H_k,\,Q_k$ be the geometric
data for $f_0$ and the $g_k$ respectively.
%
Let $\Phi_0$ and $\Phi$ be the respective solutions to
$d\Phi_0 = \Phi A dz/z$ and $d\Phi = \Phi\xi$ which respectively
induce $f_0$ and $f$ via the generalized Weierstrass representation.
By \autoref{thm:delasym-summary-1},
the metrics are $v=e^u$ and $v_k=e^{u_k}$.
we have $C^1$-convergence of $g_k$ to some rigid motion of $f_0\circ\phi$
on $\calN_\epsilon(\Sigma_{\rm Rect})$ as $x\to\infty$,
and $v/v_k\converge{C^0}1$ as $x\to\infty$.
Hence $u_k\converge{C^0}u$
on
$\calN_\epsilon(\Sigma_{\rm Rect})$.
%
We have that $Q_k\converge{C^\infty}Q$, and $H_k=H$ for all $k$,
so \autoref{thm:conformal-immersion-asymptotics}
and \autoref{thm:embedded} imply the theorem.
\end{proof}

\begin{remark}
\label{thm:exponential-convergence}
\theoremname{Exponential convergence}
When $f_0$ is embedded, the rigid motion of $f$ converges exponentially
to $f_0$ in the sense of~\cite{Korevaar_Kusner_Solomon_1989}, as proven there.
\end{remark}

\Depend{
\begin{proof}
\end{proof}
}

\begin{remark}
The above Delaunay asymptotics result for Euclidean space $\bbR^3$
most likely also holds for the spaceform $\bbS^3$,
with appropriate changes in the Sym and moving frame formulas.
For hyperbolic space $\bbH^3$, the value $\lambda_0$ in the Sym formula is
in the interior of the unit disk, so if the asymptotics result is
to hold, its proof would require the results of Part~2.
\end{remark}

\part{Dressed Delaunay asymptotics}
\label{part2}

\begin{xtext}
Part~2 generalizes the results in Part~1
to the setting of $\LoopSL{r}$ for arbitrary $r$.
It shows for this larger class of frames
that an immersion constructed from a perturbed
Delaunay $r$-potential via $r$-dressing is asymptotic
to a Delaunay surface.
\end{xtext}

\typeout{=================outline2}
\section*{Outline of results}
\label{sec:outline2}

\setcounter{COUNTER}{1}

\begin{xtext}

Dressing by simple factors performs a Bianchi-B\"acklund
transformation on a CMC surface. The formula for dressing by simple
factors is known explicitly~\cite{Terng_Uhlenbeck_2000}.
In \autoref{sec:delaunay-dressing},
we show that a dressed Delaunay frame
$C\exp(A\log z)$ which has
unitary monodromy is a multibubbleton frame,
that is, a Delaunay frame dressed by a finite product of simple factors
(\autoref{thm:bub}).
This result is used in \autoref{sec:delframe-asymptotics2}.

In \autoref{sec:simple-factor-asymptotics},
we show that
dressing by a finite product $G$ of simple factors with
distinct singularities preserves asymptotics
(\autoref{thm:bridge}).

In \autoref{Sec:delaunay-dressing-asymptotics},
we show that
if an immersion $f$ converges to a Delaunay immersion $f_0$,
then $f$ dressed by a product of simple factors with distinct singularities
also converges to a Delaunay immersion with the same necksize as $f_0$
(\autoref{thm:simple-factor-asymptotics-theorem}).
Together with the methods of \autoref{Sec:immersion},
this implies that dressing by finitely many simple factors
preserves convergence to a Delaunay surface
(\autoref{thm:shimpei}.)

Given a holomorphic perturbation
\[
\xi = A\frac{dz}{z} + \Order(z^{0}) dz
\]
of the Delaunay potential $A\,dz/z$
which represents a closed once-wrapped Delaunay surface,
let $\Phi$ satisfy $d\Phi = \Phi\xi$
with unitary monodromy at $z=0$.
In \autoref{sec:delframe-asymptotics2}
we show that the unitary factor of $\Phi$
is asymptotic to a Delaunay frame
and the metric of the surface induced by $\Phi$ is asymptotic
to the corresponding Delaunay metric
(\autoref{thm:frame2}).

These asymptotic frames results are applied
to construct CMC surfaces in $\bbR^3$.
The surface induced by $\Phi$
is asymptotic to a half-Delaunay surface
induced by $A\,dz/z$
(\autoref{thm:main2}).
\end{xtext}

\section{Simple factor dressing}
\label{Sec:simple-factor-dressing}

\typeout{=================simplefactor}
\subsection{Simple factors}
\label{sec:simple-factors}

\begin{xtext}
We recall the
notion of a \emph{simple factor}~\cite{Terng_Uhlenbeck_2000}.
\end{xtext}
\begin{definition}
\label{def:simple-factor}
\theoremname{Simple factors}
Given $r\in(0,\,1]$ and $\lambda_0\in\calA_{r,1}$,
let $f$ be the unique rational map on $\CPone$ with degree one
whose pole is at $\lambda_0$, satisfying $f^\ast = f^{-1}$ and $f(1)=1$.
(Here, $f^{-1}$ denotes the multiplicative inverse).
\Comment{See \autoref{thm:f}}
%
Let $L\in\CPone$.
An \emph{unnormalized simple factor}
$\simplefactorUN{\lambda_0}{L}\in\LoopposSL{r}$ is a loop
\begin{equation*}
\simplefactorUN{\lambda_0}{L} = f^{1/2}\pi_L + f^{-1/2}\pi_{L^\perp}\spacecomma
\end{equation*}
where $\pi_L$ denotes the orthogonal projection to $L$.
%
A \emph{normalized simple factor}
$\simplefactor{\lambda_0}{L}\in\LooppSL{r}$ is a map of the form
$\simplefactor{\lambda_0}{L} = U^{-1}\simplefactorUN{\lambda_0}{L}$, where
$U$ is the unitary factor of the QR-decomposition
of $\left.\simplefactorUN{\lambda_0}{L}\right|_{\lambda=0}$.
%
A \emph{general simple factor} is a map
of the form $Ug\in\LoopposSL{r}$,
where $U\in\matSU{2}{}$ and $g$ is a normalized simple
factor.
\end{definition}

By Proposition 4.2 in \cite{Terng_Uhlenbeck_2000}, dressing by 
simple factors is explicit:
given $r\in(0,\,1)$, an extended $r$-unitary frame $F(z,\,\lambda)$,
and a normalized simple factor $\simplefactor{\lambda_0}{L}$
with $\lambda_0\in\calA_{r,1}$,
we have the formula~\cite{Terng_Uhlenbeck_2000}
(see also~\cite{Kilian_Schmitt_Sterling_2004, Kobayashi_2004})
\begin{equation}
\label{eq:simpledress}
  \Uni{r}{\simplefactor{\lambda_0}{L}\,F} = 
  \simplefactor{\lambda_0}{L}\, F \, 
  \simplefactor{\lambda_0}{\overline{F(z,\lambda_0)}^t\,L}^{-1}\spaceperiod
\end{equation}
\Comment{See \autoref{thm:backlund} for a proof of formula~\ref{eq:simpledress}.}

\begin{xtext}
While simple factors are positive $r$-loops,
the product of two simple factors with the same singularity
extends to a meromorphic map on $\CPone$.
%
\end{xtext}

\begin{lemma}
\label{thm:simple-factor-conj}
\theoremname{Simple factor conjugation}
For $k\in\{1,\,2\}$,
let $g_k=U_k^{-1}\simplefactorUN{\lambda_0}{L_k}$
be general simple factors, where
$\lambda_0\in\calD_1^\ast$,
$L_k\in\CPone$, and $U_k\in\matSU{2}{}$.
%
Let $X:\calR\to\mattwo{\bbC}$ be a meromorphic map
on a domain $\calR\subset\bbC$.
Then $g_1 X g_2^{-1}$ extends meromorphically to a map $Y$ on $\calR$.
Moreover, suppose $\calR$ is invariant under the map
$\lambda\to 1/\ol{\lambda}$. Then
\begin{enumerate}
\item
If $g_1=g_2$ and $X^\ast = X$, then
$Y^\ast = Y$ away from its poles.
\item
If $X$ is invertible and $X^\ast = X^{-1}$, then
$Y^\ast = Y^{-1}$ away from its poles.
\end{enumerate}
\end{lemma}

\begin{proof}
By the definition of simple factors (\autoref{def:simple-factor}),
for $k\in\{1,\,2\}$
we can write $g_k = U_k f^{-1/2}h_k$, where
$U_k\in\matSU{2}{}$,
$h_k=f\pi_{L_k} + \pi_{L_k^\perp}$,
$f$ has a simple pole at $\lambda_0$, $f^\ast = f^{-1}$,
and $h_k^\ast = h_k^{-1}$.
Then
$g_1 X {g_2}^{-1} = U_1h_1Xh_2^{-1}U_2^{-1}$
extends meromorphically to $\calR$ because all its components do.
Statements (i) and (ii) follow from the symmetries
$f^\ast=f^{-1}$, $h_k^\ast = h_k^{-1}$ and $U_k^\ast = U_k^{-1}$.
\end{proof}

\typeout{=================bubbleton-dressing}
\subsection{The dressed Delaunay frame}
\label{sec:delaunay-dressing}

\begin{xtext}
We shall prove that if a 
dressing matrix preserves the unitarity of a Delaunay monodromy, then the 
dressing matrix is a product of simple factors and a loop $V$ which
conjugates the Delaunay residue to a Delaunay residue.
The effect of $V$ is a coordinate 
change and a rigid motion of the surface. Thus, for Delaunay surfaces, the 
class of dressing matrices that preserve unitarity of monodromy coincide with 
the Bianchi-B\"acklund transformations.

\end{xtext}

\begin{lemma}
\label{thm:CAC}
\theoremname{Rational residue}
Let $A$ be a Delaunay residue and let $M=\exp(2\pi i A)$.
With $r\in(0,\,1)$,
let $\CPLUS\in\LooppSL{r}$
and suppose $\CPLUS M\CPLUS^{-1}\in\LoopuSL{r}$.
Then $\CPLUS A\CPLUS^{-1}$ extends to a rational map on $\CPone$
all of whose poles are simple and lie in
$(\delSA\cap(\calA_{\RRR}\setminus\bbS^1))\cup\{0,\infty\}$.
\end{lemma}

\begin{proof}
Let $\mu$ be an eigenvalue of $A$ as in \autoref{def:delres}.
Let $x,\,y:\bbC^\ast\to\bbC$ be the holomorphic functions
$x=\cos(2\pi\mu)$ and $y=\mu^{-1}\sin(2\pi\mu)$.
Note that $x$ and $y$ are always holomorphic on $\bbC^\ast$,
even if $\mu$ is not.
Then $M = x \id + i y A$.
\Comment{By \autoref{lem:exp-xy}.}
Writing $A_1=\CPLUS A\CPLUS^{-1}$,
on $\calD_{\RRR}^\ast$ we have
\begin{equation*}
\CPLUS M\CPLUS^{-1} = x \id + i y A_1\spaceperiod
\end{equation*}
Since by hypothesis $\CPLUS M\CPLUS^{-1}$
extends holomorphically to $\bbC^\ast$,
and since $x\id$ is holomorphic on $\bbC^\ast$, then
$yA_1$ extends holomorphically to $\bbC^\ast$.
Since $y$ is holomorphic on $\bbC^\ast$, then
$A_1$ extends meromorphically to $\bbC^\ast$, and
for all $\lambda\in\bbC^\ast$
\begin{equation}
\label{eq:ord-CAC}
\ord_\lambda A_1 \ge -\ord_\lambda y\spaceperiod
\end{equation}

We compute the orders of the poles of $A_1$
on $\bbC^\ast$.
These occur only at the zeros of $y$,
all of which are in $\delSA$.
Let $\lambda_0\in\bbC^\ast$ be a zero of $y$.
Then $\mu(\lambda_0) \in \half\bbZ^\ast$.
\Comment{By \autoref{lem:exp-xy}.}

\STARTSTEPS
\CASE
If $\lambda_0=0$, then $\CPLUS$ is holomorphic and $A$ has
a simple pole, so $A_1$ is meromorphic
with at worst a simple pole at $0$.

\CASE
If $0<\abs{\lambda_0}\le r$, then
$A_1$ is holomorphic at $\lambda_0$ because
$\CPLUS$ and $A$ are.


\CASE
$r<\abs{\lambda_0}<1$.
It can be shown that
$\{\lambda\in\bbC^\ast\suchthat \mu'(\lambda)=0\}\subset\bbS^1$.
Hence $\mu'(\lambda_0)\ne 0$,
\Comment{See \autoref{lem:mu-facts}{ii} for a proof.}
from which it follows that
$\ord_{\lambda_0}y=1$.
By~\eqref{eq:ord-CAC}
$\ord_{\lambda_0} A_1 \ge -1$,
so $A_1$ has at most a simple pole at $\lambda_0$.

\CASE
$\abs{\lambda_0}=1$.
In this case, $A_1$ is holomorphic at $\lambda_0$
as in the proof of \autoref{thm:triv-del-dress}.

Since $\exp(2\pi i A_1)$
takes values in $\matSU{2}{}$ on $\bbS^1$,
and $A_1$ is tracefree,
then $A_1$ has the hermitian symmetry
${A_1}^\ast=A_1$.
Hence the poles of $A_1$ on $\bbC^\ast$ lie in
$(\calA_r\setminus\bbS^1)\cap\delSA$.
\end{proof}

\begin{theorem}
\label{thm:bub}
\theoremname{Delaunay dressing}
Let $A$ be a Delaunay residue and let $M=\exp(2\pi i A)$.
With $r\in(0,\,1)$,
let $\CPLUS\in\LooppSL{r}$
and suppose $\CPLUS M\CPLUS^{-1}\in\LoopuSL{r}$.
Then there exists a loop $G\in\LooppSL{\RRR}$ which is
a product of normalized simple factors (or $G=\id$),
the singularities of the simple factors
are distinct and lie in $\delSA\cap(\calA_{\RRR}\setminus\bbS^1)$,
and, with $V=G^{-1}C_+$,
$A_1 = \VVV A\VVV^{-1}$ is a Delaunay residue.
\end{theorem}

\begin{proof}
The proof is by induction on the number of poles of $\CPLUS A\CPLUS^{-1}$
in $\calA_{r,\,1}$.
We prove the following induction step, which constructs a single
simple factor $g$.

\begin{nestedtheorem}
With $A$, $M$, $r$ and $\CPLUS$ as in the statement of the theorem,
there exists a normalized simple factor $g\in\LooppSL{r}$
such that, with $\CPLUSP = g^{-1}\CPLUS\in\LooppSL{r}$,
we have $\CPLUSP M{\CPLUSP}^{-1}\in\LoopuSL{r}$,
and $\CPLUSP A{\CPLUSP}^{-1}$ has one fewer pole
in $\calD_1$ than $\CPLUS A\CPLUS^{-1}$.
\end{nestedtheorem}

Proof of the induction step:

By \autoref{thm:CAC},
$\CPLUS A\CPLUS^{-1}$ extends to a rational map on $\CPone$
all of whose poles are simple and lie in
$(\delSA\cap(\calA_{\RRR}\setminus\bbS^1))\cup\{0,\infty\}$.
Let
\begin{equation*}
M_1 = \CPLUS M\CPLUS^{-1}
\AND
A_1 = \CPLUS A\CPLUS^{-1}
\spaceperiod
\end{equation*}

Let $x=\cos(2\pi\mu)$ and $y=\mu^{-1}\sin(2\pi\mu)$
be the holomorphic functions on $\bbC^\ast$ as in the proof of
\autoref{thm:CAC},
so that
$M = x \id + i y A$.
\Comment{By \autoref{lem:exp-xy}.}
Then
\begin{equation}
\label{eq:expCMC-long}
\CPLUS M\CPLUS^{-1} = x \id + i y A_1\spaceperiod
\end{equation}

Choose a pole $\lambda_0\in \delSA\cap(\calA_{\RRR}\setminus\bbS^1)$
of $A_1$,
so $\ord_{\lambda_0}A_1=-1$.
Since $\lambda_0\in\delSA$,
then $\mu(\lambda_0)\in\half\bbZ^\ast$,
so $\iota:=x(\lambda_0)\in\{1,\,-1\}$.
As in case 3 of the proof of \autoref{thm:CAC}, $\ord_{\lambda_0}y=1$.
Hence $\ord_{\lambda_0}yA_1=0$.
Then $M_1(\lambda_0)\ne\iota\id$, by \eqref{eq:expCMC-long}.
Hence $M_1(\lambda_0)$ is not semisimple,
so the eigenspace of $M_1(\lambda_0)$ associated to the 
eigenvalue $\iota$ is one-dimensional.
Let $v\in\bbC^2\setminus\{0\}$ be an eigenline of $M_1(\lambda_0)$
and let $L=[v]\in\CPone$.
Define the unnormalized simple factor
$\gfact=\simplefactorUN{\lambda_0}{L}$
on $\calD_{\RRR}$.

We show that $\gfact^{-1} A_1 \gfact$
is holomorphic at $\lambda_0$.
Let $f$ be as in \autoref{def:simple-factor} for $\gfact$,
so $\gfact = f^{1/2}\pi_{L} + f^{-1/2} \pi_{L^\perp}$
and $\ord_{\lambda_0}f = -1$.
Working in the basis $\{v,\,v^\perp\}$,
define the functions $a,\,b,\,c,\,d$ by
$A_1v=av+cv^\perp$ and $A_1v^\perp=b v+ d v^\perp$.
Since $\Span_{\bbC}\{v\}$ is the kernel and
image of $\left.yA_1\right|_{\lambda_0}$, then
\begin{equation*}
\ord_{\lambda_0}a \ge 0\spacecomma\quad
\ord_{\lambda_0}b = -1\spacecomma\quad
\ord_{\lambda_0}c \ge 0\spacecomma\quad
\ord_{\lambda_0}d \ge 0\spaceperiod
\end{equation*}
Then $\det A_1 = ad-bc$, so $\ord_{\lambda_0}(ad-bc) \ge 0$.
Hence $\ord_{\lambda_0}c\ge 1$.
On the other hand,
\begin{equation*}
h^{-1}A_1hv = av+fcv^\perp
\AND
h^{-1}A_1hv^\perp = f^{-1}b v + d v^\perp\spaceperiod
\end{equation*}
Hence $h^{-1}A_1hv$ and $h^{-1}A_1hv^\perp$ are holomorphic at $\lambda_0$,
so $h^{-1}A_1h$ is holomorphic at $\lambda_0$.

Let $U\in\matSU{2}{}$ be the unitary factor
in the QR-decomposition of $\gfact(0)$,
and define $g = \gfact U^{-1}$. Then $g\in\LooppSL{r}$
is a normalized simple factor.
Define $\CPLUSP = g^{-1}\CPLUS$.
Then $\CPLUS$ and $\CPLUSP$ have the same poles
on $\CPone\setminus\{\lambda_0,\,1/\ol{\lambda_0}\}$.
Since $\CPLUS A{\CPLUS}^{-1}$ has a pole at $\lambda_0$, and
as shown above, ${\CPLUSP}A{\CPLUSP}^{-1}$ does not, then
${\CPLUSP}A{\CPLUSP}^{-1}$ has one fewer pole that
${\CPLUS}A{\CPLUS}^{-1}$ in $\calD_1$.

Since $A_1$ is hermitian away from its poles,
then so is $\gfact^{-1}A_1\gfact$ by \autoref{thm:simple-factor-conj}(i).
By conjugating~\eqref{eq:expCMC-long} by $\gfact^{-1}$, we obtain that
$\gfact^{-1}M_1\gfact\in\LoopuSL{r}$.
Then
$\CPLUSP M{\CPLUSP}^{-1} = U\gfact^{-1}M_1\gfact U^{-1}\in\LoopuSL{r}$.
This proves the induction step.

To prove the theorem,
let $n$ be the number of poles of $\CPLUS A\CPLUS^{-1}$ in $\calA_{r,1}$.
If $n=0$,
then by \autoref{thm:triv-del-dress},
$\CPLUS A\CPLUS^{-1}$ is a Delaunay residue,
and the theorem follows with $G=\id$ and $V=\CPLUS$.
Otherwise, repeated application of the induction step
produces $n$ normalized simple factors $g_1,\dots,g_n\in\LooppSL{r}$
with distinct singularities
in $\delSA\cap(\calA_{\RRR}\setminus\bbS^1)$,
such that, with $G=g_1\cdots g_n$
and $V=G^{-1}\CPLUS\in\LooppSL{r}$,
$VAV^{-1}$ has no poles in $\calA_{r}$,
and hence in $\calD_1$,
and $VMV^{-1}\in\LoopuSL{r}$.
By \autoref{thm:triv-del-dress},
$VAV^{-1}$ is a Delaunay residue.
This proves the theorem.
\end{proof}

\typeout{=================delaunay-growth2}
\subsection{Dressed Delaunay growth}
\label{sec:dressed-delaunay-growth}

\begin{xtext}
\autoref{thm:delgrowth} showed that under the assumption of
unitary monodromy,
dressing a Delaunay holomorphic frame $\exp(A\log z)$
by a loop $C\in\LoopSL{r}$ does not affect the growth behavior
of its positive factor,
provided that $C\in\LoopUpSL{r}$.
\autoref{thm:delgrowth2} generalizes this result by removing this provision.
\end{xtext}

\begin{lemma}
\label{thm:simple-factor-bound}
\theoremname{Simple factor boundedness}
Let $h(z) = U(z)^{-1}\simplefactorUN{\lambda_0}{L(z)}$
be a map from a domain $\Sigma\subset\bbC$
into the space of general simple factors.
Then for every set $\calR\subset\CPone$ bounded away from
$\{\lambda_0,\,1/\ol{\lambda_0}\}$,
$\supnorm{h(z)}{\calR}$ is bounded as a function of $z$.
\end{lemma}

\begin{proof}
With $f$ and $U$ as in \autoref{def:simple-factor},
write
\begin{equation*}
h = U^{-1}
\left({f(\lambda)}^{1/2}\pi_{L(z)}
 +
 {f(\lambda)}^{-1/2}\pi_{L(z)^\perp}\right)
\spaceperiod
\end{equation*}
Then $\abs{h} = \max\{\abs{f^{1/2}},\,\abs{f^{-1/2}}\}$.
\Comment{by \autoref{thm:proj-norm},}
On any set $\calR\subset\CPone$ bounded away from
$\{\lambda_0,\,1/\ol{\lambda_0}\}$,
the function $\abs{f}$ is bounded away from $0$ and $\infty$,
and so $\abs{f^{1/2}}$ and $\abs{f^{-1/2}}$ are bounded on $\calR$.
Hence for all $z\in\Sigma$,
\begin{equation*}
\supnorm{h(z)}{\calR}
 \le \max\{
 {\lvert\!\lvert f^{1/2}\rvert\!\rvert}_{\calR},\ %
 {\lvert\!\lvert f^{-1/2}\rvert\!\rvert}_{\calR}
 \}
\spaceperiod
\qedhere 
\end{equation*}
\end{proof}

\begin{theorem}
\label{thm:delgrowth2}
\theoremname{Dressed Delaunay growth II}
Let $A$ be a Delaunay residue,
and let $\mu$ be its eigenvalue as in \autoref{def:delres}.
Let $\tau$ be as in \autoref{thm:tau}.
Let $r\in(0,\,1)$,
and assume
$\calC_r\cap \delSA=\emptyset$, where $\delSA$ is defined by~\eqref{eq:S_A}.
Let $C\in\LoopSL{r}$ and assume the monodromy
$C\exp(2\pi i A)C^{-1}$ is in $\LoopuSL{r}$.
Then there exists $r_1,\,r_2$ with $0<r_1<r<r_2<1$, and
a continuous function $c:\calA_{r_1,r_2}\to\bbR_+$,
such that for all $(z,\,\lambda)\in\{0<\abs{z}<1\}\times\calA_{r_1,r_2}$,
\begin{equation*}
\matnorm{\Pos{r}{C\exp(A\log z)}}
 \le c\abs{z}^{-\tau}
 \le c\abs{z}^{-\Real\mu}\spaceperiod
\end{equation*}
\end{theorem}

\begin{proof}

\STARTSTEPS
\STEP
Let $\Sigma^\ast = \{0<\abs{z}<1\}$.
Let $C = C_u\cdot C_+$ be the $r$-Iwasawa factorization of $C$.
By \autoref{thm:bub},
there exist loops $G,\,V\in\LooppSL{r}$ such that
$C_+ = GV$,
$G$ is a product of normalized simple factors (or $G=\id$),
$A_1:=VAV^{-1}$ is a Delaunay residue,
and
\begin{equation*}
X:= C\exp(A \log z) = C_uG(\exp A_1\log z)V
\spaceperiod
\end{equation*}

\STEP
Since $\calC_r\cap \delSA=\emptyset$, and the points of $\delSA$
are isolated on $\bbC^\ast$, there exist $r_1$ and $r_2'$ satisfying
$0<r_1<r<r_2'<1$ such that
$\calA_{r_1,r_2}\cap \delSA=\emptyset$,
and $V$ extends holomorphically to $\calD_{r_2'}$.
Let $r_2\in(r,\,r_1)$ and let $\calA:=\calA_{r_1,r_2}$.
%
Since $V$ extends holomorphically to $\calD_{r_2'}$,
then $c_1 := \supnorm{V}{\calA}$ is finite
on $\Sigma^\ast\times\calA$,
so
\begin{equation}
\label{eq:delgrowth-dressed0}
\matnorm{\Pos{r}{X}} = c_1\matnorm{\Pos{r}{G\exp(A_1\log z)}}
\spaceperiod
\end{equation}

\STEP
%
Write $G=g_n\cdots g_1\cdot\id$
as a product of normalized simple factors
and for $k\in\{1,\,\dots,n\}$, write
$g_k = \simplefactor{\lambda_k}{L_k}$,
with $\lambda_k\in\calA_{r,1}$ and $L_k\in\CPone$.
By repeated use of the simple factor dressing formula~\eqref{eq:simpledress},
we obtain
maps $h_1,\dots,h_n:\{0<\abs{z}<1\}\to\LooppSL{r}$ defined by
\begin{equation*}
h_k = \simplefactor{\lambda_k}{\transpose{\ol{F_k(\lambda_k)}}L_k},
\quad
F_k = \Uni{r}{g_{k-1}\dots g_1 X}
\end{equation*}
such that
\begin{equation*}
\Pos{r}{g_n\dots g_1 X} = h_n\cdots h_1\Pos{r}{X}\spaceperiod
\end{equation*}
Since the singularities of $G$ are in $\delSA$, which is bounded
away from $\calA$,
then by \autoref{thm:simple-factor-bound},
$\supnorm{h_k}{\calA}$ is finite
for each $k\in\{1,\dots,n\}$.
so let $c_2$ be their product.
Then using~\eqref{eq:delgrowth-dressed0},
we have on $\Sigma^\ast\times\calA$,
\begin{equation}
\label{eq:delgrowth-dressed1}
\matnorm{\Pos{r}{X}} \le c_1 c_2 \matnorm{\Pos{r}{\exp(A_1\log z)}}\spaceperiod
\end{equation}

\STEP
It follows from $\det A = \det A_1$ that
the function $\tau$ in \autoref{thm:tau} for $A$
is the same as that for $A_1$.
By \autoref{thm:delgrowth},
there exists a continuous function $c_3:\calA\to\bbR_+$
such that
on $\{0<\abs{z}<1\}\times\calA$,
\begin{equation}
\label{eq:delgrowth-dressed2}
\matnorm{\Pos{r}{\exp(A_1\log z)}}
\le
 c_3\abs{z}^{-\tau}
 \le
 c_3\abs{z}^{-\Real\mu}.
\end{equation}
The result follows from~\eqref{eq:delgrowth-dressed1}
and~\eqref{eq:delgrowth-dressed2} with $c=c_1c_2c_3$.
\end{proof}

\section{Delaunay dressing asymptotics}

\typeout{=================simplefactor-asymptotics}
\subsection{Simple factor asymptotics}
\label{sec:simple-factor-asymptotics}

\begin{xtext}
The next preliminary lemma shows a basic
convergence property of simple factors:
two simple factors with the same fixed singularities converge
if the lines defining them converge.
\Comment{%
This result is used in \autoref{thm:bubbleton-asymptotics}
and \autoref{thm:bridge}.
The theorem is true for normalized or unnormalized simple factors,
but false for general simple factors.

The Riemann sphere $\CPone\cong\bbS^2\subset\bbR^3$
has the subspace topology of $\bbS^2\subset\bbR^3$
and is compact.
}
\end{xtext}


\begin{lemma}
\label{thm:simple-factor-limit}
\theoremname{Simple factor limit}
For $k\in\{1,\,2\}$, let
$h_k = \simplefactor{\lambda_0}{L_k}$,
where $\lambda_0\in\calD_1^\ast$
and $L_k:\Sigma^\ast\to\CPone$ are continuous maps
on a punctured neighborhood $\Sigma^\ast\subset\bbC$ of $0\in\bbC$.
Note that the map $h_1 h_2^{-1}$ extends meromorphically
to $\CPone$
(see \autoref{thm:simple-factor-conj})
and this extension is holomorphic on
$\CPone\setminus\{\lambda_0,\,1/\ol{\lambda_0}\}$.
If $L_1$ and $L_2$ converge as $z\to 0$,
then on every region $\calR\subseteq\CPone$ bounded away from
$\{\lambda_0,\,1/\ol{\lambda_0}\}$,
\begin{equation}
\label{eq:bub21}
\LIMIT{h_1 h_2^{-1}}{\id}{\calR}{z\to 0}
\spaceperiod
\end{equation}
\Comment{This convergence is
in the standard topology of $\CPone\cong\bbS^2\subset\bbR^3$.
See \autoref{def:line-convergence}.}

\end{lemma}

\begin{proof}
Let $\calR\subset\CPone$
be a region bounded away from $\{\lambda_0,\,1/\ol{\lambda_0}\}$.

\STARTSTEPS
\STEP
For $k\in\{1,\,2\}$, let $h_k = U_k^{-1}\simplefactorUN{\lambda_0}{L_k}$.
Since $h_1$ and $h_2$ share the same singularity $\lambda_0$,
we have $f$ as in \autoref{def:simple-factor} such that
\begin{equation*}
h_k^\circ = f^{1/2}\pi_{L_k} + f^{-1/2}\pi_{L_k^\perp},
\quad
k\in\{1,\,2\}\spaceperiod
\end{equation*}
We show that $h_1^\circ {h_2^\circ}^{-1}\to \id$.
We have
\begin{equation}
\label{eq:bub22}
h_1^\circ  {h_2^\circ}^{-1} - \id
 = \left(\pi_{L_1}\pi_{L_2} + \pi_{L_1^\perp}\pi_{L_2^\perp} - \id\right)
   + f \pi_{L_1}\pi_{L_2^\perp}
   + f^{-1} \pi_{L_1^\perp}\pi_{L_2}\spaceperiod
\end{equation}
Then
\begin{equation}
\label{eq:proj-limit-3}
\lim_{z\to 0}
 \bigl(\pi_{L_1}\pi_{L_2} + \pi_{L_1^\perp}\pi_{L_2^\perp} -\id\bigr)
 = 0\spacecomma
\quad
\lim_{z\to 0}
 \pi_{L_1}\pi_{L_2^\perp} = 0\spacecomma\\
\quad
\lim_{z\to 0}
 \pi_{L_1^\perp}\pi_{L_2} = 0\spaceperiod
\end{equation}
\Comment{See \autoref{thm:proj-limit}.}
Moreover, since $f$ and $f^{-1}$ are holomorphic on $\calR$, which is bounded
away from $\{\lambda_0,\,1/\ol{\lambda_0}\}$.
then $f$ and $f^{-1}$ are bounded on $\calR$.
Then \eqref{eq:bub22} and~\eqref{eq:proj-limit-3}
imply
\begin{equation}
\label{eq:bub23}
\lim_{z\to 0}\supnorm{h_1^\circ {h_2^\circ}^{-1}- \id}{\calR} = 0\spaceperiod
\end{equation}

\STEP
We show that $h_1 h_2^{-1}\to \id$.
We have
\begin{equation*}
h_1 h_2^{-1} - \id = U_1^{-1}h_1^\circ{h_2^\circ}^{-1}U_2 - \id
=
U_1^{-1}(h_1^\circ{h_2^\circ}^{-1}- U_1U_2^{-1})U_2\spacecomma
\end{equation*}
so
\begin{equation}
\label{eq:bub29}
\supnorm{h_1h_2^{-1} - \id}{\calR}
\le \supnorm{h_1^\circ{h_2^\circ}^{-1}- U_1U_2^{-1}}{\calR}
\le \supnorm{h_1^\circ{h_2^\circ}^{-1}-\id}{\calR}
 + \supnorm{U_1U_2^{-1}-\id}\spaceperiod
\end{equation}
By~\eqref{eq:bub23}, the first term
of the right-hand side of~\eqref{eq:bub29}
converges to $0$ as $z\to 0$.

We have $U_1U_2^{-1}=h_1^\circ(0) {h_2^\circ(0)}^{-1}$.
Since $\{0\}$ is bounded away from $\{\lambda_0,\,1/\ol{\lambda_0}\}$,
then by step 1, \eqref{eq:bub23} holds at $0$. Hence
\begin{equation*}
\lim_{z\to 0}
\matnorm{h_1^\circ(x,\,0) {h_2^\circ(x,\,0)}^{-1}-\id} = 0\spaceperiod
\end{equation*}
Hence
the second term
of the right-hand side of~\eqref{eq:bub29}
converges to $0$ as $z\to 0$.
This implies the result~\eqref{eq:bub21}.
\end{proof}



\begin{xtext}
The following theorem shows that
dressing by a simple factor preserves convergence.
For a meromorphic map $U$,
the set of poles of $U$ is denoted by $\sing{U}$.
\end{xtext}

\begin{theorem}
\label{thm:bridge}
\theoremname{Bridge theorem}
With $\RRR\in(0,\,1)$,
let $\Phi_1,\,\Phi_2:\widetilde{\Sigma^\ast}\to\LoopSL{\RRR}$ be analytic maps
on the universal cover
$\widetilde{\Sigma^\ast}\to\Sigma^\ast$
of a punctured neighborhood $\Sigma^\ast\subset\bbC$ of $0\in\bbC$.
Let $U\in\LoopuMeroSL{\RRR}$ and $V\in\LooppSL{\RRR}$.
Suppose that for every region $\calA\subseteq\calA_{\RRR}$ bounded away
from $\sing{U}$,
\CONVERGE{\Phi_1}{\Phi_2}
\begin{subequations}
\label{eq:bridge1}
\begin{align}
\label{eq:bridge1-F}
&\LIMIT{ {\left(\Uni{\RRR}{\Phi_1}\right)}^{-1} U \Uni{\RRR}{\Phi_2}}{\id}{\calA}{z\to 0}\spacecomma\\
\label{eq:bridge1-B}
&\LIMIT{ \Pos{\RRR}{\Phi_2} V^{-1} {\left(\Pos{\RRR}{\Phi_1}\right)}^{-1}}{\id}{\calD_{\RRR}}{z\to 0}
\spaceperiod
\end{align}
\end{subequations}


For $k\in\{1,\,2\}$,
let $g_k=W\simplefactor{\lambda_0}{L_k}$
be general simple factors,
where
$W\in\matSU{2}{}$,
$\lambda_0\in\calA_{\RRR,1}\setminus\sing{U}$,
and $L_1,\,L_2\in\CPone$ are related by
$L_2 = {\transpose{\ol{U(\lambda_0)}}}L_1$.
Then for every region $\calA'\subseteq\calA_{\RRR}$ bounded away
from $\sing{U}\cup\{\lambda_0,\,1/\ol{\lambda_0}\}$,
\CONVERGE{g_1\Phi_1}{g_2\Phi_2}
\begin{equation}\begin{split}
\label{eq:bridge2}
&\LIMIT{ {\left(\Uni{\RRR}{g_1\Phi_1}\right)}^{-1} (g_1Ug_2^{-1}) \Uni{\RRR}{g_2\Phi_2}}{\id}{\calA'}{z\to 0}\spacecomma\\
&\LIMIT{ \Pos{\RRR}{g_2\Phi_2} V^{-1} {\left(\Pos{\RRR}{g_1\Phi_1}\right)}^{-1}}{\id}{\calD_{\RRR}}{z\to 0}
\spaceperiod
\end{split}\end{equation}
\end{theorem}

\begin{proof}
For $k\in\{1,\,2\}$, let $F_k = \Uni{\RRR}{\Phi_k}$
and $B_k=\Pos{\RRR}{\Phi_k}$.
For $k\in\{1,\,2\}$
we have by the simple factor formula~\eqref{eq:simpledress}
\begin{equation*}
\Uni{\RRR}{g_k\Phi_k} = g_kF_kh_k^{-1}
\AND
\Pos{\RRR}{g_k\Phi_k} = h_k B_k,\quad k\in\{1,\,2\}
\spacecomma
\end{equation*}
where the $h_k:\Sigma^\ast\to\LooppSL{\RRR}$ are defined by
\begin{equation*}
h_k = \simplefactor{\lambda_0}{\transpose{\ol{F_k(\lambda_0)}}L_k},\quad
k\in\{1,\,2\}.
\end{equation*}

Let $\calA'\subset\calA_{\RRR}$ be a region bounded away from
$\sing{U}\cup\{\lambda_0,\,1/\ol{\lambda_0}\}$.
Since $\calA'$ is bounded away from $\sing{U}$,
then \eqref{eq:bridge1} hold for $\calA'$:
\begin{equation}\begin{split}
\label{eq:bridge3}
&\LIMIT{F_1^{-1}UF_2}{\id}{\calA'}{z\to 0}\spacecomma\\
&\LIMIT{B_2 V^{-1} B_1^{-1}}{\id}{\calD_{\RRR}}{z\to 0}\spaceperiod
\end{split}\end{equation}
By \autoref{thm:simple-factor-conj},
$h_1h_2^{-1}$
extends to a map $\Sigma^\ast\times\calA'_{\RRR}\to\matSL{2}{\bbC}$
which is meromorphic in the second variable,
and $h_2h_1^{-1}$
extends to a holomorphic map $\Sigma^\ast\to\LooppSL{\RRR}$.
For $k\in\{1,\,2\}$, define $M_k=\transpose{\ol{F_k(\lambda_0)}}$,
and define $U_0=\transpose{\ol{U(\lambda_0)}}$.
Since $\lambda_0\not\in\sing{U}$,
by~\eqref{eq:bridge1-F} we have
$\lim_{z\to 0}M_2U_0M_1^{-1}=\id$.
Then
\begin{equation*}
M_2 L_2 = (M_2 U_0 M_1^{-1})M_1(U_0^{-1}L_2) = (M_2 U_0 M_1^{-1})M_1L_1
\spacecomma
\end{equation*}
so $M_1L_1$ and $M_2L_2$ converge as $z\to 0$.
Since $\calA'$ is bounded away from $\{\lambda_0,\,1/\ol{\lambda_0}\}$,
we have by \autoref{thm:simple-factor-limit},
\begin{equation}\begin{split}
\label{eq:bridge4}
&\LIMIT{h_1h_2^{-1}}{\id}{\calA'}{z\to 0}\spacecomma\\
&\LIMIT{h_2h_1^{-1}}{\id}{\calD_{\RRR}}{z\to 0}
\spaceperiod
\end{split}\end{equation}

Since $\calA'$ and $\calD_r$ are
bounded away from $\{\lambda_0,\,1/\ol{\lambda_0}\}$,
then by \autoref{thm:simple-factor-bound},
each of $h_1$, $h_1^{-1}$, $h_2$ and $h_2^{-1}$ is bounded
on $\Sigma^\ast\times\calA'$ and on $\Sigma^\ast\times\calD_r$.
By~\eqref{eq:bridge3} and~\eqref{eq:bridge4} it follows that
\begin{align*}
&\LIMITZ{ h_1(F_1^{-1}UF_2-\id)h_2^{-1} + (h_1h_2^{-1}-\id) }{0}{\calA'}{z\to 0}\spacecomma\\
&\LIMITZ{ h_2(B_2V^{-1}B_1^{-1}-\id)h_1^{-1} + (h_2h_1^{-1}-\id) }{0}{\calD_{\RRR}}{z\to 0}
\spaceperiod
\end{align*}
This gives the result~\eqref{eq:bridge2}.
%
\Comment{
That is,
\begin{equation}\begin{split}
&{\left(\Uni{\RRR}{g_1\Phi_1}\right)}^{-1} (gUg^{-1}) \Uni{\RRR}{g_2\Phi_2}
 = {(gF_1h_1^{-1})}^{-1} (gUg^{-1}) (gF_2h_2^{-1})\\
&\Pos{\RRR}{g_2\Phi_2} V {\left(\Pos{\RRR}{g_1\Phi_1}\right)}^{-1}
 = (h_2 B_2) V {(h_1 B_1)}^{-1}\spaceperiod
\end{split}\end{equation}
}
\end{proof}

\label{Sec:delaunay-dressing-asymptotics}
\typeout{=================bubbleton-asymptotics}
\subsection{Bubbleton asymptotics}

\begin{xtext}
In this section, we show that
if an immersion $f$ converges to a Delaunay immersion $f_0$,
then $f$ dressed by a product of simple factors with distinct singularities
converges to a Delaunay immersion with the same necksize as $f_0$
(\autoref{thm:simple-factor-asymptotics-theorem}).
The proof of \autoref{thm:simple-factor-asymptotics-theorem}
is outlined as follows.

For each singularity $\lambda_0$ there is a family of simple factors
sharing this singularity.
Given a Delaunay residue $A$, dressing the Delaunay frame
$F = \Uni{}{\exp( (x+iy) A)}$ by this family
generically produces a bubbleton frame $\dress{}{g}{F}$.
However, for one special member $g_0$ of this family of simple
factors, $\dress{}{g_0}{F}$ is not a bubbleton frame as expected, but rather
a Delaunay frame (\autoref{thm:special-del-dressing}).

On the other hand, any two Delaunay frames
dressed by simple factors
are asymptotic (\autoref{thm:bubbleton-asymptotics}).
In particular, a generic bubbleton frame $\dress{}{g}{F}$
constructed by dressing a Delaunay frame $F$ by a simple factor $g$,
is asymptotic to the corresponding special Delaunay frame $\dress{}{g_0}{F}$.
This implies that any bubbleton frame is asymptotic
to a Delaunay frame.

By the general result in \autoref{sec:simple-factor-asymptotics},
that dressing by simple factors preserves end convergence
(\autoref{thm:bridge}),
step 3 above implies that dressing by a finite product of simple
factors with distinct singularities preserves Delaunay
asymptotics (\autoref{thm:simple-factor-asymptotics-theorem}).
\end{xtext}


\begin{xtext}
The following lemma shows that
given a unitary Delaunay frame $F$ and
simple factors $g_1$ and $g_2$ with the same singularity,
then generically, $\dress{}{g_1}{F}$ and $\dress{}{g_2}{F}$
are asymptotic modulo a unitary factor.
A similar result holds for the positive factors.
\end{xtext}

\begin{lemma}
\label{thm:exp-limit}
\theoremname{exp limit}
Let $A\in\matsl{2}{\bbC}$, let $\mu$ be an eigenvalue of $A$
and suppose $\mu\in\bbR_+$.
Let $(\mu,\,L_+),\,(-\mu,\,L_-)\in\bbR^\ast\times\CPone$
be the eigenvalue-eigenline pairs for $A$.
Then for all $L\in\CPone\setminus\{L_+\}$,
$\exp(x A)L\to L_-$ as $x\to-\infty$.
\end{lemma}

\Comment{
\begin{proof}
\textit{Proof 1.}
We have
\begin{align*}
\exp(x A) &= \cosh(x\mu)\id + \mu^{-1}\sinh(x\mu)A\\
&= \half(e^{x\mu}+e^{-x\mu})\id + \half\mu^{-1}(e^{x\mu}-e^{-x\mu})A\\
&= \half e^{x\mu}(\id + \mu^{-1}A) + \half e^{-x\mu}(\id-\mu^{-1}A)
\spaceperiod
\end{align*}

The kernel of $\id-\mu^{-1}A$ is $L_+$.

The image of $\id-\mu^{-1}A$ is $L_-$.
Proof:
Let $v\in\bbC^2$ and let $w=(\id-\mu^{-1}A)v$
Then using $A^2 = (-\det A)\id = \mu^2A$,
\[
Aw = A(\id-\mu^{-1}A)v = Av-\mu^{-1}A^2v = Av - \mu v = -\mu(v-\mu^{-1}A)v = -\mu w
\spaceperiod
\]
Hence $[w]=L_-$.
(Note that $\exp(xA)L_{\pm} = L_{\pm}$.)

Hence for $L\in\CPone\setminus\{L_+\}$,
$\lim_{x\to-\infty}\exp(xA)L=(\id-\mu^{-1}A)L = L_-$.

\textit{Proof 2.}
Suppose first that $A=\diag(\mu,\,-\mu)$ is diagonal.
Then $L_+ = [\transpose{(1,0)}]$ and $L_- = [\transpose{(0,1)}]$.
Then $\exp(xA) = \diag(e^{x\mu},\,e^{-x\mu})$.
Let $(u,\,v)\in\bbC^2\setminus\{0\}$ with $v\ne 0$, so
$L = [\transpose{(u,\,v)}]\in\CPone\setminus\{L_+\}$.
Then $\exp(xA)L = [\transpose{(ue^{x\mu},\,ve^{-x\mu})}]$.
Hence $\lim_{x\to-\infty}\exp(xA)L = L_-$.

For nondiagonal $A$, let $CAC^{-1}$ be diagonal, etc.
\end{proof}
}

\begin{lemma}
\label{thm:bubbleton-asymptotics}
\theoremname{Bubbleton dressing}
Let
$A$ be a Delaunay residue and let $\mu$ be its eigenvalue
as in \autoref{def:delres}.
%
%
With $\RRR\in(0,\,1)$,
let $\lambda_0\in\calA_{\RRR,1}\cap\{\mu\in\bbR_+\}$,
and let $E$ be the eigenline of $\ol{\transpose{A(\lambda_0)}}$
corresponding to the eigenvalue $\ol{\mu(\lambda_0)}$.
For $k\in\{1,\,2\}$, let $L_k\in\CPone\setminus\{E\}$,
and let $g_k=\simplefactor{\lambda_0}{L_k}$ be normalized simple factors.
Let $\Phi_0=\exp(A\log z)$.
\Comment{$\Phi_0$ is on a universal cover.}
Then for every region $\calA\subseteq\calA_{\RRR}$ bounded away
from $\{\lambda_0,\,1/\ol{\lambda_0}\}$,
\CONVERGE{g_1\Phi_0}{g_2\Phi_0}
\begin{equation}\begin{split}
\label{eq:bublimit}
&\LIMIT{\left( \Uni{\RRR}{g_1\Phi_0} \right)^{-1}(g_1g_2^{-1})\Uni{\RRR}{g_2\Phi_0}}{\id}{\calA}{z\to 0}\spacecomma
\\
&\LIMIT{\Pos{\RRR}{g_2\Phi_0}{\left( \Pos{\RRR}{g_1\Phi_0} \right)^{-1}}}{\id}{\calD_{\RRR}}{z\to 0}
\spaceperiod
\end{split}\end{equation}

\end{lemma}

\begin{proof}
Let $F=\Uni{\RRR}{\Phi_0}$ and $B=\Pos{\RRR}{\Phi_0}$.
By the simple factor dressing formula~\eqref{eq:simpledress}
applied to $g_1$ and $g_2$, we have
\begin{equation*}
\Uni{\RRR}{g_k\Phi_0} = g_k F h_k^{-1}
\AND
\Pos{\RRR}{g_k\Phi_0} = h_k B\spacecomma\quad
k\in\{1,\,2\}\spacecomma
\end{equation*}
where
\begin{equation*}
h_k = \Pos{\RRR}{g_kF} =
 \simplefactor{\lambda_0}{GL_k}\spacecomma
\quad
G = \transpose{\ol{F(\lambda_0)}}\spaceperiod
\end{equation*}

By \autoref{thm:simple-factor-conj}, $h_1h_2^{-1}$ and $h_2h_1^{-1}$
are holomorphic on $\CPone\setminus\{\lambda_0,\,1/\ol{\lambda_0}\}$.
We compute in the coordinates $x+iy = \log z$.
We show that the lines $GL_1$ and $GL_2$ converge as $x\to -\infty$.
By~\eqref{eq:periodF},
\[
F(x + n\rho,\,y) = C^nF(x,\,y)\spacecomma\quad
C = \exp( (\rho-\dels)A )\spaceperiod
\]
Let $C_0=\transpose{\ol{C(\lambda_0)}}$.
By \autoref{thm:exp-limit},
for any $L\in\CPone\setminus\{E\}$,
the sequence $C_0^n L$ converges to $L_0$ as $n\to\infty$.
Let $\calF$ be the space of continuous functions $[0,\,\rho]\to\CPone$
with the $C^0$-norm.
Define the map $\calP:\calF\to\calF$ by
$f(x)\mapsto \transpose{\ol{F(x,\,y,\,\lambda_0))}}f(x)$.
Then $\calP$ is continuous,
\begin{Comment}
Because the multiplication map is continuous.
\end{Comment}
so $\transpose{\ol{F(x,\,y,\,\lambda_0))}}C_0^nL$ converges
to $\transpose{\ol{F(x,\,y,\,\lambda_0))}}L_0$ in $\calF$
as $n\to\infty$.
Hence
\[
\lim_{x\to\infty}
\dist\left(
  \transpose{\ol{F(x,\,y,\,\lambda_0)}}L_1,\,
  \transpose{\ol{F(x,\,y,\,\lambda_0)}}L_0\right) = 0
\spaceperiod
\]
Hence
$GL_1$ and $GL_2$ converge as $x\to -\infty$.
By \autoref{thm:simple-factor-limit},
for every region $\calA\subseteq\calA_{\RRR}$ bounded
away from $\{\lambda_0,\,1/\ol{\lambda_0}\}$,
\begin{equation}\begin{split}
\label{eq:bub120}
&\LIMIT{h_1h_2^{-1}}{\id}{\calA}{z\to 0}\spacecomma\\
&\LIMIT{h_2h_1^{-1}}{\id}{\calD_{\RRR}}{z\to 0}
\spaceperiod
\end{split}\end{equation}
%
Since
\begin{equation}\begin{split}
\label{eq:bub121}
&\left(\Uni{\RRR}{g_1\Phi_0}\right)^{-1}
(g_1 g_2^{-1})
\Uni{\RRR}{g_2\Phi_0}
 =(g_1 F h_1^{-1})^{-1}(g_1g_2^{-1})(g_2 F h_2^{-1}) = h_1 h_2^{-1}\spacecomma\\
&\Pos{\RRR}{g_2 \Phi_0} {\left(\Pos{\RRR}{g_1\Phi_0}\right)}^{-1}
 =
h_2 B \left({h_1 B}\right)^{-1} = h_2 h_1^{-1}
\spacecomma
\end{split}\end{equation}
the result~\eqref{eq:bublimit} follows
from~\eqref{eq:bub120} and~\eqref{eq:bub121}.
\end{proof}

\Comment{
Note that the bubbleton frame in \autoref{thm:bubbleton-asymptotics}
$\Uni{\RRR}{g_1\Phi_0}$
is non-singular at $\lambda_0$, while the Delaunay
frame $U\,\Uni{\RRR}{g_0\Phi_0}$ to which it is asymptotic
is singular at $\lambda_0$ since $U$ is.
}

\subsection{Delaunay dressing asymptotics}

\begin{xtext}
We show that for a certain special simple factor $g$,
the dressed Delaunay frame
$\dress{}{g}{\exp( (x+iy) A )}$ is a unitary Delaunay frame,
rather than a multibubbleton frame as is generically the case.
For each choice of $\lambda_0$, there are generically
two special Delaunay dressings, one for each eigenline of
$\transpose{\ol{A(\lambda_0)}}$.
\Comment{The geometry of the special Delaunay dressing is not well understood.
As a dressing approaches the special Delaunay dressing,
is the bubble going off to infinity? In which direction?
}
\end{xtext}

\begin{lemma}
\label{thm:special-del-dressing}
\theoremname{Special Delaunay dressing}
With $A$ a Delaunay residue,
let $g=\simplefactor{\lambda_0}{L}$ be the normalized simple
factor defined with
$\lambda_0\in\calD_1^\ast$,
and $L$ an eigenline of $\transpose{\ol{A(\lambda_0)}}$.
Then
$g A g^{-1}$ is a Delaunay residue.
\end{lemma}

\begin{proof}
By \autoref{thm:simple-factor-conj}(i),
$gAg^{-1}$ extends meromorphically to
$\bbC^\ast$
and this extension satisfies 
$(gAg^{-1})^\ast = gAg^{-1}$ away from its poles.

We show that $gAg^{-1}$ is holomorphic at $\lambda_0$.
As in \autoref{def:simple-factor}, let $U\in\matSU{2}{}$
and $f$ with $\ord_{\lambda_0}f=-1$ be
such that $g = U^{-1}(f^{1/2}\pi_L + f^{-1/2}\pi_{L^\perp})$.
Then
\begin{equation}\begin{split}
\label{eq:del-special-dressing1}
gAg^{-1} 
 &= U^{-1}(f^{1/2}\pi_L + f^{-1/2}\pi_{L^\perp})A
    (f^{-1/2}\pi_L + f^{1/2}\pi_{L^\perp})U\\
 &=
  U^{-1}\left( \pi_{L}A\pi_{L}
  + f\pi_{L}A\pi_{L^\perp}
  + f^{-1} \pi_{L^\perp}A\pi_{L}
  + \pi_{L^\perp}A\pi_{L^\perp} \right)U
\spaceperiod
\end{split}\end{equation}
%
At $\lambda_0$, the image of $\pi_{L^\perp}$ is $L^\perp$.
Since $L$ is an eigenline of $\transpose{\ol{A(\lambda_0)}}$,
it follows that $L^\perp$ is an eigenline of $A(\lambda_0)$.
Hence the image of $L^\perp$ under $A(\lambda_0)$ is $L^\perp$.
Since the image of $L^\perp$ under $\pi_L$ is $0$,
then $\ord_{\lambda_0}\pi_LA\pi_{L^\perp} \ge 1$.
Hence $f\pi_{L}A\pi_{L^\perp}$ is holomorphic at $\lambda_0$.
Since none of the other three terms in~\eqref{eq:del-special-dressing1}
has a pole at $\lambda_0$, then $gAg^{-1}$ is holomorphic at
$\lambda_0$.
By its hermitian symmetry, $gAg^{-1}$ is holomorphic at $1/\ol{\lambda_0}$,
and hence is holomorphic on $\bbC^\ast$.

We consider $gAg^{-1}$ at $\lambda=0$.
Since $g$ is holomorphic and upper-triangular at $\lambda=0$,
and the only pole of $A$ at $\lambda=0$ appears in the upper-right
entry, with order $-1$, the same is true of
$gAg^{-1}$.
Hence $gAg^{-1}$ is a Delaunay residue by \autoref{thm:delres}.
\end{proof}


\begin{xtext}
\autoref{thm:bubbleton-asymptotics} and \autoref{thm:bridge}
come together in the following theorem,
which shows that if the Iwasawa factors
 of a holomorphic frame $\Phi$ converge to those of a 
holomorphic Delaunay frame respectively,
then the same holds after dressing $\Phi$ by a finite product
of simple factors.
\end{xtext}

\begin{theorem}
\label{thm:simple-factor-asymptotics-theorem}
\theoremname{Simple factor dressing}
With $\RRR\in(0,\,1)$,
let $\Phi:\Sigma^\ast\to\LoopSL{\RRR}$ be an analytic map
on a punctured neighborhood $\Sigma^\ast\subset\bbC$ of $0\in\bbC$.
Let $A$ be a Delaunay residue.
Suppose that
\CONVERGE{\Phi}{\Phi_1 = \exp(A \log z)}
\begin{equation}
\label{eq:simple-factor-dressing-asym0}
\begin{split}
&\LIMIT{\left(\Uni{\RRR}{\Phi_1}\right)^{-1}\Uni{\RRR}{\Phi}}{\id}{\calA_r}{z\to 0}\spacecomma
\\
&\LIMIT{\Pos{\RRR}{\Phi}\left(\Pos{\RRR}{\Phi_1}\right)^{-1}}{\id}{\calD_{\RRR}}{z\to 0}
\spaceperiod
\end{split}
\end{equation}
%
%
Let $G\in\LooppSL{r}$ be a product of general simple factors (or $G=\id$)
with distinct singularities in
$\calA_{\RRR,\,1}\cap\{\mu\in\bbR_+\}$.
\Comment{Note that this set depends only on the eigenvalues of $A$,
not on $A$ itself, and hence doesn't change as $A$ gets conjugated
in the induction.}
Then there exist
\begin{itemize}
\item
a loop $U_2\in\LoopuMeroSL{\RRR}$
satisfying $\sing{U_2}=\sing{G}$,
\item
and a loop $V_2\in\LooppSL{\RRR}$ for which
$A_2:=V_2AV_2^{-1}$ is a Delaunay residue,
\end{itemize}
such that,
for every region $\calA\subseteq\calA_{\RRR}$ bounded away from
$\sing{U_2}$,
\CONVERGE{G\Phi}{\Phi_2:=\exp(A_2 \log z)}
\begin{align*}
&\LIMIT{\left(\Uni{\RRR}{\Phi_2}\right)^{-1}U_2\Uni{\RRR}{G\Phi}}{\id}{\calA}{z\to 0}\spacecomma
\\
&\LIMIT{\Pos{\RRR}{G\Phi}{V_2}^{-1}\left(\Pos{\RRR}{\Phi_2}\right)^{-1}}{\id}{\calD_{\RRR}}{z\to 0}
\spaceperiod
\end{align*}
\end{theorem}

\begin{proof}
The proof is by induction on the factors of $G$.
We first prove the following induction step,
which shows convergence after dressing by a single simple factor $g$.

\begin{nestedtheorem}
With $\RRR\in(0,\,1)$,
let $\Phi:\widetilde{\Sigma^\ast}\to\LoopSL{\RRR}$ be an analytic map
on a universal cover $\widetilde{\Sigma^\ast}\to\Sigma^\ast$
of a punctured neighborhood $\Sigma^\ast\subset\bbC$ of $0\in\bbC$.
Let $A$ be a Delaunay residue, and let $\mu$ be its eigenvalue
as in \autoref{def:delres}.
\begin{itemize}
\item
Let $U_1\in\LoopuMeroSL{\RRR}$.
\item
Let $V_1\in\LooppSL{\RRR}$
for which $A_1:=V_1AV_1^{-1}$ is a Delaunay residue.
\end{itemize}
Suppose
that for every region $\calA\subseteq\calA_{\RRR}$ bounded away
from $\sing{U_1}$,
\CONVERGE{\Phi}{\Phi_1 = \exp(A_1 \log z)}
\begin{equation}\begin{split}
\label{eq:bub1}
&\LIMIT{\left(\Uni{\RRR}{\Phi_1}\right)^{-1}U_1\Uni{\RRR}{\Phi}}{\id}{\calA}{z\to 0}\spacecomma
\\
&\LIMIT{\Pos{\RRR}{\Phi}V_1^{-1}\left(\Pos{\RRR}{\Phi_1}\right)^{-1}}{\id}{\calD_{\RRR}}{z\to 0}
\spaceperiod
\end{split}\end{equation}
%


Let $g$ be a general simple factor
with singularity $\lambda_0\in\calA_{\RRR,\,1}\cap\{\mu\in\bbR^+\}$.
%
Then there exist
\begin{itemize}
\item
a loop $U_2\in\LoopuMeroSL{\RRR}$
satisfying $\sing{U_2} = \sing{U_1}\cup\{\lambda_0,\,1/\ol{\lambda_0}\}$,
\item
and a loop $V_2\in\LooppSL{\RRR}$
for which $A_2:=V_2A{V_2}^{-1}$ is a Delaunay residue,
\end{itemize}
such that
for every region $\calA'\subseteq\calA_{\RRR}$ bounded away from
$\sing{U_2}$,
\CONVERGE{g\Phi}{\Phi_2:=\exp(A_2 \log z)}
\begin{equation}\begin{split}
\label{eq:bub2}
&\LIMIT{\left(\Uni{\RRR}{\Phi_2}\right)^{-1}U_2\Uni{\RRR}{g\Phi}}{\id}{\calA'}{z\to 0}\spacecomma
\\
&\LIMIT{\Pos{\RRR}{g\Phi}{V_2}^{-1}\left(\Pos{\RRR}{\Phi_2}\right)^{-1}}{\id}{\calD_{\RRR}}{z\to 0}
\spaceperiod
\end{split}\end{equation}
%

\end{nestedtheorem}

Proof of the induction step:

Let $\calA'\subset\calA_{r}$
be a region bounded away from $\sing{U_1}\cup\{\lambda_0,\,1/\ol{\lambda_0}\}$.

\Comment{
Summary:
\begin{enumerate}[label=\arabic*.]
\item
$\Phi\to\Phi_1$ (given).
\item
$g\Phi\to\galt\Phi_1$ (\autoref{thm:bridge}).
\item
$\galt\Phi_1\to \gdel\Phi_1$ (\autoref{thm:bubbleton-asymptotics}).
\item
$g\Phi\to \gdel\Phi_1$ (Steps 1 and 2).
\end{enumerate}
}

\STARTSTEPS
\STEP
%
Write $g=W\simplefactor{\lambda_0}{L}$
and
let $\galt = W\simplefactor{\lambda_0}{{\transpose{\ol{U(\lambda_0)}}}^{-1}L}$.
By \autoref{thm:bridge} applied to~\eqref{eq:bub1},
\CONVERGE{g\Phi}{\galt\Phi_1}
\begin{equation}\begin{split}
\label{eq:bub3}
&\LIMIT{\left(\Uni{\RRR}{\galt\Phi_1}\right)^{-1} (\galt U_1 g^{-1}) \Uni{\RRR}{g\Phi}}{\id}{\calA'}{z\to 0}\spacecomma
\\
&\LIMIT{\Pos{\RRR}{g\Phi}V_1^{-1}\left(\Pos{\RRR}{\galt\Phi_1}\right)^{-1}}{\id}{\calD_{\RRR}}{z\to 0}
\spaceperiod
\end{split}\end{equation}

\STEP
Let $\mu$ be the eigenvalue function of $A$
as in \autoref{def:delres}.
Let $(\ol{\mu(\lambda_0)},\,E_+)$ and $(-\ol{\mu(\lambda_0)},\,E_-)$
be the eigenvalue-eigenline pairs of $\transpose{\ol{A_1(\lambda_0)}}$.
Let $\gdel = \simplefactor{\lambda_0}{E_-}$
be the normalized simple factor as in \autoref{thm:special-del-dressing}.
By that lemma,
$A_2 := \gdel A_1\gdel^{-1}$ is a Delaunay residue.
Since $\mu(\lambda_0)\ne 0$, then $E_-\ne E_+$.
Hence by \autoref{thm:bubbleton-asymptotics} applied to $\gdel$ and $\galt$,
the loops
\CONVERGE{\galt\Phi_1}{\gdel\Phi_1}
\begin{equation}\begin{split}
\label{eq:bub4a}
&\LIMIT{\left(\Uni{\RRR}{\gdel\Phi_1}\right)^{-1}\gdel\galt^{-1}\Uni{\RRR}{\galt\Phi_1}}{\id}{\calA'}{z\to 0}\spacecomma
\\
&\LIMIT{\Pos{\RRR}{\galt\Phi_1}\left(\Pos{\RRR}{\gdel\Phi_1}\right)^{-1}}{\id}{\calD_{\RRR}}{z\to 0}
\spaceperiod
\end{split}\end{equation}
%
Let $\Phi_2 = \exp(A_2\log z)$.
Then $\gdel\Phi_1 = \Phi_2 \gdel$,
so by~\eqref{eq:bub4a},
\CONVERGE{\galt\Phi_1}{\Phi_2}
\begin{equation}\begin{split}
\label{eq:bub4}
&\LIMIT{\left(\Uni{\RRR}{\Phi_2}\right)^{-1}\gdel\galt^{-1}\Uni{\RRR}{\galt\Phi_1}}{\id}{\calA'}{z\to 0}\spacecomma
\\
&\LIMIT{\Pos{\RRR}{\galt\Phi_1}\gdel^{-1}\left(\Pos{\RRR}{\Phi_2}\right)^{-1}}{\id}{\calD_{\RRR}}{z\to 0}
\spaceperiod
\end{split}\end{equation}

\STEP
By \autoref{thm:simple-factor-conj}(ii),
$U_2:=\gdel U_1 g^{-1}$ is an element of $\LoopuMeroSL{r}$
and $\sing{U_2} = \sing{U_1}\cup\{\lambda_0,\,1/\ol{\lambda_0}\}$.
The loop $V_2:=\gdel V_1$ is an element of $\LooppSL{\RRR}$
and $V_2 A V_2^{-2} = \gdel A_1 \gdel^{-1}$ is a Delaunay residue.
By~\eqref{eq:bub3} and \eqref{eq:bub4},
\CONVERGE{g\Phi}{\Phi_2} given by~\eqref{eq:bub2}.
%
%
This proves the induction step.

To prove the theorem,
first note that if $G=\id$, the theorem is trivially true
with $U_2=V_2=\id$.
Otherwise, write $G=g_n\cdots g_1$ as a product of simple factors,
and starting with $U_1=V_1=\id$,
apply the induction step to $g_1,\dots,g_n$ in turn.
\end{proof}

\begin{xtext}
By the methods of \autoref{Sec:immersion},
\autoref{thm:simple-factor-asymptotics-theorem}
implies the following
(see also~\cite{Kobayashi_2006}):
\end{xtext}

\begin{theorem}
\label{thm:shimpei}
\theoremname{Delaunay convergence}
If a CMC surface $f$
converges to a half Delaunay surface $f_0$, then
the surface produced from $f$ by dressing by finitely many
simple factors converges to a rigid motion of $f_0$.
\end{theorem}

\Depend{
\begin{proof}
The theorem depends on
\autoref{thm:simple-factor-asymptotics-theorem},
\autoref{thm:conformal-immersion-asymptotics}, and
\autoref{thm:embedded}.
\end{proof}
}

\typeout{=================frame-asymptotics2}
\subsection{Dressed Delaunay asymptotics}
\label{sec:delframe-asymptotics2}

\begin{xtext}
We show that the $r$-Iwasawa factors of a holomorphic frame
obtained from a perturbed Delaunay potential
converges to those of a Delaunay frame.
This theorem is a generalization of \autoref{thm:frame1},
removing the restriction that the holomorphic frame
extend to $\calA_{r,1}$.
\end{xtext}

\Comment{
The setup for \autoref{thm:frame2} is as follows:
\begin{itemize}
\item
Let $A$ be a Delaunay residue as in \autoref{def:delres}.
\item
Let $\mu$ be an eigenvalue of $A$ as in \autoref{def:delres}.
\item
Suppose
\begin{equation*}
\min_{\lambda\in\calC_r} \Real\mu(\lambda) \le 1/2
\spaceperiod
\end{equation*}
\item
Let $r\in(0,\,1)$.
\item
Let $\Sigma\subset\bbC$ be a neighborhood of $0\in\bbC$.
\item
Let $\xi$ be an $r$-potential (\autoref{def:setup}) on $\Sigma^\ast$
satisfying
\begin{equation*}
\xi=Az^{-1}dz + \Order(z^{n})dz\spaceperiod
\end{equation*}
\item
Let $\widetilde{\Sigma^\ast}\to\Sigma^\ast$ be a universal cover of
$\Sigma^\ast$.
\item
Let $\Phi:\widetilde{\Sigma^\ast}\to\LoopSL{r}$
satisfy $d\Phi=\Phi\xi$.
\item
Assume that the monodromy $M$ of $\Phi$ at $z=0$
satisfies $M\in\LoopuSL{r}$.
\end{itemize}
}

\begin{theorem}
\label{thm:frame2}
\theoremname{Frame asymptotics II}
Let $r\in(0,\,1)$ and assume $\calC_r\cap\delSA = \emptyset$.
Let $\xi$ be a perturbed Delaunay $r$-potential
\begin{equation}
\label{eq:delasym-xi-2}
\xi=Az^{-1}dz + \Order(z^{n})dz\spaceperiod
\end{equation}
Suppose
\begin{equation}
\label{eq:delasym-max-2}
\max_{\lambda\in\calC_r} \Real\mu(\lambda) < (n+1)/2
\spaceperiod
\end{equation}
Let $\Phi:\widetilde{\Sigma^\ast}\to\LoopSL{r}$
satisfy $d\Phi=\Phi\xi$
on the universal cover $\widetilde{\Sigma^\ast}\to\Sigma^\ast$
of $\Sigma^\ast$,
and assume that the monodromy $M$ of $\Phi$ around $z=0$
satisfies $M\in\LoopuSL{r}$.

Then there exists
a loop $U\in\LoopuMeroSL{\RRR}$
with $\sing{U}\subset\calA_{r,1}\cap\delSA$
and only simple poles,
and a loop $V\in\LooppSL{\RRR}$
for which $A_1:=VAV^{-1}$
is a Delaunay residue,
such that
for every region $\calA\subseteq\calA_{\RRR}$ bounded away
from $\sing{U}$,
\CONVERGE{\Phi}{\Phi_1:=\exp(A_1\log z)}
\begin{equation}\begin{split}
\label{eq:delasym2-FB}
&\LIMIT{\left(\Uni{\RRR}{\Phi_1}\right)^{-1}U\Uni{\RRR}{\Phi}}{\id}{\calA}{z\to 0}\spacecomma
\\
&\LIMIT{\Pos{\RRR}{\Phi}V^{-1}\left(\Pos{\RRR}{\Phi_1}\right)^{-1}}{\id}{\calD_{\RRR}}{z\to 0}\spaceperiod
\end{split}\end{equation}
\end{theorem}

\begin{proof}
\Comment{
Summary:
\begin{enumerate}[label=\arabic*.]
\item
$\Phi\to C\Phi_0$ (Generalization of \autoref{thm:frame1}.)
\item
$\Phi\to G\Phi_0$ (Replacing $C$ with $G$)
\item
$G\Phi_0\to\Phi_1$ (\autoref{thm:simple-factor-asymptotics-theorem}).
\item
$\Phi\to\Phi_1$ (Steps 2 and 3).
\end{enumerate}
}
\STARTSTEPS
\STEP
Let $\Phi_0=\exp(A\log z)$ and let
$\Phi = C\Phi_0 P$ be the $z^AP$-decomposition of $\Phi$.
By \autoref{thm:delgrowth2},
there exist $0<s_1<r<s_2<1$ and
a continuous function $c:\calA_{s_1,s_2}\to\bbR_+$ such that
$\Pos{r}{C\Phi_0}$ extends continuously to $\calA_{s_1,s_2}$ and
\[
\supnorm{\Pos{r}{C\Phi_0}}{\calA_{s_1,s_2}} \le c \abs{z}^{-\Real\mu}
\spaceperiod
\]

\STEP
By \autoref{thm:frame1},
\CONVERGE{\Phi}{C\Phi_0}
\begin{equation}\begin{split}
\label{eq:delasym1}
&\LIMIT{ \left(\Uni{\RRR}{C\Phi_0}\right)^{-1}\Uni{\RRR}{\Phi}}{\id}{\calA_{\RRR}}{z\to 0}
\spacecomma\\
&\LIMIT{ \Pos{\RRR}{\Phi}\left(\Pos{\RRR}{C\Phi_0}\right)^{-1} }{\id}{\calD_{\RRR}}{z\to 0}
\spaceperiod
\end{split}\end{equation}

\STEP
Let $C=C_u\cdot C_+$ be the $\RRR$-Iwasawa factorization of $C$.
Using that $C\exp(2\pi i A)C^{-1}\in\LoopuSL{\RRR}$,
by \autoref{thm:bub},
there exist loops $G,\,V_1\in\LooppSL{\RRR}$ such that
$\CPLUS = GV_1$,
$G$ is a product of normalized simple factors (or $G=\id$),
$A_0:=V_1AV_1^{-1}$ is a Delaunay residue,
and the singularities of the simple factors are distinct
and in $\delSA\cap(\calA_{\RRR}\setminus\bbS^1)$.
Then, with $\Phi_0 = \exp(A_0\log z)$,
\begin{equation}
\label{eq:delasym2}
C \Phi_0
 = C_u GV_1 \exp(A \log z)
 = C_u G \exp(A_0 \log z) V_1
 = C_u G\Phi_0V_1 \spaceperiod
\end{equation}
By~\eqref{eq:delasym1} and~\eqref{eq:delasym2},
\CONVERGE{\Phi}{G\Phi_0}
\begin{equation}\begin{split}
\label{eq:main10}
&\LIMIT{\left(\Uni{\RRR}{G\Phi_0}\right)^{-1}C_u\Uni{\RRR}{\Phi}}{\id}{\calA_{\RRR}}{z\to 0}\spacecomma
\\
&\LIMIT{\Pos{\RRR}{\Phi}V_1^{-1}\left(\Pos{\RRR}{G\Phi_0}\right)^{-1}}{\id}{\calD_{\RRR}}{z\to 0}\spaceperiod
\end{split}\end{equation}

\STEP
%
The singularities of $G$ in $\calD_1$ are distinct and lie in
$\calA_{r,1}\cap\{\mu\in\bbR_+\}$.
By \autoref{thm:simple-factor-asymptotics-theorem},
where in~\eqref{eq:simple-factor-dressing-asym0},
the $\Phi$ and $\Phi_0$ are both replaced by $\Phi_0$,
so that~\eqref{eq:simple-factor-dressing-asym0}
holds vacuously,
there exists a loop $V_2\in\LooppSL{\RRR}$
for which $A_1 := V_2AV_2^{-1}$
is a Delaunay residue,
and a loop $U_2\in\LoopuMeroSL{\RRR}$
such that
for every region $\calA\subset\calA_{\RRR}$ bounded away from $\sing{U_2}$,
\CONVERGE{G\Phi_0}{\Phi_1:=\exp(A_1\log z)}
\begin{equation}\begin{split}
\label{eq:main111}
&\LIMIT{\left(\Uni{\RRR}{\Phi_1}\right)^{-1}U_2\Uni{\RRR}{G\Phi_0}}{\id}{\calA}{z\to 0}\spacecomma
\\
&\LIMIT{\Pos{\RRR}{G\Phi_0}V_2^{-1}\left(\Pos{\RRR}{\Phi_1}\right)^{-1}}{\id}{\calD_{\RRR}}{z\to 0}\spaceperiod
\end{split}\end{equation}

\STEP
The result~\eqref{eq:delasym2-FB}
follows from~\eqref{eq:main10} and \eqref{eq:main111},
with $U=U_2C_u$ and $V=V_2V_1$.
\end{proof}

\typeout{=================main2}

\begin{xtext}
We conclude with the main theorem of the paper,
showing that a CMC end obtained from a perturbed Delaunay
potential is asymptotic to a half-Delaunay surface.
This theorem is a generalization of \autoref{thm:main1},
removing the conditions $\calA_{r,1}\cap\delSA=\emptyset$
and $\Phi\in\LoopUpSL{r}$.
Hence it applies to surfaces obtained by $r$-dressing,
such as the $n$-noids with bubbles constructed
in~\cite{Kilian_Schmitt_Sterling_2004}.
The proof is similar
to that of \autoref{thm:main1},
taking the appearance of the loops $U$ and $V$ into consideration.
In the case of embedded ends,
we obtain exponential convergence
(see \autoref{thm:exponential-convergence}).
\end{xtext}

\begin{theorem}
\label{thm:main2}
\theoremname{Main theorem II}
\Comment{Need to make sure $\lambda_0\not\in\sing{U}$, where $\lambda_0$
is the Sym point.}
Let $A$ be a Delaunay residue as in~\eqref{eq:delres}
satisfying condition~\eqref{eq:delasym-max}.
On the punctured unit disk $\Sigma^\ast=\{z\in\bbC\suchthat 0<\abs{z}<1\}$,
let
\begin{equation*}
\xi = A\frac{dz}{z} + \Order(z^0)dz
\end{equation*}
be a perturbed Delaunay $r$-potential, $r\in(0,\,1]$.
Let $f_0$ and $f$ be the immersions of $\Sigma^\ast$
induced by the generalized Weierstrass representation at $\lambda=1$
by $A dz/z$ and $\xi$ respectively,
so $f_0$ is a Delaunay immersion.
Assume that $f$ is obtained from a holomorphic $r$-frame
whose monodromy at $z=0$ is in $\LoopuSL{r}$.

With $\Sigma_{\rm Cyl}$ as in ~\eqref{eq:cylinder},
let $\phi:\Sigma_{\rm Cyl}\to\Sigma^\ast$ be the map $\phi(x,\,y)=e^{x+iy}$.
Then some rigid motion of $f\circ\phi$ converges to $f_0\circ\phi$
in the $C^\infty$-topology of $\Sigma_{\rm{Cyl}}$ as $x\to -\infty$.
Furthermore, if $f_0$ is embedded, then $f$ is properly embedded.
\end{theorem}

\Depend{
\begin{proof}
The theorem depends on
\autoref{thm:delasym-summary-1},
\autoref{thm:embedded},
\autoref{thm:frame2},
\autoref{thm:conformal-immersion-asymptotics},
\autoref{thm:main1}.
\end{proof}
}

\bibliographystyle{amsplain}
\bibliography{asymptotics}

\end{document}